\documentclass[10pt]{amsart}
\usepackage[utf8]{inputenc}
\usepackage{amsmath, amsthm, amsfonts, amssymb}
\usepackage{tikz-cd}
\usepackage{dsfont} % indicator function
\usepackage{zref-clever}
\usepackage{hyperref}
\usepackage{enumitem}

\hypersetup{
  colorlinks   = true, %Colors links instead of ugly boxes
  urlcolor     = blue, %Color for external hyperlinks
  linkcolor    = blue, %Color of internal links
  citecolor    = blue  %Color of citations
}

\let\oldtocsection=\tocsection
\let\oldtocsubsection=\tocsubsection
\renewcommand{\tocsection}[2]{\hspace{0em}\oldtocsection{#1}{#2}}
\renewcommand{\tocsubsection}[2]{\hspace{1em}\oldtocsubsection{#1}{#2}}

\AddToHook{env/theorem/begin}{\zcsetup{countertype={equation=theorem}}}
\AddToHook{env/lemma/begin}{\zcsetup{countertype={equation=lemma}}}
\AddToHook{env/claim/begin}{\zcsetup{countertype={equation=claim}}}
\AddToHook{env/corollary/begin}{\zcsetup{countertype={equation=corollary}}}
\AddToHook{env/proposition/begin}{\zcsetup{countertype={equation=proposition}}}
\AddToHook{env/definition/begin}{\zcsetup{countertype={equation=definition}}}
\AddToHook{env/assumption/begin}{\zcsetup{countertype={equation=assumption}}}
\AddToHook{env/question/begin}{\zcsetup{countertype={equation=question}}}
\AddToHook{env/conjecture/begin}{\zcsetup{countertype={equation=conjecture}}}
\AddToHook{env/examples/begin}{\zcsetup{countertype={equation=examples}}}
\AddToHook{env/remark/begin}{\zcsetup{countertype={equation=remark}}}
\AddToHook{env/remarks/begin}{\zcsetup{countertype={equation=remarks}}}
\AddToHook{env/example/begin}{\zcsetup{countertype={equation=example}}}

\numberwithin{equation}{section}

\theoremstyle{definition}

\newtheorem{theoremx}{Theorem}

\zcsetup{countertype={theoremx=theorem}}
\newtheorem{theorem}[equation]{Theorem}
\newtheorem{lemma}[equation]{Lemma}

\newtheorem{corollary}[equation]{Corollary}
\newtheorem{proposition}[equation]{Proposition}
\newtheorem{definition}[equation]{Definition}

\newtheorem{question}[equation]{Question}
\newtheorem{conjecture}[equation]{Conjecture}

\newtheorem{remark}[equation]{Remark}

\renewcommand{\div}{{\mathrm{div}}}

\newcommand{\inj}{{\mathrm{inj}}}
\newcommand{\loc}{{\mathrm{loc}}}
\newcommand{\Ric}{{\mathrm{Ric}}}
\newcommand{\Rm}{{\mathrm{Rm}}}
\newcommand{\Spec}{{\mathrm{Spec}}}
\newcommand{\supp}{{\mathrm{supp}}}
\newcommand{\Hess}{{\mathrm{Hess}}}

\newcommand{\Vol}{{\mathrm{Vol}}}

\newcommand{\ps}{{\mathrm{ps}}}

\newcommand{\norm}[1]{{\lVert #1 \rVert}}

\newcommand{\RR}{{\mathbb{R}}}

\newcommand{\cA}{{\mathcal{A}}}

\newcommand{\cC}{{\mathcal{C}}}

\newcommand{\cE}{{\mathcal{E}}}

\newcommand{\cL}{{\mathcal{L}}}

\newcommand{\cP}{{\mathcal{P}}}

\newcommand{\cR}{{\mathcal{R}}}

\newcommand{\cU}{{\mathcal{U}}}

\newcommand{\cW}{{\mathcal{W}}}

\title{Blowdowns along Ricci flow in dimension four}
\author{Dain Kim \and Tristan Ozuch}
\address{MIT, Dept. of Math., 77 Massachusetts Avenue, Cambridge, MA 02139-4307}
\email{dain0327@mit.edu}
\address{MIT, Dept. of Math., 77 Massachusetts Avenue, Cambridge, MA 02139-4307}
\email{ozuch@mit.edu}

\begin{document}
\begin{abstract}
We show that the blowdown shrinking Ricci soliton $(M,\overline g, f)$ of Feldman, Ilmanen and Knopf is strictly dynamically stable in the following strong local sense:
\begin{enumerate}
    \item The Ricci flow from every metric on a closed four-manifold containing a region sufficiently close to a rescaling of a sufficiently large compact subset of $(M,\overline g)$ in either orientation and with controlled geometry elsewhere develops a singularity with unique tangent flow $(M,\overline g, f)$.
    \item The singular time slice has an isolated conical singularity, and, by work of Gianniotis--Schulze, the flow can be continued through the singularity by blowing down a sphere of self-intersection $\pm1$. 
    \item On a closed manifold $N^4$, the set of initial metrics whose Ricci flows have their first singularity modeled on $(M,\overline g, f)$ is $C^{14}$-open. It is nonempty if and only if $N^4$ contains a sphere of self-intersection $\pm1$. 
    \item The soliton $(M,\overline g, f)$  is strongly rigid among shrinking Ricci solitons. 
\end{enumerate}
  As applications, we construct $C^{14}$-open sets of initial data whose Ricci flows perform any prescribed ordered sequence of blowdowns. Because these initial data sets are open, the following exotic phenomena must be accounted for in any theory of generic Ricci flow: 
  \begin{enumerate}
      \item Every closed simply connected smooth four-manifold, including spin and exotic examples, can be reached by Ricci flow from a standard $k\mathbb{CP}^2\#\ell\overline{\mathbb{CP}}{}^2$, through finitely many blowdown singularities.
      \item On the fixed standard $2\mathbb{CP}^{2}\#9\overline{\mathbb{CP}}{}^{2}$, there are infinitely many  open sets of initial metrics whose Ricci flows, after a single blowdown, reach pairwise nondiffeomorphic manifolds.
  \end{enumerate}
The proof uses the spectral analysis of Naff and the second author and deals with the notoriously difficult gauge-fixing problem along Ricci flow through a Ważewski argument inspired by Stolarski. This is combined with and simplified by ideas borrowed from Bamler and Lai's PDE-ODI principle. 

All of the results, except the continuation of the flow (which requires nonnegativity of the curvature of the asymptotic cone) in this article apply verbatim to linearly stable nonflat asymptotically conical shrinking solitons in all dimensions.
\end{abstract}

\maketitle

\section{Introduction}

Ricci flow, introduced by Hamilton \cite{Ham82}, has had profound implications in three-dimensional geometry and topology. The success of this program relied on a thorough understanding of singularity formation, built on the Hamilton-Ivey pinching estimate \cite{Ham82,Ive93}, Hamilton's Harnack inequality \cite{Ham93} and Perelman's noncollapsing estimate \cite{Per02}. With these tools, Perelman classified three-dimensional shrinking solitons and proved his \emph{canonical neighborhood theorem}: regions of high curvature look like regions of ancient $\kappa$-solutions, such as round spheres, cylinders and their quotients \cite{Per02,Per03}. This allowed him to continue the flow through singularities by surgery and to complete Hamilton's program. Three dimensional ancient $\kappa$-solutions have since been classified \cite{Bre20,BDS21,BK21}, and the flow through singularities has been made canonical \cite{KL17,BK22}. In arbitrary dimension, Type I singularities admit gradient shrinking solitons as blowup limits \cite{Nab10,EMT11}.
More generally, Bamler's entropy and heat kernel estimates \cite{Bam21a}
provide the analytic foundation of a compactness and partial regularity
theory for Ricci flow \cite{Bam23,Bam21c}.

In dimension four, a \emph{classification of shrinking solitons} seems out of reach, since it contains the classification of Einstein four-manifolds with positive scalar curvature. This classification is far from complete: Einstein metrics with positive scalar curvature were recently constructed on infinitely many topologies \cite{Ozu26b,LS26, BL26,CS26}. However, infinitely many of these metrics are unstable by \cite{Ozu26b}, and all of these metrics are expected to be. It has been proposed to focus instead on \textit{stable} solitons, which should be much less numerous and are expected to be the only singularity models of a \textit{generic} Ricci flow, that is, for initial data in an open dense set. They should correspond to the topological surgeries performed by Ricci flow starting at an open dense set of metrics at its singularities. The second author and Naff showed that the blowdown soliton constructed by Feldman, Ilmanen and Knopf \cite{FIK03} is \emph{linearly stable} \cite{NO25a}. This suggested that the blowdown of spheres of self-intersection $\pm1$ is such a generic surgery.

In this article, we prove that this blowdown soliton $(M,\overline{g},f)$ is indeed \emph{dynamically stable} in the strongest sense. Every Ricci flow that has a region close enough to a rescaling of a large compact of $(M,\overline{g})$  and bounded geometry elsewhere develops a singularity uniquely modeled on $(M,\overline{g},f)$. In other words, the blowdown soliton belongs to the list of \emph{canonical neighborhoods of generic Ricci flow} in dimension four, and the blowdown is a \emph{generic topological surgery} of Ricci flow through the continuation of \cite{GS18}, which we show is applicable. This also lets us easily construct metrics from which Ricci flow has \textit{prescribed singularities}. Some topological applications include the construction of Ricci flows through singularities from standard $k\mathbb{CP}^2\#l\overline{\mathbb{CP}}{}^2$ to arbitrary exotic $4$-manifolds and open sets of metrics on $2\mathbb{CP}^2\#9\overline{\mathbb{CP}}{}^2$ whose Ricci flows reach infinitely many exotic $\mathbb{CP}^2\#9\overline{\mathbb{CP}}{}^2$ after their first singularity.

This is in stark contrast with the mean curvature flow of hypersurfaces, where spheres and cylinders are the only stable singularity models \cite{CM12}. The stability of cylinders is also much weaker. Indeed, a Ricci flow which comes close to a cylinder on a large compact set need not develop a neckpinch there, even in dimension three: it can instead develop a degenerate neckpinch modeled on the Bryant soliton \cite{AIK15}, and the same happens for the mean curvature flow of surfaces \cite{AV97}. By contrast, a Ricci flow which comes close enough to a large compact of the blowdown soliton always develops a blowdown singularity.

The main difficulty is to find the right gauge and scale to compare a Ricci flow on a closed manifold with the noncompact soliton $(M,\overline{g},f)$ when the two are only close on a compact region at the initial time. We develop a new way to gauge-fix based on a Wa\.zewski box argument simplified by the PDE-ODI approach. While Stolarski perturbs the metric in \cite{Sto26}, our parameters change only the diffeomorphism and scale. This  is possible because the nonnegative eigenmodes of the soliton are all generated by scaling and diffeomorphisms \cite{NO25a}. For cylinders, the linearized operator has neutral modes which are not generated by symmetries, and their analysis uses a \L{}ojasiewicz inequality \cite{CM15} or the PDE-ODI principle \cite{BL25b}. \zcref[S]{subsec:strategy} describes the argument and its relation to \cite{Sto26,BL25b,CM25}.

\subsection{Background and motivation}

\subsubsection{The blowdown soliton}

The blowdown soliton $(M,\overline g, f)$ of Feldman, Ilmanen and Knopf \cite{FIK03}, which we also call FIK, lives on the total space $M$ of $\mathcal O(-1)\to\mathbb{CP}^1$, that is, on the blow-up of $\mathbb C^2$ at a point, or $\mathbb R^4\#\overline{\mathbb{CP}}{}^2$. It solves $\mathrm{Ric}(\overline{g})+\nabla^{\overline{g}, 2}  f=\frac12\overline{g}$, it is K\"ahler and $U(2)$-invariant, and it is asymptotic to a cone $\mathcal C$ over a Berger sphere, with nonnegative curvature operator. Its zero section $\overline\Sigma$ is a sphere of self-intersection $-1$, which shrinks to a point at the singular time.

Song and Weinkove established the contraction of exceptional curves
by K\"ahler-Ricci flow \cite{SW11,SW13}.
M\'aximo showed that the blowdown soliton arises as a singularity model
of K\"ahler-Ricci flows on compact surfaces \cite{Max14}. In fact, finite time singularities of the K\"ahler-Ricci flow on compact K\"ahler surfaces with volume bounded below are of Type I \cite{CHM25}, hence modeled on the blowdown soliton \cite{CCD26}. Up to scaling and isometry, it is the unique complete K\"ahler
shrinking soliton on its underlying smooth manifold \cite{CDS24,LW23}. Its central density, about $0.672$, is the lowest among the known stable shrinking solitons \cite{CHI04,NO25a}. Isenberg, Knopf and \v{S}e\v{s}um conjectured the dynamical stability
of the blowdown soliton and gave partial evidence within a class
of $U(2)$-invariant metrics
\cite[Section 1.3, Main Conjecture]{IKS19}.

\subsubsection{Generic Ricci flow and canonical neighborhoods}

In dimension four, a generic version of Perelman's canonical neighborhood theorem has been proposed \cite{NO25a}.

\begin{conjecture}\label{conj:generic}
For an open and dense set of initial metrics on a closed four-manifold, every region of high curvature along the Ricci flow is, after rescaling, close to a region of a stable shrinking or steady Ricci soliton. Moreover, there are only a few such solitons.
\end{conjecture}

This conjecture is motivated by the apparent scarcity of stable solitons. Apart from $S^4$ and $\mathbb{CP}^2$, all the Einstein metrics with positive scalar curvature known before \cite{Ozu26b,LS26} are unstable \cite{BO23}, and $\mathbb{CP}^2$ is dynamically unstable \cite{Kro15,KS19}, its perturbations developing blowdown singularities \cite{GIKW24}, at least numerically. Compact K\"ahler shrinkers \cite{HM11} are unstable, as is the shrinker constructed by Bamler, Cifarelli, Conlon and Deruelle \cite[Theorem A]{BCCD22}, whose instability was proved in \cite{NO25b}. Even among Ricci-flat metrics, all known non hyperk\"ahler examples are unstable \cite{BO23, KO25, Ozu26a}.

Up to quotients, the nonflat shrinking solitons that are expected to be stable correspond to a generic topological surgery:
\begin{itemize}
\item $S^4$, which corresponds to a component vanishing in a round point,
\item $S^3\times\mathbb R$, the neckpinch, which undoes a connected sum, $S^3\times B^1\to B^4\times S^0$,
\item $S^2\times\mathbb R^2$, the bubble sheet, which is believed to correspond to a handle surgery, $S^2\times B^2\to B^3\times S^1$,
\item the blowdown soliton on $\mathcal O(-1)$, or $\mathcal O(1)$ in the opposite orientation, which corresponds to the blowdown, $B^2(\mathcal O(-1))\to B^4$, where $B^2(\mathcal O(-1))$ is a tubular neighborhood of the zero section.
\end{itemize}

The results of the present article show that the blowdown soliton indeed belongs to this list: Corollary \ref{cor:open} shows that blowdown singularities form from an open set of initial data: no generic version of the canonical neighborhood theorem can avoid them. \zcref[S]{thm:B} describes the corresponding canonical neighborhood. Near the singular point, the flow is modeled on the blowdown soliton at every scale before the singular time. The singular time slice has a conical singularity over a Berger sphere, and the flow continues through the singularity following the geometric surgery proposed in \cite{FIK03} and constructed in \cite{GS18}. 

Whether the above \emph{flow through a singularity is unique} is a central question. Uniqueness holds in dimension three \cite{BK22}, fails from dimension five \cite{AK22}, and remains open in dimension four. For the blowdown soliton, it was asked by Feldman, Ilmanen and Knopf \cite[Question 4]{FIK03}, and partial answers are known in the K\"ahler setting \cite{Che25,CHL26,CHL26b}. We conjecture that the continuation of \zcref[S]{thm:B} is unique in a natural class (Conjecture \ref{conj:canonical}). 

\subsubsection{Generic surgeries and four-manifolds}

After connected sums and handle surgeries, the next simplest topological operation on $4$-manifolds is the complex blow-up, the connected sum with $\mathbb{CP}^2$ or $\overline{\mathbb{CP}}{}^2$. By Freedman \cite{Fre82} and Donaldson \cite{Don83}, every closed, simply connected, smooth, non-spin four-manifold is homeomorphic to a connected sum of copies of $\mathbb{CP}^2$ and $\overline{\mathbb{CP}}{}^2$. \zcref[S]{thm:C} below shows that Ricci flow can perform any prescribed sequence of blowdowns. This shows in Corollary \ref{cor:exotic-stab} that Ricci flow can flow through singularities from a standard connected sum of $\mathbb{CP}^2$ and $\overline{\mathbb{CP}}{}^2$ to arbitrary closed, simply connected $4$-manifolds, even exotic and spin ones. Similarly, infinitely many open sets of initial metrics on $2\mathbb{CP}^2\#9\overline{\mathbb{CP}}{}^2$ have their Ricci flow reach infinitely many exotic $\mathbb{CP}^2\#9\overline{\mathbb{CP}}{}^2$ after the blowdown of a $(+1)$-sphere, see Corollary \ref{cor:exotic}.
\subsubsection{Notions of stability}

Ricci flow is the gradient flow of Perelman's entropy, and a shrinking soliton is a fixed point of the rescaled flow. The soliton is \emph{variationally stable} if it locally maximizes Perelman's entropy, and \emph{dynamically stable} if the rescaled flows starting near it flow back to it, up to the symmetries given by scaling and diffeomorphisms. The corresponding notion of linear stability is defined by the nonpositivity of the second variation of Perelman's entropy
\cite{CHI04,CZ12}.

\zcref[S]{thm:A} and \zcref[S]{thm:B}(1) establish the dynamical stability of the blowdown soliton, and Corollary \ref{cor:entropy} its variational stability, in a very strong form. Only a compact region needs to look like the soliton, at some scale. The rest of the manifold does not matter, as long as it has controlled geometry at a slightly larger scale, the flow develops a singularity modeled on the soliton, the rescaled flow converges back to it, and the entropy at the singular scale is below that of the soliton.

\subsubsection{Comparison with mean curvature flow}

For the mean curvature flow of hypersurfaces, spheres and cylinders are the only linearly stable shrinkers \cite{CM12}. Huisken's conjecture that a generic flow only has spherical and cylindrical singularities \cite{Ilm03} has since been established in many cases, notably for closed surfaces in $\mathbb R^3$ and for hypersurfaces with low entropy \cite{CCMS24a,CCMS24b,CMS25,BK23}. Cylinders are moreover rigid: a shrinker close to a cylinder on a large compact set is a cylinder \cite{CIM15}. Cylindrical tangent flows are also unique \cite{CM15}, and so are asymptotically conical ones \cite{CS21}, which is the analogue of \zcref[S]{thm:B}(1). The analogous rigidity of cylinders among Ricci shrinkers was proved by Colding and Minicozzi \cite{CM25}, Li and Wang \cite{LW24}, and Li and Zhang \cite{LZ25}. The latter also show that tangent flows are unique when one of them is a generalized cylinder. Near neck singularities, the mean convex neighborhood conjecture of Ilmanen \cite{Ilm03} describes the flow by a canonical model. It was proved by Choi, Haslhofer and Hershkovits for surfaces in $\mathbb R^3$ \cite{CHH22}, and in general by Bamler and Lai \cite{BL25a}, through their PDE-ODI principle \cite{BL25b}. Sun, Wang, and Xue \cite{SWX25} showed that the topological changes through nondegenerate cylindrical singularities can be described by level set changes of a Morse function. See also the recent paper by Daniels-Holgate--Hershkovits--Sz\'ekelyhidi \cite{DHHS26} which shows that a generic compact mean curvature flow in $\RR^3$ encounters only finitely many singularities.

For the blowdown soliton, our results provide analogues of strong rigidity and of a canonical description of the flow near a singularity. \zcref[S]{cor:rigidity} is its strong rigidity, and \zcref[S]{thm:A, thm:B} describe the flow near its singularities. The analysis is simpler than for cylinders, because all the neutral and destabilizing modes of the blowdown soliton are generated by scaling and diffeomorphisms \cite{NO25a}.

\subsubsection{Harmonic map heat flow and gauge fixing}
Gauge fixing is a major difficulty in the analysis of Ricci flow. Bamler--Kleiner \cite{BK22} use harmonic map heat flow to prove uniqueness and stability of three-dimensional Ricci flows through singularities, while Choi--Lai \cite{CL25} use it to establish exponential convergence of ancient Ricci flows asymptotic to compact integrable shrinkers. We likewise use harmonic map heat flow to fix the gauge, obtaining a strictly parabolic Ricci-DeTurck flow that allows us to describe the behavior of the flow at the appropriate scale and under a suitable identification.

\subsection{Statement of the main results}
Let $(M, \overline{g}, f)$ be the blowdown soliton. 
For $l>0$, let $M_l:=\{2\sqrt{f} < l\}\subset M$, and let $\omega_4$ be the volume of the Euclidean unit ball in $\mathbb R^4$.

\subsubsection{Strong stability and genericity of blowdowns}

\begin{theoremx}\label{thm:A}
There exist $\Gamma>0$, $\epsilon>0$, and $0 < c < 1/4$ such that the following holds. Let $T_0 > 0$. Suppose that a closed Riemannian manifold $(N, G_0)$ admits an open set $U\subset N$ and a diffeomorphism $\iota\colon U\to M_{100\Gamma}$ for which
\begin{equation}\label{eq:A1}
\|T_0^{-1}G_0-\iota^*\overline{g}\|_{C^{13}(S^2T^\ast U;\iota^*\overline{g})}<\epsilon
\end{equation}
and, for every $x\in N\setminus \iota^{-1}(M_\Gamma)$,
\begin{equation}\label{eq:A2}
\sup_{B_{\sqrt{T_0}r_0}^{G_0}(x)}|\Rm_{G_0}|_{G_0}< \dfrac{1}{T_0r_0^2} \quad\text{and}\quad\Vol_{G_0}B^{G_0}_{\sqrt{T_0}r_0}(x)>\dfrac{1}{2}\omega_4T_0^2 r_0^4
\end{equation}
where $r_0 = c\Gamma$.
Then the maximal Ricci flow $(G(t))_{t\in[0,T_G)}$ starting from $G_0$ satisfies $|T_G-T_0| < \frac{T_0}{10}$, and it develops a finite time singularity modeled on the blowdown soliton.
\end{theoremx}
\noindent The hypotheses are invariant under parabolic rescaling, and they are local. Condition \eqref{eq:A1} asks the metric to be close to the blowdown soliton at scale $\sqrt{T_0}$ on a large compact set, while \eqref{eq:A2} only asks it to have controlled geometry at the slightly larger scale $\sqrt{T_0}\Gamma$ elsewhere. The conical region $M_{100\Gamma} \setminus M_{8\Gamma}$ serves as a buffer: by pseudolocality, nothing outside of it can affect the formation of the singularity.

We say that a Ricci flow develops a \emph{blowdown singularity} at time $T$ if its curvature blows up at time $T$ and every tangent flow there is the blowdown soliton. Since \eqref{eq:A1} and \eqref{eq:A2} are open conditions, and by continuity of Ricci flow with respect to initial data (as long as it remains smooth), \zcref[S]{thm:A} has the following consequence.

\begin{corollary} \label{cor:open}
Let $(N^4, G_0)$ be a closed Riemannian manifold. Let $(N, G(t))_{t \in [0, T_G)}$ be the Ricci flow starting from $G_0$, with the first singular time $T_G < \infty$. Suppose that there exist $t_i \nearrow T_G$ and $p_i \in N$ such that 
\[
    (N, (T_G-t_i)^{-1}G(t_i), p_i) \longrightarrow (M, \overline{g}, o)
\]
in the smooth pointed Cheeger--Gromov sense. Suppose that there is no secondary curvature blowup at $T_G$, in the following sense:
\[
    \lim_{R \to +\infty} \limsup_{i \to +\infty} (T_G-t_i) \sup_{N \setminus B_i} |\Rm_{G(t_i)}|_{G(t_i)} = 0
\]
where $B_i = B_{\sqrt{T_G-t_i}R}^{G(t_i)}(p_i)$. Then $G_0$ has a $C^{14}$-open neighborhood $\cU$ such that every metric in $\cU$ develops a blowdown singularity at its first singular time.
\end{corollary}

\subsubsection{Strong rigidity of the blowdown soliton}

As a corollary, we obtain strong rigidity of the blowdown soliton. Strong rigidity was previously known for cylinders, in mean curvature flow \cite{CIM15} and in Ricci flow \cite{CM25,LW24,LZ25}.

\begin{corollary}\label{cor:rigidity}
There exist $\Gamma>0$ and $\epsilon>0$ such that the following holds. If a complete shrinking Ricci soliton $(M',g',f')$, normalized as the blowdown soliton, admits an open set $U'\subset M'$ and a diffeomorphism $\iota\colon U'\to M_{100\Gamma}$ with $\|g'-\iota^*\overline{g}\|_{C^{13}(S^2T^\ast U';\iota^*\overline{g})}<\epsilon$, then $(M',g')$ is isometric to $(M,\overline{g})$.
\end{corollary}
\noindent As for cylinders, \zcref[S]{cor:rigidity} only assumes closeness on a large compact set. Kotschwar and Wang, \cite{KW15}, prove a rigidity from infinity: the asymptotic cone determines the soliton, and \zcref[S]{cor:rigidity} can be seen as a complementary rigidity from a compact region. Using an estimate of Li and Wang \cite{LW20}, we show that such a soliton is asymptotically conical, with asymptotic cone close to $\mathcal C$, so that \zcref[S]{thm:A} applies with $T_0=1$. Its self-similar flow then converges to the blowdown soliton after rescaling by \zcref[S]{thm:B}(1), while all its rescalings are isometric to $(M',g')$.

Let $\overline\Sigma\subset M=\mathcal{O}(\pm1)$ be the zero section, and $(\mathcal C,g_{\mathcal C})$ the asymptotic cone of $(M,\overline{g})$. The flows of \zcref[S]{thm:A} have the following canonical structure near their singularity.

\begin{theoremx}\label{thm:B}
Let $(N,G_0)$, $U$ and $\iota$ be as in \zcref[S]{thm:A}, and let $T_G$ be the singular time.
\begin{enumerate}
\item \emph{Uniqueness of the tangent flow.} All of the rescaled flows $(T_G-t)^{-1}G(t)$, pointed at the maximum of curvature, converge to $\overline{g}$ in the pointed Cheeger--Gromov sense as $t\nearrow T_G$.
\item \emph{Conical singular time slice.} As $t \nearrow T_G$, the metrics $G(t)$ converge smoothly on every compact subset of the limit, under the identifications constructed in the proof. Its metric completion $G(T_G)$ is a compact space $Z$ with one isolated conical singularity $x_0$, modeled on $\cC$. 
\item \emph{Continuation by a blowdown.} Suppose that $N$ is closed, and let $N'$ be obtained by replacing a neighborhood of the sphere $\iota^{-1}(\overline\Sigma)$ by a ball, so that $N\approx N'\#\overline{\mathbb{CP}}{}^2$ or $N\approx N'\#\mathbb{CP}^2$, according to the sign of its self-intersection. Gluing in the blowdown expanding soliton of \cite{FIK03} as in \cite{GS18}, one obtains Ricci flows on $N'$ coming out of $G(T_G)$. They have the following properties:
\begin{itemize}
\item their curvature is bounded by $C(t-T_G)^{-1}$,
\item they converge to $G(T_G)$ as $t\searrow T_G$, smoothly away from the singular point and in the Gromov-Hausdorff sense globally,
\item their forward tangent flow at the singular point is the blowdown expanding soliton, and
\item if the singular time slices converge, these continuations subconverge.
\end{itemize}
\end{enumerate}
\end{theoremx}

\begin{conjecture}\label{conj:canonical}
The continuation in \zcref[S]{thm:B}(3) does not depend on the choices made in the gluing. It is the only Ricci flow coming out of $G(T_G)$ with $\mathrm{Rm}\ge -o(|t-T_G|^{-1})$ near the singular point.
\end{conjecture}
\noindent Rescalings of such a flow at the singular point subconverge to a Ricci flow coming out of $\mathcal C$ with nonnegative curvature operator, which is an expanding soliton by Chan, Lee and Peachey \cite{CLP24}. The expanding soliton asymptotic to $\mathcal C$ with nonnegative curvature operator is unique \cite{Der16,GS18}. Moreover, the blowdown expanding soliton is linearly stable \cite{DL17,NO25b,LT26}. The difficulty lies in passing from the rescalings to the flow itself, which already requires substantial work in the K\"ahler case \cite{Che25,CHL26}.

\begin{corollary}\label{cor:entropy}
Let $(N,G_0)$ be closed and as in \zcref[S]{thm:A}, with singular time $T_G$. Then $\mu(G_0,T_G)<\mu(\overline{g},1)$, where $\mu$ denotes Perelman's entropy.
\end{corollary}
\noindent In this sense, the blowdown soliton is a \emph{strict local maximizer of Perelman's entropy}, among metrics that are close to it on a large compact set and have controlled geometry elsewhere.

\begin{remark}\label{rem:CHI}
\zcref[S]{thm:A, thm:B}(1)-(2) only use that the blowdown soliton is asymptotically conical and linearly stable in the sense of \cite{CHI04}, as shown in \cite{NO25a}. They extend to all such shrinkers, but the only currently known example is that of \cite{FIK03}. The extension \zcref[S]{thm:A, thm:B}(3) using \cite{GS18} additionally requires that the asymptotic cone has nonnegative Riemannian curvature. 
\end{remark}

\subsubsection{Ricci flows with prescribed ordered surgeries and topological consequences}

 \zcref[S]{thm:A} also gives a very simple way to construct Ricci flows with prescribed singularities: it suffices to glue small copies of the blowdown soliton into a given metric, and \zcref[S]{thm:A} describes their evolution. This contrasts with previous constructions of Ricci flows with prescribed behavior, which rely on delicate parabolic gluing arguments and require controlling the whole flow \cite{BK17,DO24,Sto26,Lot25}.

\begin{theoremx}\label{thm:C}
Let $X$ be a closed oriented four-manifold, and let $k,l \ge 0$. Fix an ordering of the $k+l$ summands of $X\#k\mathbb{CP}^2\#l\overline{\mathbb{CP}}{}^2$. Then there is a metric on $X\#k\mathbb{CP}^2\#l\overline{\mathbb{CP}}{}^2$ whose Ricci flow, continued as in \zcref[S]{thm:B}(3) at each singular time, for any choice of these continuations, has the prescribed first $k+l$ singular times after which the flow lives on $X$.
\end{theoremx}
\noindent By the results of Freedman and Donaldson recalled above, every closed, simply connected, smooth, non-spin four-manifold is homeomorphic to a connected sum $S^4\#k\mathbb{CP}^2\#l\overline{\mathbb{CP}}{}^2$. By \zcref[S]{thm:C} with $X=S^4$, every such connected sum carries Ricci flows that blow down all its summands, in any order, and end on $S^4$.

Taking instead for $X$ an exotic manifold which becomes standard after one specific connected sum with $\mathbb{CP}^2$, Ricci flow creates exotic smooth structures.

\begin{corollary}\label{cor:exotic}
Let $Y:=2\mathbb{CP}^2\#9\overline{\mathbb{CP}}{}^2$. There exist closed four-manifolds $X_j$, $j\in \mathbb{N}$, all homeomorphic to $\mathbb{CP}^2\#9\overline{\mathbb{CP}}{}^2$ but pairwise non-diffeomorphic, and pairwise disjoint nonempty $C^{14}$-open sets $\mathcal U_j$ of metrics on $Y = X_j\#\mathbb{CP}^2$, such that for every $g_j\in\mathcal U_j$, the Ricci flow starting at $g_j$ develops a blowdown singularity at its first singular time, and its continuations given by \zcref[S]{thm:B}(3) live on $X_j$. 
\end{corollary}
\noindent The $X_j$ are the Dolgachev surfaces $E(1)_{2,2j+1}$, obtained from the rational elliptic surface $\mathbb{CP}^2\#9\overline{\mathbb{CP}}{}^2$ by logarithmic transforms of orders $2$ and $2j+1$ on two regular fibers \cite{Dol81}. By Freedman's classification \cite{Fre82}, they are all homeomorphic to $\mathbb{CP}^2\#9\overline{\mathbb{CP}}{}^2$. However, they are pairwise non-diffeomorphic \cite{Don87,FM88}. In other words, $Y$ contains infinitely many spheres of self-intersection $+1$ whose blowdowns are pairwise non-diffeomorphic, and Ricci flow blows down each of them from an open set of initial metrics. The smooth manifold produced by the surgery is thus determined by the initial metric, not by the initial manifold. There are infinitely many other such examples $Y$ by \cite{Gom91}.

\begin{corollary} \label{cor:exotic-stab}
Let $Z^4$ be a closed simply connected smooth 4-manifold possibly with an exotic smooth structure. Then there exist integers $k, l \ge 0$ and a Riemannian metric on the standard connected sum $k \mathbb{CP}^2 \# l\overline{\mathbb{CP}}{}^2$ whose Ricci flow ends on a manifold diffeomorphic to $Z^4$ after finitely many blowdown singularities which were continued by \zcref[S]{thm:B}(3) at each singular time, regardless of the choices.
\end{corollary}
\noindent Wall's Stabilization Theorem \cite{Wal64a, Wal64b} gives that for any closed simply connected 4-manifold $Z^4$, there exists $m \ge 0$ such that $Z^4 \# m\mathbb{CP}^2 \# m \overline{\mathbb{CP}}{}^2$ is diffeomorphic to $(m + b_2^+(Z)) \mathbb{CP}^2 \# (m + b_2^-(Z)) \overline{\mathbb{CP}}{}^2$. Hence, \zcref[S]{cor:exotic-stab} follows from \zcref[S]{thm:C}.

\begin{remark}\label{rem:Ric}
By M\'aximo \cite{Max14}, $\mathrm{Ric}_{\overline{g}}$ has a negative eigenvalue. Hence the Ricci curvature of the flows of \zcref[S]{thm:A} is unbounded below as $t\nearrow T_G$, which gives an open set of initial data whose Ricci flows develop Ricci curvature unbounded below.
\end{remark}

\begin{remark}\label{rem:Kahler}
By \cite{CHM25,CCD26}, the K\"ahler-Ricci flow on compact K\"ahler surfaces provides many flows satisfying the hypotheses of \zcref[S]{thm:A} shortly before their singular time. By \zcref[S]{thm:A}, all small perturbations of these flows, K\"ahler or not, still develop blowdown singularities.
\end{remark}

\begin{remark}\label{rem:infinite}
 Applying \zcref[S]{thm:C} to $X=\mathbb T^4$ and lifting the flow to the universal cover gives Ricci flows with bounded curvature on $\mathbb R^4\#\infty\mathbb{CP}^2\#\infty\overline{\mathbb{CP}}{}^2$, which blow down infinitely many spheres, simultaneously at each of finitely many singular times, and end on $\mathbb R^4$.
\end{remark}

We conclude this section with a question regarding unstable asymptotically conical shrinking solitons.

\begin{question}\label{q:codim}
Let $(M,g,f)$ be an asymptotically conical shrinker with $d$ destabilizing modes that are not generated by scaling and diffeomorphisms. Is the set of initial data developing a singularity modeled on $(M,g,f)$ a submanifold of codimension $d$?
\end{question}

A related dynamical picture for mean curvature flow was shown by Colding and Minicozzi \cite{CM21}.

Stolarski \cite{Sto26} constructed closed Ricci flows developing singularities modeled on asymptotically conical shrinkers, by tuning one parameter for each unstable mode, and \cite{BL25a,BL25b} develop tools to handle neutral and destabilizing modes for cylindrical singularities of mean curvature flow.

\subsection{Strategy of the proof} \label{subsec:strategy}
As usual with intrinsic geometric evolution equations, a key issue is gauge-fixing, especially here when comparing a compact flow and a noncompact singularity model $(M,\overline g, f)$. The proof selects the singular time and initial diffeomorphism by requiring the decay of the error between a region of the compact flow and the soliton $(M,\overline g, f)$. Harmonic map heat flow makes the comparison parabolic, but an arbitrary initial gauge need not yield decay \cite[Section 1.1, footnote 2]{CM25}. By \cite{NO25a}, recalled in \zcref[S]{prop:fiksa}, all nine nonnegative modes of $L_f:=\Delta_f+2\Rm_{\overline g}$ are generated by scaling and diffeomorphisms, which we use as degrees of freedom parametrized by $p\in \mathbb R^9$ for our initial gauge choice (\zcref[S]{prop:geolc}). 

More concretely, suppose that we are given an initial metric $G(0)$ satisfying the assumptions of \zcref[S]{thm:A}, namely $C^{13}$ closeness to a rescaled copy of $\overline g$ on a large compact region, together with curvature and volume bounds outside the core. We compare our Ricci flow $G(t)$ with
\[
    \overline G_p(t):=(T_p-t)\Psi_{p,t}^*\overline g,
    \qquad \Psi_{p,0}=\Xi_p,
\]
where $\Psi_{p,t}$ is generated by $(T_p-t)^{-1}\nabla^{\overline g}f$. Thus $T_p$ is the candidate singular time and $\Xi_p$ fixes the initial identification. We eventually select the parameter $p$ by adapting Stolarski's Wa\.zewski argument \cite[Section 8.3]{Sto26}. It is important to note that Stolarski's parameters perturb the initial metric, whereas ours change only the diffeomorphism and scale used in the comparison with the flow $\overline{G}(t)$ generated by the shrinking soliton $(M,\overline g, f)$. 

We glue $\overline G_p(t)$ to $G(t)$ in this region and use the resulting metric as the target of harmonic map heat flow. The asymptotic conicality of the blowdown soliton creates a buffer region that completely localizes the analysis of the soliton by pseudolocality, contrarily to the cylindrical case. 

The Wa\.zewski box requires both spectral components to remain below $\mu e^{-\gamma\tau}$, with $\lambda_{10}<-\gamma<0$ is determined by the spectral gap obtained in \cite{NO25a}. For the spectral controls, we adapt the localization mechanism underlying Bamler and Lai's PDE-ODI principle \cite[Theorem 4.17 and (4.19)]{BL25b}: the cutoff radius is chosen according to the decay rate being proved. Concretely, the choice \eqref{eqn:rtau}, with $R(\tau)^2$ affine in $\tau$, balances two effects: 1) the Gaussian cutoff error decays faster than $e^{-\gamma\tau}$, while 2) the loss in passing from weighted $L^2$ estimates to pointwise estimates  is absorbed by this decay. 

The proof relies on estimates which exclude every other exit. After the spectral parametrization in \zcref[S]{sec:spectral}, \zcref[S]{sec:gauge} defines the comparison and the box. \zcref[S]{sec:B4} controls the comparison on all of $N$, and \zcref[S]{sec:B3l,sec:B3e} improve the local $C^{12}$ bound at later and earlier times. \zcref[S]{sec:B2} excludes exits through the negative eigenmodes and proves outward crossing through the nonnegative eigenmodes. \zcref[S]{sec:parameters} identifies the admissible initial comparisons. Finally, \zcref[S]{sec:proofmain} shows that singularities of Ricci flow or of the harmonic map heat flow cannot cause an exit and completes the topological argument showing the existence of a parameter $p^*$ such that the flow remains in the box until it develops a singularity modeled on $(M,\overline g, f)$.

The rescaled flow converges exponentially on compact subsets to $(M,\overline g)$, and $T_{p^*}$ is its singular time. Thus a single choice of singular time and initial identification yields a decaying comparison throughout the evolution. \zcref[S]{sec:rigidity} proves strong rigidity, and \zcref[S]{sec:consequences} deduces the remaining results.

The proof applies verbatim to arbitrary linearly stable asymptotically conical shrinking Ricci solitons. 

\subsection*{Acknowledgements}

The authors would like to thank Hokuto Konno for discussions about the exotic examples that can be obtained by blowing down  spheres of self-intersection $\pm 1$. DK would like to thank Bill Minicozzi for constant support and invaluable discussions. DK also thanks Richard Bamler for inspiring discussions. The authors would like to thank Tang-Kai Lee for bringing several references to their attention. DK is partially supported by Simons Dissertation Fellowship and NSF Grant DMS-2304684.
TO is supported by the National Science Foundation under Grant No. DMS-2405328.

\subsection*{Use of AI tools}

The authors only used LLMs at the very last stage of the project for proofreading purposes: to verify constant dependencies and to fix a few mathematical and grammatical typos.

\section{Weighted Spectral Analysis} \label{sec:spectral}

Throughout the paper, let $(M, \overline{g}, f)$ be the FIK shrinking Ricci soliton, with the potential function $f$ normalized by
\[
    \Ric_{\overline{g}} + \Hess_f = \dfrac{1}{2} \overline{g} \quad \text{and} \quad R_{\overline{g}} + |\nabla^{\overline{g}} f|_{\overline{g}}^2 - f = 0.
\]
It is known that $f > 0$ on $M$.
For $x \in M$, define
\[
    \overline{b}(x) := 2\sqrt{f(x)}.
\]
For each $l>0$, denote the corresponding sublevel set by 
\[
    M_l := \{x \in M \colon \overline{b}(x) < l \} = \{x \in M \colon f(x) < l^2/4\}.
\]

\subsection{Normalized main theorem}
Let $\omega_4$ denote the volume of the unit ball in $\mathbb{R}^4$. Since $(M, \overline{g}, f)$ is asymptotically conical, we may also fix a constant $0<c<1/4$ such that for any $l > 0$, $|\Rm_{\overline{g}}|_{\overline{g}} < (2cl)^{-2}$ on $M \setminus M_l$.

To avoid adding unnecessary complexity, we fix the scale of \zcref[S]{thm:A} and instead prove the following statement. Note that the two statements are equivalent. 
\begin{theorem} \label{thm:mainthm}
    There exist $\Gamma > 0$ and $\epsilon > 0$ with the following property. Suppose that $(N, G_0)$ is a closed Riemannian manifold admitting an open set $U \subset N$ and a diffeomorphism $\iota \colon U \to M_{100\Gamma}$ for which 
   \begin{align}
        &\norm{G_0 - \iota^\ast \overline{g}}_{C^{13}(S^2 T^\ast U; \iota^\ast\overline{g})} < \epsilon, \label{eqn:cl1}\\
        &\sup_{B_{r_0}^{G_0}(x)} |\Rm_{G_0}|_{G_0} \le \dfrac{1}{r_0^2} \quad \text{and} \quad \Vol_{G_0} B_{r_0}^{G_0}(x) \ge \dfrac{1}{2}\omega_4r_0^4 \quad \text{for } x\in N \setminus \iota^{-1} (M_\Gamma) \text{ and } r_0 = c\Gamma. \label{eqn:cl2}
    \end{align} 
    Then the maximal Ricci flow $(N, G(t))_{t \in [0, T_G)}$ with initial data $(N, G_0)$ develops a finite-time singularity modeled on the blowdown soliton $(M, \overline{g}, f)$ and $|T_G-1| < 1/10$. More precisely, there exist $x_t \in N$ for $t \in [0, T_G)$ and $o \in M$ such that 
    \[
        (N, (T_G-t)^{-1}G(t), x_t) \to (M, \overline{g}, o)
    \]
    as $t \nearrow T_G$ in the smooth pointed Cheeger--Gromov sense. 
\end{theorem}

Denote by $L_f^2 := L^2(S^2T^\ast M, e^{-f}d\mu_{\overline{g}})$ the space of symmetric 2-tensors with finite weighted $L^2$ norm. A standard spectral argument (cf. \cite[Theorem 2.18]{Sto26}) shows that in the distributional sense, the Lichnerowicz Laplacian $L_f := \Delta_f + 2\Rm_{\overline{g}}$ has discrete spectrum $\lambda_1 \ge \lambda_2 \ge \cdots$ and finite dimensional eigenspaces. 

For the FIK shrinking Ricci soliton, the nonnegative part of this spectrum consists entirely of modes generated by rescaling and diffeomorphisms.
\begin{proposition}[{\cite[Remark 2]{NO25a}}] \label{prop:fiksa}
On the FIK shrinking Ricci soliton, there are nine nonnegative eigenmodes of the Lichnerowicz Laplacian $L_f$, with eigenvalues and eigentensors given by
\begin{itemize}
    \item $\lambda_1 = 1$ with eigentensor $\overline{h}_1=\Ric_{\overline{g}}$,
    \item $\lambda_2 = \cdots =\lambda_5 = 1 - \dfrac{1}{\sqrt 2}$ with orthogonal eigentensors $\overline{h}_i = \nabla^{\overline{g}, 2} \overline{w}_i$ for some functions $\overline{w}_i \in C^\infty(M)$ for $2 \le i \le 5$,
    \item $\lambda_6 = \cdots = \lambda_9 = 0$ with orthogonal eigentensors $\overline{h}_i=\nabla^{\overline{g}, 2} \overline{w}_i$ for some functions $\overline{w}_i \in C^\infty(M)$ for $6 \le i \le 9$,
\end{itemize}
All remaining eigenvalues are negative. We fix $\overline{w}_9 = f$. Moreover, all eigenmodes $\overline{h}_1, \ldots, \overline{h}_9$ are bounded.

\begin{remark} \label{rem:ssct}
To prove \zcref[S]{thm:A}, we only require $(M, \overline{g}, f)$ to be a nonflat asymptotically conical shrinking Ricci soliton whose nonnegative eigenspace is spanned by the Ricci tensor and Lie derivatives of the metric along vector fields.
\end{remark}

\end{proposition}
We denote by $V_+$ the space of unstable modes and $V_-$ its orthogonal complement. For our purposes, unstable modes include zero eigenmodes. 
\[
    V_+ := \mathrm{Span}\{\overline{h}_1, \ldots, \overline{h}_9\} \quad \text{and} \quad V_- := V_+^\perp \subset L_f^2.
\]
Thus, $L_f^2 = V_+ \oplus V_-$. We define $\mathcal{P}_+ \colon L_f^2 \to V_+$ and $\mathcal{P}_- \colon L_f^2 \to V_-$ to be the projections onto the space of unstable modes $V_+$ and stable modes $V_-$, respectively.

Throughout the paper, we fix $\gamma > 0$ so that 
\begin{equation} \label{eqn:gamma}
    \lambda_{10} < -\gamma < 0.
\end{equation}

\subsection{Geometric Variations}
We perturb the background metric in unstable directions through parameters $p \in \mathfrak{P} \subset \mathbb{R}^9$. The set $\mathfrak{P}$ will be specified later. For $p = (p_1, \ldots, p_9) \in \mathfrak{P}$, define the coefficients by
\[
    a(p) = -\dfrac{1}{2}p_1,\; b_i(p) = -p_i \; (2 \le i \le 8), \; b_9(p) = p_1-p_9.
\]
In terms of the eigentensors $\overline{h}_i$, these coefficients satisfy
\begin{equation} \label{eqn:lcpm}
    a(p) \overline{g} + \sum_{i=2}^9 b_i(p) \nabla^{\overline{g}, 2} \overline{w}_i = -\sum_{j=1}^9 p_j\overline{h}_j
\end{equation}
since $\overline{g} = 2(\Ric_{\overline{g}} + \Hess_f)$.
Let $X_p := \frac12\sum_{i=2}^9 b_i(p)\nabla^{\overline g}\overline w_i$, so that $\mathcal L_{X_p}\overline g=\sum_{i=2}^9 b_i(p)\nabla^{\overline g,2}\overline w_i$, and let $\Xi_p$ be the time-one flow of $X_p$, which is complete since $|X_p|$ grows at most linearly. Let $T_p := 1 + a(p)$.

The initial metric $G_0$ remains fixed throughout the argument. The parameter $p$ changes only the comparison soliton: its first component changes the candidate singular time, and the other components change the identification by diffeomorphisms. The coefficients above are chosen so that the linearized comparison error is $\sum_{j=1}^9 p_j \overline{h}_j$. The next proposition quantifies the quadratic error in this parametrization.

\begin{proposition}\label{prop:geolc}
    There exists $C_{\cR} > 0$ with the following property.
    Let $0 < p_0 \le \overline{p_0}$ for a small enough $0 < \overline{p_0} \ll 1$, $\Gamma \ge \underline{\Gamma}$ for a large enough $\underline{\Gamma} \gg 1$,  and $0 < \epsilon \le \overline{\epsilon}$ for a small enough $0 < \overline{\epsilon} \ll 1$.
    Suppose that $(N, G_0)$ is a closed Riemannian manifold satisfying \eqref{eqn:cl1} and \eqref{eqn:cl2}. Suppose $p \in \mathbb{R}^9$ with $|p| \le p_0$. Then 
    \[
        T_p^{-1}(\Xi_p^{-1})^\ast ((\iota^{-1})^\ast G_0) - \overline{g} = \sum_{j=1}^9 p_j\overline{h}_j + \cR(p)
    \]
    where $\norm{\cR(p)}_{C^{13}(M_\Gamma)} \le C_{\cR}(|p|^2+\epsilon)$ and $\norm{D_p\cR(p)}_{\mathbb{R}^9 \to C^0(M_\Gamma)} \le C_{\cR}(|p|+\epsilon\Gamma)$.
\end{proposition}

\begin{proof}
    For simplicity, we write $(\iota^{-1})^\ast G_0 = \overline{g} + q$, where $q$ is regarded as an error term. Write $\phi_{p,s}$ for the time-$s$ flow of $-X_p$ so that $\phi_{p,1} = \Xi_p^{-1}$. Since $|X_p| \le C|p|(1+\overline{b})$, for sufficiently small $\overline{p_0}$ and sufficiently large $\underline{\Gamma}$, 
    \[
        \phi_{p,s}(M_\Gamma) \subset M_{2\Gamma} \subset M_{100\Gamma}
    \]
    for $0 \le s \le 1$, $|p| \le p_0 \le \overline{p_0}$, and $\Gamma \ge \underline{\Gamma}$.

    Write $\cR = \cR_s + \cR_e$ where 
    \[
        \cR_s := T_p^{-1}\phi_p^\ast \overline{g} - \overline{g} - \sum_{j=1}^9 p_j\overline{h}_j \quad \text{and} \quad \cR_e := T_p^{-1}\phi_p^\ast q.
    \]
    By Taylor's formula, 
    \[
        \cR_s(p) = \int_0^1 (1-t)D^2(\cR_s(tp))(p,p)dt,
    \]
    which gives $\norm{\cR_s(p)}_{C^0(M_\Gamma)} \le C_{\cR}|p|^2$. Higher derivatives follow similarly using $|\nabla^k X_p| \le C_k|p|(1+\overline{b})^{1-k}$. 
    For the error term $\cR_e$, the assumption \eqref{eqn:cl1} gives $\norm{\cR_e}_{C^{13}(M_\Gamma)} \le C_{\cR}\epsilon$.

    Note also that 
    \begin{align*}
        (D_p\cR_s(p))(v) &= \int_0^1 D_p^2 \cR_s(tp)(p,v)dt\\
        (D_p\cR_e(p))(v) &= -a(v)T_p^{-2}\phi_p^\ast q + T_p^{-1}\phi_p^\ast \cL_{Z_{p,v}}q,
    \end{align*}
    where $Z_{p,v} := (D_p\phi_p(v))\circ \phi_p^{-1} = -\int_0^1 (\phi_{p,s})_\ast X_vds$. Since $|Z_{p,v} | \le C_k|v|(1+\overline{b})$, one gets 
    the second assertion as well.
\end{proof}

Throughout the paper, $p_0 \le 1/2$ so that $T_p < 2$ always holds. Let $\Psi_{p,t}$, $t \in [0, T_p)$, be a 1-parameter family of diffeomorphisms on $M$ solving 
\[
    \partial_t \Psi_{p,t} = \dfrac{1}{T_p - t} \nabla^{\overline{g}}f \circ \Psi_{p,t} \quad \text{with} \quad \Psi_{p,0} = \Xi_{p}.
\]

\section{Harmonic Map Heat Flow Gauge and the Wa\.zewski Box} \label{sec:gauge}
Let $(N, G_0)$ be a closed Riemannian manifold, and denote its maximal Ricci flow by $(G(t))_{t \in [0, T_G)}$. To compare the Ricci flow $G(t)$ to the Ricci flow generated by the FIK shrinking Ricci soliton, we use a cutoff supported inside the region of $N$ identified with the soliton. 
Fix a smooth cutoff $\zeta \colon M \to [0,1]$ with
\[
    \zeta|_{\{x \colon \overline{b}(x) < 4\Gamma\}} = 1, \quad \zeta|_{\{x : \overline{b}(x) > 8\Gamma\}} = 0.
\]
Note that we may choose $\zeta$ so that 
\[
|\nabla^{\overline{g}, j} \zeta|_{\overline{g}} \le \dfrac{C_\zeta}{\Gamma^j}
\] 
for $1 \le j \le k$, where $C_\zeta = C_\zeta(n, \overline{g}, k) > 0$ is independent of $\Gamma$. Pulling back by $\iota$ and extending by zero outside $U$ defines the cutoff $\zeta^N$ on $N$.

For a parameter $p$, the soliton Ricci flow is $\overline{G}_p(t) := (T_p-t)\Psi_{p,t}^\ast \overline{g}$. We interpolate between $\overline{G}_p(t)$ and $G(t)$ to construct the time-dependent target metric $\beta_p(t)$.
\begin{align*}
    \beta_p(t) := \begin{cases}
        \zeta^N \iota^\ast \overline{G}_p(t) + (1-\zeta^N)G(t) &\text{on } U,\\
        G(t) &\text{outside } U.
    \end{cases}
\end{align*}
The target metric agrees with the soliton flow where $\zeta^N = 1$ and with the actual Ricci flow when $\zeta^N = 0$. It is generally not a Ricci flow in the transition annulus. Its Ricci flow defect will therefore enter the exterior estimates. Harmonic map heat flow supplies a common gauge in which the difference between the actual flow and this target can be estimated.

We fix the gauge using harmonic map heat flow from $G(t)$ to $\beta_p(t)$. Let $\Phi_{p} \colon N \times [0,T_\Phi) \to N$ solve
\begin{equation}
    \partial_t \Phi_p(\cdot, t) = \Delta_{G(t), \beta_p(t)}\Phi_p(\cdot,t ) \quad \text{with} \quad \Phi_p(\cdot, 0) = \mathrm{Id}_N.
    \label{eqn:hmhf}
\end{equation}
Write $\Phi_{p,t} := \Phi_p(\cdot, t)$. Let $T_\Phi$ be the maximal time such that $\Phi_p$ exists and $\Phi_{p,t}$ is a diffeomorphism for every $t \in [0,T_\Phi)$. The time $T_\Phi$ depends on $p$, but we suppress this dependence in the notation.
The pullback metric $A_p(t) := (\Phi_{p,t}^{-1})^\ast G(t)$ is the metric in this gauge, and $H_p(t) := A_p(t) - \beta_p(t)$ is its difference from the target metric on $N$. 

\subsection{Identifying the soliton region through cutoffs}
Since perturbation by rescaling lies in $V_+$, the candidate singular time depends on the parameter. We normalize the time as 
\[
    \tau_p(t) = -\log \left( \dfrac{T_p-t}{T_p} \right),
\]
and define $t_p(\tau)$ to be the inverse: $\tau_p(t_p(\tau)) = \tau$. Note that the time normalization depends on $p$.
Then we rescale $H_p(t)$, the difference between the gauged metric and the target metric, by 
\begin{equation} \label{eqn:hp}
    h_p(\tau) := \dfrac{1}{T_p-t} (\Psi_{p,t}^{-1})^\ast (\iota^{-1})^*H_p(t).
\end{equation}
The renormalization places the candidate singular time at infinite $\tau$. Although $H_p(t)$ is a tensor on $N$, the identification through $\iota$ allows a truncation of $h_p(\tau)$ to be considered on $M_{4R(\tau)}$, provided $\Gamma$ is sufficiently large and $R(\tau)$ grows sufficiently slowly.

The radius used in the estimates is fixed by 
\begin{equation} \label{eqn:rtau}
    R(\tau)^2 = R_0^2 + \dfrac{8\gamma}{\theta^2}\tau,
\end{equation}
where $0 < \theta<1$ will be chosen later in \eqref{eqn:theta-exp}. The boundary moves under $\Psi_{p,t}$ at rate $e^{\tau/2}$. Consequently, this choice of radius grows slowly enough to leave a comparison region $M_{4R(\tau)}$ on which the truncated rescaled perturbation is defined.
The choice \eqref{eqn:rtau} gives 
\[
    e^{-R(\tau)^2/8} = e^{-R_0^2/8}e^{-\gamma\tau/\theta^2}.
\]
Since $0<\theta<1$, this Gaussian decay is faster than the decay $e^{-\gamma \tau}$ of the spectral threshold, which allows the cutoff error to be absorbed uniformly in $\tau$.

To truncate the 2-tensor $h_p(\tau)$, we fix a cutoff $w_0 \colon [0, \infty) \to [0,1]$ satisfying
\begin{align*}
    w_0(x) &= 1 \quad \text{for all } 0 \le x \le 1 \quad \text{and}\\
    w_0(x) &= 0 \quad \text{for all } x \ge 2.
\end{align*}
Then take a family of cutoffs 
\begin{equation} \label{eqn:wtau}
    w(x, \tau) := w_0 \left( \dfrac{\overline{b}(x)}{R(\tau)} \right).
\end{equation}
Note that $\supp w \subset M_{2R(\tau)}$ and $w \equiv 1$ on $M_{R(\tau)}$.
We will often write $w_{R(\tau)}$ for the spatial cutoff $w(\cdot, \tau)$.

\subsection{Wa\.zewski Box}
We now define the Wa\.zewski box $\mathfrak{B}$. 
\begin{definition}[Box $\mathfrak{B}$] Given $0 < p_0 < 1$, $R_0 > 0$, $\Gamma > 0$, $0 < \theta < 1$, $\mu > 0$, and $\eta > 0$, we say $H_p(t_p(\tau)) \in \mathfrak{B}$ for $p \in \mathbb{R}^9$ with $|p| \le p_0$ and $\tau \ge 0$ if
\begin{enumerate}[label=(B\arabic*), ref=(B\arabic*)]
    \item\label{box:B1} $\norm{\mathcal{P}_+(w_{R(\tau)}h_p(\tau))}_{L_f^2} \le \mu e^{-\gamma \tau}$,
    \item\label{box:B2} $\norm{\mathcal{P}_-(w_{R(\tau)}h_p(\tau))}_{L_f^2} \le \mu e^{-\gamma \tau}$,
    \item\label{box:B3} $\norm{h_p(\tau)}_{C^{12}(M_{2R(\tau)})} \le \eta$,
    \item\label{box:B4} $(1-8\eta) \beta_p(t_p(\tau)) \le A_p(t_p(\tau)) \le (1+8\eta)\beta_p(t_p(\tau))$,
\end{enumerate}
where $t_p(\tau)$ is given by 
\begin{equation}
    \tau = -\log \left( \dfrac{T_p-t_p(\tau)}{T_p} \right).
    \label{eqn:tptau}
\end{equation}
Given $\tau_0 \in [0, \infty]$, we say $H_p(t_p(\tau))$ \emph{exits} the box $\mathfrak{B}$ at time $\tau = \tau_0$ if $\tau_0$ is the maximal time for which $H_p(t_p(\tau)) \in \mathfrak{B}$ for $\tau \in [0, \tau_0)$.
\end{definition}

The four conditions serve different purposes. Conditions \ref{box:B1} and \ref{box:B2} control the truncated perturbation in the weighted spectral decomposition. Condition \ref{box:B3} provides local regularity on the growing soliton region, while \ref{box:B4} keeps the gauged metric uniformly comparable to its target on the whole manifold. The argument will improve the last three bounds and show that equality in \ref{box:B1} forces an outward crossing.

The parameters $p$ change the scale and gauge of the comparison soliton while keeping the initial metric fixed. This is one distinction between our setup and Stolarski's \cite{Sto26}. The corresponding initial unstable modes form a closed ball. We seek to show that any finite exit occurs through the unstable mode condition \ref{box:B1} and that the crossing is strictly outward. If every parameter exited, the first-exit map would give a continuous retraction of the ball onto its boundary, contradicting the no retraction theorem. A parameter that never exits would give the desired soliton limit.

\section{{\ref{box:B4}} cannot cause an exit} \label{sec:B4}
The condition \ref{box:B4} controls the metric on all of $N$, including the region that is not identified with the FIK shrinking Ricci soliton. 

For $x \in U$ and $t \ge 0$, set
\[
    f^N_{p,t}(x) := f(\Psi_{p,t}(\iota(x))) \quad \text{and} \quad b^N_{p,t}(x) := 2\sqrt{f^N_{p,t}(x)}.
\]
The function $b_{p,t}^N(x)$ is well-defined on this range and increases in $t$ for each $x \in U$.

We fix constants
\begin{equation}
    \theta = \dfrac{7}{8}, \quad \theta' = \dfrac{5}{8}, \quad \text{and} \quad \delta_\theta = \dfrac{1}{2}.
    \label{eqn:theta-exp}
\end{equation}
These are chosen so that 
\begin{equation}
    0 < \theta' < \theta < 1 \quad \text{and} \quad (4+\delta_\theta)e^{-7/2} < \frac{1}{2}\theta'.
    \label{eqn:theta}
\end{equation}

For $R_0 > 0$ and $0 < T < T_p$, we define a space-time domain
\begin{equation} \label{eqn:nprt}
    \mathfrak{N}_{p, R_0, T} := \left\{(x,t) \in U \times [0,T] | b^N_{p,t}(x) \ge R(\tau_p(t)) \right\} \sqcup (N \setminus U) \times [0,T],
\end{equation}
where $t_p(\tau)$ is given by \eqref{eqn:tptau}.
Take $\underline{R_0} \gg 1$ large enough so that $|\nabla^{\overline{g}}f|_{\overline{g}}^2 > 2\gamma/\theta^2$ on $M\setminus M_{\underline{R_0}}$. Then for $R_0 \ge \underline{R_0}$, we have
\[
    \partial_t (f^N_{p,t})(x) = \dfrac{|df|_{\overline{g}}^2(\Psi_{p,t}(\iota(x)))}{T_p - t} > \dfrac{1}{T_p-t}
\]
for $x \in U$ with large enough $b^N_{p,t}(x)=\overline{b}(\Psi_{p,t}\iota(x))$. Therefore, if $(x, t) \in \mathfrak{N}_{p, R_0, T}$ then $(x, t') \in \mathfrak{N}_{p, R_0, T}$ for any $t' \in [t,T]$. Moreover, by taking $\Gamma = \Gamma(R_0) \gg 1$ large enough,
\[
    b^N_{p,0}(x) \ge 2R_0 \quad \text{on} \quad \partial U,
\]
which implies that $\partial U \times (0,T)$ does not meet the boundary of $\mathfrak{N}_{p, R_0, T}$. Hence we define the parabolic boundary $\partial \mathfrak{N}_{p, R_0, T}$ as
\begin{equation}
    \partial \mathfrak{N}_{p, R_0, T} := \left\{(x, t) \in U \times [0,T] | b^N_{p,t}(x) = R(\tau_p(t)) \right\} \cup \{(x,0) \in U \times \{0\} | b^N_{p,t}(x) \ge R_0 \} \cup (N \setminus U) \times \{0\}.
\end{equation}

The following is the main proposition of this section.
\begin{proposition}[{\ref{box:B4}} cannot cause an exit] \label{prop:aprhmhf}
    Let $R_0 \ge \underline{R_0}$ for large enough $\underline{R_0} \gg 1$, $0 < \eta \le \overline{\eta}$ for small enough $0 < \overline{\eta} \ll 1$, $\Gamma \ge \underline{\Gamma}$ for large enough $\underline{\Gamma} = \underline{\Gamma}(R_0, \eta) \gg 1$, and $0 < \epsilon \le \overline{\epsilon}$ for small enough $0 < \overline{\epsilon} = \overline{\epsilon}(\Gamma, R_0, \eta) \ll 1$.
    Suppose that $(N, G_0)$ is a Riemannian manifold admitting an open set $U \subset N$ and a diffeomorphism $\iota \colon U \to M_{100\Gamma}$ satisfying \eqref{eqn:cl1} and \eqref{eqn:cl2}, and that for some $T < \min\{2, T_G, T_p, T_\Phi\}$, $\Phi_p \colon N \times [0,T] \to N$ is a smooth solution to Harmonic Map Heat Flow \eqref{eqn:hmhf} satisfying 
    \begin{align*}
        (1-8\eta) \beta_p(t) \le (\Phi_{p,t}^{-1})^\ast G(t) \le (1+8\eta) \beta_p(t) \quad &\text{on } \mathfrak{N}_{p, R_0, T} \quad \text{and}\\
        \norm{(\Phi_{p,t}^{-1})^\ast G(t) - \beta_p(t)}_{\beta_p(t)} < \eta \quad &\text{on } \partial \mathfrak{N}_{p, R_0, T},
    \end{align*}
    then 
    \begin{equation}
        \norm{(\Phi_{p,t}^{-1})^\ast G(t) - \beta_p(t)}_{\beta_p(t)} < 4\eta \quad \text{on } \mathfrak{N}_{p, R_0, T}.
        \label{eqn:bh}
    \end{equation}
    Consequently, the comparison improves to
    \[
        (1-4\eta)\beta_p(t) \le (\Phi_{p,t}^{-1})^\ast G(t) \le (1+4\eta) \beta_p(t) \quad \text{on } \mathfrak{N}_{p, R_0, T}.
    \]
\end{proposition}

\subsection{A parabolic inequality for $H_p$}
Recall the notation $A_p(t) := (\Phi_{p,t}^{-1})^\ast G(t)$ and $H_p(t) := A_p(t) - \beta_p(t)$ on $N$. The harmonic map heat flow puts $A_p$ in Ricci-DeTurck gauge; see \cite[Chapter 3.4]{CK04}. To apply a maximum principle to $|H_p|_{\beta_p}^2$, we first obtain a parabolic inequality that $|H_p|_{\beta_p}^2$ satisfies.

\begin{lemma} \label{lem:rdtfc}
    For the DeTurck vector field $W_p^k = (A_p)^{ij} \left( (\Gamma_{A_p})_{ij}^k - (\Gamma_{\beta_p})_{ij}^k \right)$,
    \[
    -2(\Ric_{A_p} - \Ric_{\beta_p}) + \mathcal{L}_{W_p}A_p = (A_p)^{ab}\nabla^{\beta_p}_a \nabla^{\beta_p}_b H_p + \Rm_{\beta_p} \ast H_p + (A_p)^{-1} \ast (A_p)^{-1} \ast \nabla^{\beta_p} H_p \ast \nabla^{\beta_p} H_p.
    \]
\end{lemma}

\begin{proof}
    Write $A = A_p$, $\beta = \beta_p$, $W = W_p$, and $\nabla = \nabla^{\beta_p}$.
    The difference of the two connections is
    \[
        C_{ij}^k := (\Gamma_A)_{ij}^k - (\Gamma_\beta)_{ij}^k = \dfrac{1}{2} A^{kl} ( \nabla_i H_{jl} + \nabla_j H_{il} - \nabla_l H_{ij}).
    \]
    Then as $(\Ric_A - \Ric_{\beta})_{ij} = \nabla_k C_{ij}^k - \nabla_j C_{ik}^k + (C \ast C)_{ij}$, we get 
    \[
    -2(\Ric_A - \Ric_{\beta})_{ij} = A^{kl} (-\nabla_k \nabla_i H_{jl} - \nabla_k \nabla_j H_{il} + \nabla_k \nabla_l H_{ij} + \nabla_j \nabla_i H_{kl}) + (A^{-1} \ast A^{-1} \ast \nabla H \ast \nabla H)_{ij}.
    \]  
    The term contributed by the DeTurck vector field is
    \begin{align*}
        (\mathcal{L}_W A)_{ij}
        &= W^k \nabla_k A_{ij} + A_{kj}\nabla_i W^k + A_{ik}\nabla_j W^k\\
        &= A^{ab} \left(\nabla_i \nabla_a H_{bj} + \nabla_j \nabla_a H_{bi} - \frac{1}{2} (\nabla_i \nabla_j + \nabla_j \nabla_i) H_{ab}\right) + (A^{-1} \ast A^{-1} \ast \nabla H \ast \nabla H)_{ij}.
    \end{align*}
    Adding the two expressions cancels the mixed derivatives. Hence,
    \begin{align*}
        -2(\Ric_A - \Ric_{\beta}) + \mathcal{L}_W A_p
        &= A^{ab} \nabla_a \nabla_b H + \Rm_\beta \ast H + A^{-1} \ast A^{-1} \ast \nabla H \ast \nabla H. \qedhere
    \end{align*}
\end{proof}

The target metric need not be a genuine Ricci flow in the gluing region. Hence its defect $S_p := \partial_t \beta_p + 2\Ric_{\beta_p}$ must appear when differentiating the norm of $H_p$.
\begin{proposition} \label{prop:pi}
    There exist $0 < \overline{\eta} \ll 1$ and $C_R > 0$ such that, if $0 < \eta \le \overline{\eta}$ and 
    \[
        (1-8\eta)\beta_p(t) < A_p(t) < (1+8\eta)\beta_p(t)
    \]
    on $N \times [0,T]$, then
    \begin{equation}
        (\partial_t - (A_p)^{ab} \nabla_a^{\beta_p}\nabla_b^{\beta_p})|H_p|_{\beta_p}^2 \le C_R|\Rm(\beta_p)|_{\beta_p} |H_p|_{\beta_p}^2 + 2|S_p|_{\beta_p} |H_p|_{\beta_p}^2 + 2|S_p|_{\beta_p}|H_p|_{\beta_p}
        \label{eqn:pi}
    \end{equation}
    on $N \times [0,T]$, where $S_p := \partial_t \beta_p + 2\Ric_{\beta_p}$.
\end{proposition}

\begin{proof}
    Write $A = A_p$, $\beta = \beta_p$, $H = H_p$, $W = W_p$, $S = S_p$, and $\nabla^\beta = \nabla^{\beta_p}$. Since $\partial_t A = -2\Ric_A + \mathcal{L}_W A$ and $\partial_t \beta = -2\Ric_\beta + S$, 
    \begin{align}
    \partial_t |H|_{\beta}^2 
    &= 2\langle \partial_t H, H \rangle_{\beta} - 2\langle \partial_t \beta, H \circ H \rangle_{\beta_p} \nonumber\\
    &= 2\langle -2(\Ric_{A} - \Ric_{\beta}) + \mathcal{L}_WA - S, H \rangle_{\beta_p} - 2\langle -2\Ric_{\beta} + S, H \circ_{\beta} H \rangle_{\beta}, \label{eqn:pit}\\
    A^{ab}\nabla_a^{\beta}\nabla_b^{\beta} |H|_\beta^2 &= 2\langle A^{ab}\nabla_a^\beta \nabla_b^\beta H, H \rangle_\beta + 2A^{ab}\langle \nabla_a^\beta H, \nabla_b^\beta H \rangle_\beta. \label{eqn:pil}
    \end{align}
    The last term in \eqref{eqn:pil} is positive and coercive. For $\overline{\eta}$ sufficiently small, it absorbs the pairing of $H$ with the quadratic gradient term in \zcref[S]{lem:rdtfc}. The remaining curvature terms are bounded by $|\Rm_{\beta}||H|_\beta^2$ up to a multiplicative constant, while the two defect pairings contribute the last two terms of \eqref{eqn:pi}. This proves the desired inequality.
\end{proof}

\subsection{A supersolution on the three exterior regions}
Now we construct a barrier. Let $\chi$ be a smooth cutoff function on $M$ such that $\chi \equiv 1$ on $M_\Gamma$ and $\chi \equiv 0$ on $M \setminus M_{2\Gamma}$. We may choose $\chi$ so that 
\[
    |\nabla^{\overline{g}, j} \chi|_{\overline{g}} \le \dfrac{C_\chi}{\Gamma^j}
\]
for $0 \le j \le k$, where $C_\chi = C_\chi(k) > 0$ does not depend on $\Gamma$. Then denote by $\chi^N$ the cutoff $\iota^\ast \chi$ extended to zero outside $U$.

The auxiliary cutoff $\chi$ has its transition region strictly inside the region where $\zeta = 1$. Consequently, the target remains an exact soliton flow whenever $\chi$ is nonzero, and the Ricci flow defect of the target is supported in the outer region $\{\chi = 0\}$. This separation allows the barrier to handle the soliton curvature near the inner boundary and the gluing defect by different estimates.

\begin{lemma} \label{lem:rde}
    Given $C_S' > 0$, there exists $C_S > 0$ with the following property. Let $\Gamma \ge \underline{\Gamma}$ for large enough $\underline{\Gamma} \gg 1$. Let $(N, G_0)$ be a closed Riemannian manifold and $(g_1(t))_{t \in [0,T)}$ and $(g_2(t))_{t \in [0,T)}$ Ricci flows on $N$. Suppose that on $U \subset N$ the metrics $g_1(t)$ and $g_2(t)$ are identified with the metrics on $\{3\Gamma < \overline{b} < 9\Gamma\} \subset M$ and that 
    \[
        \dfrac{\overline{g}}{C_S'} \le g_i(0) \le C_S' \overline{g} \quad \text{and} \quad \Gamma^j|\nabla^{\overline{g}, j}g_i(0)|_{\overline{g}} \le C_S'
    \]
    for $i=1, 2$ and $j=1, 2$. Define $g_0 = \zeta g_1 + (1-\zeta) g_2$. Then 
    \[
        |\Rm_{g_0(t)}|_{g_0(t)} + |S(t)|_{g_0(t)} \le \dfrac{C_S}{\Gamma^2}
    \]
    for $t \in [0,2] \cap [0,T)$ on $\{4\Gamma < \overline{b} < 8\Gamma\} \subset M$.
\end{lemma}
\begin{proof}
    Throughout the proof, the constant $C$ may increase.
    Since  $g_1$ and $g_2$ are Ricci flows, we have 
    \begin{align*}
        S_p &= \partial_t g_0 + 2\Ric_{g_0}
        = \zeta \partial_t g_1 + (1-\zeta) \partial_t g_2 + 2\Ric_{g_0}
        = 2\left[ \Ric_{g_0} - \zeta \Ric_{g_1} - (1-\zeta)\Ric_{g_2} \right].
    \end{align*}
    Note that 
    \begin{equation} \label{eqn:S-est-Rm}
        (\Rm_g)_{abcd} = g_{ap}(\Rm_{\overline{g}})^p{}_{bcd} - \nabla_c^{\overline{g}}B_{bda} + \nabla_d^{\overline{g}} B_{bca} + g^{pq}(B_{acp}B_{bdq} - B_{adp}B_{bcq})
    \end{equation}
    where $B_{abc} := g_{cp}A^{p}{}_{ab}$ with $A(X,Y) = \nabla_X^g Y - \nabla_X^{\overline{g}}Y$ the connection difference. But also, due to a formula for the difference of Levi-Civita connections, 
    \[
        B_{abc} = \dfrac{1}{2}(\nabla_a^{\overline{g}}g_{bc} + \nabla_b^{\overline{g}}g_{ac} - \nabla_c^{\overline{g}}g_{ab}).
    \]
    Therefore, the first three terms in the expression \eqref{eqn:S-est-Rm} are linear in $g$. That is, in the expansion of 
    \[
        S'(t) := \Rm_{g_0(t)} - \zeta \Rm_{g_1(t)} - (1-\zeta)\Rm_{g_2(t)}
    \]
    using \eqref{eqn:S-est-Rm}, the second derivatives of $g_0$ cancel. In particular, 
    \[
        |S'| \le C(|\nabla^{\overline{g}}\zeta| |\nabla^{\overline{g}}(g_1-g_2)| + |\nabla^{\overline{g},2} \zeta| |g_1-g_2| + |\nabla^{\overline{g}}g_1|^2 + |\nabla^{\overline{g}}g_2|^2 + |\Rm_{\overline{g}}|(|g_1| + |g_2|))
    \]
    for a fixed constant $C>0$. Note that $\partial_t g_i = -2\Ric_{g_i}$. Therefore, to control $|g_1-g_2|$ and $|\nabla^{\overline{g}}(g_1-g_2)|$, it is enough to control $\int_0^{\min\{2,T\}} |\Ric_{g_i(t)}|dt$ and $\int_0^{\min\{2,T\}} |\nabla^{\overline{g}} \Ric_{g_i(t)}|dt$ for each $i = 1, 2$. By pseudolocality and Shi's estimate, 
    \[
        |\Rm_{g_i(t)}| \le \dfrac{C}{\Gamma^2} \quad \text{and} \quad |\nabla^{\overline{g}} \Rm_{g_i}| \le \dfrac{C}{\Gamma^2 \sqrt{t}}
    \]
    for $i=1, 2$ and $t \in [0,T)$. Since $t^{-1/2}$ is integrable on $(0,2)$, we get 
    \[
        |\Rm_{g_0(t)}|_{g_0(t)} \le \dfrac{C}{\Gamma^2}
    \]
    for $t \in [0, \min\{2, T\})$ as well. Hence, the desired result follows.
\end{proof}

The cancellation of second metric derivatives in the curvature remainder is essential here. Pseudolocality controls the curvature of the two flows, and the integrable $t^{-1/2}$ first-derivative estimate controls their connection differences relative to the fixed background. Thus the remaining cutoff and first-derivative terms are $O(\Gamma^{-2})$. Combining this remainder with the two curvature bounds controls the curvature of the interpolated metric. Since $t^{-1}$ is not integrable, Shi's estimate for the second-derivative does not give the desired estimate.

Next, we construct a supersolution which will be a piece of barrier to improve the box condition \ref{box:B4}.

\begin{proposition}
    There exists a constant $C_v > 0$ with the following property. Let $R_0 \ge \underline{R_0}$ for a large enough $\underline{R_0} \gg 1$, $0 < p_0 \le \overline{p_0}$ for a small enough $0 < \overline{p_0} = \overline{p_0}(R_0) \ll 1$, $0 < \eta \le \overline{\eta}$ for a small enough $0 < \overline{\eta} \ll 1$, $\Gamma \ge \underline{\Gamma}$ for a large enough $\underline{\Gamma} = \underline{\Gamma}(R_0, \eta) \gg 1$, $0 < \epsilon \le \overline{\epsilon}$ for a small enough $0 < \overline{\epsilon} = \overline{\epsilon}(R_0, \eta, \Gamma) \ll 1$. Let $p \in \mathbb{R}^9$ with $|p| \le p_0$.
    Suppose that $(N, G_0)$ is a Riemannian manifold admitting an open set $U \subset N$ and a diffeomorphism $\iota \colon U \to M_{100\Gamma}$ satisfying \eqref{eqn:cl1} and \eqref{eqn:cl2}, and that for some $0 < T < \min\{2, T_G, T_p, T_\Phi\}$,
    
    \[
        (1-8\eta) \beta_p(t) \le A_p(t) \le (1+8\eta) \beta_p(t)
    \]
    on $\mathfrak{N}_{p, R_0, T}$. Then 
    \begin{equation}
        v := 2\eta^2 e^t \left(1 - \dfrac{C_v\chi}{1+f_{p,t}^N}\right)
        \label{eqn:bar-out}
    \end{equation}
    is a supersolution 
        \begin{equation}
        (\partial_t - (A_p)^{ab}\nabla_a^{\beta_p}\nabla_b^{\beta_p})v \ge C_R |\Rm_{\beta_p}|_{\beta_p} v + 2|S_p|_{\beta_p}v + 2|S_p|_{\beta_p} \sqrt{v}.
        \label{eqn:sups}
    \end{equation}
    on $\mathfrak{N}_{p, R_0, T}$.
    \label{prop:bar}
\end{proposition}
\begin{proof}
    Write $P_{A_p} := \partial_t - (A_p)^{ab} \nabla_a^{\beta_p} \nabla_b^{\beta_p}$. We show the inequality on the three regions $\{\chi \equiv 1\}$, $\{ 0 < \chi < 1\}$, and $\{ \chi \equiv 0 \}$ separately.

    \textbf{1. Inner Region: }$\{\chi \equiv 1\} \cap \mathfrak{N}_{p, R_0, T}$.

    On this region, 
    \[
        P_{A_p}v = 2\eta^2 e^t \left[ 1 - \dfrac{C_v}{1+f_{p,t}^N} + \dfrac{C_v}{(1+f_{p,t})^2}P_{A_p}f_{p,t}^N + \dfrac{2C_v}{(1+f_{p,t})^3}|df_{p,t}|_A^2 \right].
    \]

    First note that by taking $R_0$ sufficiently large compared to $C_v$, $\dfrac{C_v}{1+f_{p,t}^N} < \dfrac{C_v}{R_0}$ can be made arbitrarily small. Also, the fourth term is positive and therefore harmless.

    To estimate the third term, we compute the following. By shrinker identities,
    \begin{align*}
    \partial_t f_{p,t}^N(x) &= df|_{\Psi_{p,t}(\iota(x))}(\partial_t \Psi_{p,t}(\iota(x)))
    = \dfrac{1}{T_p-t} |df|_{\overline{g}}^2(\Psi_{p,t}(\iota(x))) = \dfrac{f_{p,t}^N}{T_p-t}-R_{\beta_p},\\
    \Delta_{\beta_p(t)}f_{p,t}^N(x) &= \dfrac{2}{T_p-t}-R_{\beta_p}.
    \end{align*}
    Thus $(\partial_t - \Delta_{\beta_p(t)}) f_{p,t}^N = \dfrac{f_{p,t}^N - 2}{T_p-t}$, which gives 
    \[
        P_{A_p}f_{p,t}^N = \dfrac{f_{p,t}^N - 2}{T_p-t} - ((A_p)^{ab} - (\beta_p)^{ab})\nabla_a^{\beta_p}\nabla_b^{\beta_p}f_{p,t} \ge \dfrac{f_{p,t}^N + 1}{2(T_p-t)}
    \]
    for sufficiently small $0 < \overline{\eta} \ll 1$ and sufficiently large $\underline{R_0} \gg 1$. 

    On this region, $\beta_p(t) = \overline{G}_p(t)$ and FIK shrinking Ricci soliton is asymptotically conical, so
    \[
        |\Rm_{\beta_p}|_{\beta_p} \le \dfrac{C_{AC}}{(T_p-t)(1+f_{p,t})}.
    \]
    for some $C_{AC} > 0$. Then for $C_v = 4 C_R \cdot C_{AC}$, 
    \[
        P_{A_p}v \ge C_R|\Rm_{\beta_p}|_{\beta_p}v.
    \]
    Since $S_p = 0$ on $\{\chi \equiv 1\}$, \eqref{eqn:sups} holds in this region.

    \textbf{2. Middle Region: }$\{0 < \chi < 1\} \cap \mathfrak{N}_{p, R_0, T}$.
    On this region, we show that derivatives of the cutoff contribute only small errors. The full expansion is
    \begin{align*}
        P_{A_p}v &= 2\eta^2 e^t \left[ 1 - \dfrac{C_v\chi}{1+f_{p,t}^N} + \dfrac{C_v\chi}{(1+f_{p,t}^N)^2}P_{A_p}f_{p,t}^N + \dfrac{2C_v\chi}{(1+f_{p,t}^N)^3}|df_{p,t}^N|_A^2 \right.\\
         &\quad\quad\quad\quad\quad\left. + \dfrac{C_v}{1+f_{p,t}^N}((A_p)^{ab}\nabla_a^{\beta_p}\nabla_b^{\beta_p}\chi) - \dfrac{2C_v}{(1+f_{p,t}^N)^2} \langle d\chi, df_{p,t}^N\rangle_A \right].
    \end{align*}
    For $x \in \{0 < \chi < 1\}$, $\Gamma \le \overline{b}(\iota(x)) \le 2\Gamma$. Therefore, $f_{p,t}$ is comparable to $\dfrac{\Gamma^2}{T_p-t}$. For $0 < \overline{\eta} \ll 1$ sufficiently small, we therefore have a constant $C_c > 0$ independent of $\Gamma$ and $T_p-t$ for which
    \begin{align*}
        &|P_{A_p}f_{p,t}^N| \le \dfrac{2(1+f_{p,t})}{T_p-t} \le \dfrac{C_c\Gamma^2}{(T_p-t)^2},\\
        &|df_{p,t}^N|_A^2 \le \dfrac{2(1+f_{p,t})}{T_p-t} \le \dfrac{C_c\Gamma^2}{(T_p-t)^2},\\
        &|(A_p)^{ab} \nabla_a^{\beta_p} \nabla_b^{\beta_p} \chi| \le \dfrac{C_c}{\Gamma^2},\\
        &|d\chi|_A \le \dfrac{C_c}{\Gamma}.
    \end{align*}

    Since $T_p < 2$, by taking $0< \overline{\eta} \ll 1$ sufficiently small and $\underline{R_0} \gg 1$ and $\underline{\Gamma} = \underline{\Gamma}(R_0) \gg 1$ sufficiently large, we have 
    \[
        P_{A_p}v \ge 2\eta^2 e^t \cdot \dfrac{1}{2}.
    \]
    Then by taking $\underline{\Gamma}$ large enough, since $M$ is asymptotically conical, $|\Rm_{\beta_p}|_{\beta_p} < \dfrac{1}{2C_R}$ on this region for $t \in [0,T]$. Then 
    \[
        C_R|\Rm(\beta_p)|_{\beta_p} v \le C_R \dfrac{1}{2C_R} 2\eta^2 e^t = \eta^2 e^t.
    \]
    Since $S_p = 0$ in this region, \eqref{eqn:sups} holds.

    \textbf{3. Outer Region: }$\{\chi \equiv 0\} \cap \mathfrak{N}_{p, R_0, T}$. 

    By \zcref[S]{lem:rde}, by taking $\underline{\Gamma} = \underline{\Gamma}(R_0, \eta) \gg 1$ sufficiently large, we have 
    \[
        |\Rm_{\beta_p}|_{\beta_p} \le \dfrac{1}{4C_R} \quad \text{and} \quad |S_p|_{\beta_p} \le \dfrac{\eta}{8}.
    \]
    on this region for $t \in [0,T]$. Then
    \begin{align*}
    (\partial_t - (A_p)^{ab}\nabla_a^{\beta_p}\nabla_b^{\beta_p})v &= 2\eta^2 e^t\\
    C_R|\Rm(\beta_p)|_{\beta_p} v + 2|S_p|_{\beta_p}v + 2|S_p|_{\beta_p} \sqrt{v} 
    &< \dfrac{1}{4}\cdot 2\eta^2 e^t + \dfrac{\eta}{4}\cdot 2\eta^2 e^t + \dfrac{\eta}{4}\cdot2\eta e^{t/2}\\
    &< 2\eta^2 e^t \left( \dfrac{1}{4} + \dfrac{1}{4} + \dfrac{1}{4}e^{-t/2}\right) < 2\eta^2 e^t. 
    \end{align*}
    Therefore, \eqref{eqn:sups} holds.
\end{proof}

\begin{proof}[Proof of {\zcref[S]{prop:aprhmhf}}]
    Compare $|H_p|_{\beta_p}^2$ with the supersolution \eqref{eqn:bar-out}. Take $\underline{R_0} > 2C_v$ sufficiently large so that on $\partial \mathfrak{N}_{p, R_0, T}$,
    \[
        2\eta^2 e^t \left(1 - \dfrac{C_v\chi}{1+f_{p,t}^N}\right) \ge 2\eta^2 \left(1 - \dfrac{C_v}{R_0} \right) \ge \eta^2. 
    \]
    Since $v$ is a supersolution by \zcref[S]{prop:bar}, the maximum principle yields
    \[
        |H_p|_{\beta_p}^2 \le v < 16\eta^2
    \]
    on $\mathfrak{N}_{p, R_0, T}$. The desired estimate \eqref{eqn:bh} follows.
\end{proof}

\section{\ref{box:B3} cannot cause an exit at later times} \label{sec:B3l}
In the next couple of sections, the analysis will primarily be done on a rescaled perturbation $h_p$. We therefore begin by deriving its evolution equation. For the computation, see \cite[Corollary 4.17]{Sto26}.

\begin{proposition}
    On $\Psi_{p,t_p(\tau)} (M_{4\Gamma})$, the tensor $h_p(\tau)$ given by \eqref{eqn:hp} satisfies
    \begin{equation} \label{eqn:hp-pe}
        \partial_\tau h_p = L_f h_p + Q(h_p),
    \end{equation}
    where
    \begin{align}
        Q(h)_{ij} &:= -\hat{h}^{ab} \nabla_a^{\overline{g}} \nabla_b^{\overline{g}} h_{ij} + g^{ab}g^{pq}(\nabla^{\overline{g}}h \ast \nabla^{\overline{g}} h)_{abpqij} \label{eqn:qh}\\
        &\quad + (\Rm_{\overline{g}})_{ja}{}^p{}_b(\overline{g}_{ip}\widetilde{h}^{ab} + \widehat{h}^{ab} h_{ip}) + (\Rm_{\overline{g}})_{ia}{}^p{}_b(\overline{g}_{jp}\widetilde{h}^{ab} + \hat{h}^{ab}h_{jp})\nonumber
    \end{align}
    with $g = \overline{g} + h$, $\hat{h}^{ab} := -(g^{ab} - \overline{g}^{ab})$, and $\widetilde{h}^{ab} := \hat{h}^{ab} - h^{ab}$.
\end{proposition}

To simplify the notation, we write 
\begin{equation} \label{eqn:qabc}
    Q(h) = Q_A(h) \ast \nabla^2 h + Q_B(h) \ast \nabla h \ast \nabla h + Q_C(h) \ast \Rm_{\overline{g}} \ast h \ast h
\end{equation}
where $|Q_A(h)| \le C_Q|h|$, $|Q_B(h)| \le C_Q$, and $|Q_C(h)| \le C_Q$ for a uniform $C_Q > 0$ as long as $|h| \le 8\eta \le 8\overline{\eta}$ for a small enough $0 < \overline{\eta} \ll 1$. We also note that $Q_A$, $Q_B$, and $Q_C$ are analytic in $h$.

We improve the local regularity bound in three steps. First, weighted $L^2$ control and interpolation give a small $C^0$ bound on the smaller region $M_{\theta' R(\tau) - 1}$. Next, we remove the drift by following $\nabla f$ and use a barrier over a fixed time interval of length 8 to extend this smallness to $M_{4R(\tau)}$. Interior parabolic estimates then recover the required $C^{12}$ bound on $M_{2R(\tau)}$.

\subsection{$C^0$ estimate on the initial time-slice}
We use the following interpolation inequality to bootstrap from $L^2$ estimate to $C^0$ estimate.
\begin{lemma}[Gagliardo--Nirenberg inequality]
    Let $s > 0$ and $B(x,1) \subset M^n$ be a unit geodesic ball. Then there exists $C_{GN} > 0$ depending only on the geometry of $B(x,1)$ such that 
    \[
        \norm{V}_{C^0(B(x, 1/2))} \le C_{GN}\norm{V}_{L^2(B(x,1))}^{\theta_{GN}} \norm{V}_{C^s(B(x,1))}^{1 - \theta_{GN}}
    \]
    for symmetric 2-tensors $V$. Here, $\theta_{GN} = \frac{2s}{2s+n}$.
\end{lemma}

\begin{proposition} \label{prop:c0-in}
    For any $K' > 0$, there exists $K > 0$ with the following property.
    Let $R_0 \ge \underline{R_0}$ for large enough $\underline{R_0} \gg 1$ and $0 < \eta \le \overline{\eta}$ for small enough $0 < \overline{\eta} \ll 1$. Let $h_p(\tau)$ satisfy \eqref{eqn:hp-pe} on $M_{4R(\tau)}$.
    Set 
    \[
        \mu = \dfrac{\eta}{K} e^{-\frac{\theta^2}{8}R_0^2}.
    \]
    Denote by 
    \[
        k_\tau^+ := \mathcal{P}_+ (w_{R(\tau)}h_p(\tau)) \quad \text{and} \quad k_\tau^- := \mathcal{P}_-(w_{R(\tau)}h_p(\tau)),
    \]
    where we suppress $p$ for simplicity.
    Given $\tau > 0$, if $\norm{k_\tau^+}_{L_f^2} \le \mu e^{-\gamma \tau}$, $\norm{k_\tau^-}_{L_f^2} \le \mu e^{-\gamma \tau}$, and $\norm{h_p(\tau)}_{C^2(M_{2R(\tau)})} \le \eta$, then 
    \begin{equation} \label{eqn:c0-k}
        \norm{k_\tau^+}_{C^0(M_{\theta' R(\tau) - 1})} \le \dfrac{\eta}{2K'} \quad \text{and} \quad \norm{k_\tau^-}_{C^0(M_{\theta'R(\tau) - 1})} \le \dfrac{\eta}{2K'}.
    \end{equation}
\end{proposition}

\begin{proof}
    For the unstable mode $k^+_\tau$, since the space of unstable modes is finite-dimensional, the $C^0$ norm is bounded by its $L_f^2$ norm up to a multiplicative constant. Therefore, the first inequality of \eqref{eqn:c0-k} is clear by increasing $K > 0$.

    For the stable mode $k^-_\tau$, by a direct estimate of the weight we have 
    \begin{align*}
        \norm{k^-_\tau}_{L^2(M_{\theta' R(\tau)})}
        \le e^{\frac{(\theta' R(\tau))^2}{8}} \norm{k^-_\tau}_{L_f^2}
        \le \dfrac{\eta}{K} e^{-\frac{\theta^2}{8}R_0^2 - \gamma \tau + \frac{(\theta')^2}{8}R(\tau)^2} \le \dfrac{\eta}{K}.
    \end{align*}
    Then by the Gagliardo-Nirenberg inequality with $s = 2$, 
    \begin{align*}
        \norm{k^-_\tau}_{C^0(M_{\theta' R(\tau) - 1})} \le C_{GN} \norm{k^-_\tau}_{L^2(M_{\theta' R(\tau)})}^{\theta_{GN}} \norm{k^-_\tau}_{C^2(M_{\theta' R(\tau)})}^{1-\theta_{GN}}
        \le C_{GN} \left( \dfrac{\eta}{K} \right)^{\theta_{GN}} \eta^{1-\theta_{GN}}
    \end{align*}
    where $\theta_{GN} = \frac{4}{4+4} = \frac{1}{2}$.
    Hence, by taking $K$ large enough, we get the second inequality of \eqref{eqn:c0-k}.
\end{proof}

\subsection{$C^0$ estimate on the final time-tube}

We rewrite the parabolic equation \eqref{eqn:hp-pe} with $h = h_p$ as follows:
\[
    \partial_\tau h = a^{ij} \nabla^{\overline{g}}_i \nabla^{\overline{g}}_jh - \nabla_{\nabla f} h + 2\Rm_{\overline{g}} (h) + Q_0(h)
\]
where $a^{ij} = (\overline{g}+h)^{ij}$ and
\begin{equation} \label{eqn:q0h}
    Q_0(h) = Q_B(h) \ast \nabla h \ast \nabla h + Q_C(h) \ast \Rm_{\overline{g}} \ast h \ast h
\end{equation}
for $Q_B$ and $Q_C$ as in \eqref{eqn:qabc}.

Define $w := |h|_{\overline{g}}^2$. Then 
\begin{align} \label{eqn:piw}
    (\partial_\tau - a^{ij}\nabla_i \nabla_j + \nabla_{\nabla f})w
    &= -2a^{ij} \langle \nabla_i h, \nabla_j h \rangle + 4 \langle \Rm_{\overline{g}}(h), h \rangle + 2\langle Q_0(h), h \rangle\\
    &\le -|\nabla h|_{\overline{g}}^2 + 4|\Rm_{\overline{g}}||h|_{\overline{g}}^2 + 2CC_Q\eta (|\nabla h|^2 + |h|^2) \nonumber\\
    &\le C_w|h|_{\overline{g}}^2 = C_w w \nonumber
\end{align}
for a fixed $C_w > 0$ provided $0 < \overline{\eta} \ll 1$ is sufficiently small and $0 < \eta \le \overline{\eta}$.

Following the trajectories of $\nabla f$ removes the first-order drift of \eqref{eqn:piw}. The trajectories of $\nabla f$ move outward through the soliton region. The interval of length 8 is chosen so that a smaller initial sublevel set covers the larger region needed at the final times. Its quantitative inclusion is established in \zcref[S]{lem:t7-est}.

Given $\tau_0 \ge 8$, we obtain a parabolic inequality on the time interval $[\tau_0 - 8, \tau_0]$. Let $\Upsilon_\sigma$ be the time-$\sigma$ flow of $\nabla f$. 
Define several quantities by pulling back along its trajectories as follows.
\begin{align}
    \widetilde{w}(y, \sigma) &:= w(\Upsilon_\sigma y, \tau_0 - 8 + \sigma),\\
    \widetilde{g} &:= e^{-\sigma}\Upsilon_\sigma^\ast \overline{g},\label{eqn:wtg}\\
    \widetilde{a} &:= e^{\sigma}\Upsilon_\sigma^\ast a. \label{eqn:wta}
\end{align}
Then the drift is removed and we obtain $(\partial_\sigma - e^{-\sigma}\widetilde{a}^{ij}\nabla_i^{\widetilde{g}} \nabla_j^{\widetilde{g}})\widetilde{w} \le C_w \widetilde{w}$.

Fix $S_0 > 0$ so that 
\begin{equation} \label{eqn:s0}
    0 < R_{\overline{g}} \le S_0 \quad \text{on } M.
\end{equation}

\begin{lemma} \label{lem:f-fl-est}
    For any $y \in M$ and $\sigma > 0$,
    \[
        e^\sigma (f(y) - S_0) \le f(\Upsilon_\sigma y) \le e^\sigma f(y).
    \]
\end{lemma}
\begin{proof}
    Let $u(y,\sigma) = f(\Upsilon_\sigma y)$. Then 
    \[
        \partial_\sigma u(y, \sigma) = \langle \nabla f, \nabla f \rangle = u - R_{\overline{g}}.
    \]
    Hence, $\partial_\sigma (e^{-\sigma} u(y,\sigma)) = -e^{-\sigma}R_{\overline{g}}$. The desired estimate follows from \eqref{eqn:s0}.
\end{proof}

\begin{lemma} \label{lem:t7-est}
    There exists $\underline{R_0} \gg 1$ such that if $R_0 \ge \underline{R_0}$, then 
    \[
        \Upsilon_7 M_{\frac{1}{2} \theta' R(\tau_0 - 8)} \supset M_{4R(\tau_0)}.
    \]
\end{lemma}
\begin{proof}
    By \zcref[S]{lem:f-fl-est}, it is enough to show that 
    \[
        e^7 \left( \dfrac{\left( \frac{1}{2}\theta' R(\tau_0-8) \right)^2}{4} - S_0 \right) \ge \dfrac{(4R(\tau_0))^2}{4}.
    \]
    After a simplification, a sufficient condition is 
    \[
        R_0^2 \ge 16e^7 S_0 + 64 \dfrac{8\gamma}{\theta^2} \cdot 8,
    \]
    which holds by taking $\underline{R_0}$ sufficiently large.
\end{proof}

We proceed to construct a barrier. Define $\widetilde{f}(y, \sigma) := e^{-\sigma}f(\Upsilon_\sigma y)$ and $\rho_0 = \dfrac{1}{2} \theta' R(\tau_0 - 8)$. Then define 
\[
    d(y,\sigma) := \rho_0 - 2\sqrt{\widetilde{f}(y,\sigma) + S_0 + 1}
\]
on $M_{\rho_0} \times [0,8]$.

\begin{lemma}
    $d(y,\sigma) < 0$ on $\partial M_{\rho_0} \times [0,8]$. \label{lem:d-neg}
\end{lemma}
\begin{proof}
    This follows from \zcref[S]{lem:f-fl-est}.
\end{proof}

\begin{lemma}
    There exist $\underline{R_0} \gg 1$ and $c_d > 0$ with the following property.  For $R_0 \ge \underline{R_0}$, $\sigma \in [7,8]$, and $x \in M_{4R(\tau_0)}$, we have $\Upsilon_{-\sigma}(x) \in M_{\rho_0}$ and 
    \[
        d(\Upsilon_{-\sigma}x, \sigma) \ge c_d\rho_0.
    \]
\end{lemma}
\begin{proof}
    The first claim follows from \zcref[S]{lem:t7-est}. The second claim is equivalent to 
    \[
        4(e^{-\sigma}f(x) + S_0 + 1) \le (1-c_d)^2 \rho_0^2,
    \]
    hence it is enough to show that 
    \[
        e^{-7}(4R(\tau_0))^2 + 4(1+S_0) \le (1-c_d)^2 \cdot \dfrac{1}{4} (\theta')^2 R(\tau_0-8)^2.
    \]
    Recall \eqref{eqn:theta}. By setting $c_d = \frac{\delta_\theta/2}{4+\delta_\theta} = 1/18$, a sufficient condition is 
    \[
        16e^{-7} \cdot \dfrac{8\gamma}{\theta^2} \cdot 8 + 4(S_0 + 1) \le 4\delta_\theta e^{-7}R_0^2,
    \]
    which holds by taking $\underline{R_0}$ large enough.
\end{proof}

\begin{lemma}
    $|\nabla^{\widetilde{g}}d|_{\widetilde{g}} \le 1$, $|\nabla^{\widetilde{g}, 2}d|_{\widetilde{g}} \le C_f + 1$, and $0 \le \partial_\sigma d \le S_0$, where $C_f$ is a constant so that $|\nabla^2 f|_{\overline{g}} \le C_f$. 
\end{lemma}
\begin{proof}
    The normalized shrinker identity gives 
    \[
        |d\widetilde{f}|_{\widetilde{g}}^2 = e^{-\sigma} (|\nabla f|_{\overline{g}}^2 \circ \Upsilon_\sigma) = \widetilde{f} - e^{-\sigma}(R_{\overline{g}}\circ \Upsilon_\sigma) \le \widetilde{f}.
    \]
    Moreover, 
    \[
        \Hess_{\widetilde{g}}\widetilde{f} = e^{-\sigma} \Upsilon_\sigma^\ast (\Hess_{\overline{g}}f),
    \]
    which yields $|\Hess_{\widetilde{g}} \widetilde{f}|_{\widetilde{g}} \le C_f$, because the tensor and metric scaling factors cancel.

    Let $F = \widetilde{f} + S_0 + 1$. Then 
    \begin{align*}
        |\nabla^{\widetilde{g}}d|_{\widetilde{g}}^2 \le \dfrac{\widetilde{f}}{F} \le 1 \quad \text{and} \quad
        |\Hess_{\widetilde{g}} d|_{\widetilde{g}} \le \dfrac{C_f}{\sqrt{F}} + \dfrac{\widetilde{f}}{2F^{3/2}} \le C_f + 1.
    \end{align*}
    Finally, since $\Upsilon_\sigma$ is the time-$\sigma$ flow of $\nabla f$, 
    \[
        \partial_\sigma \widetilde{f} = e^{-\sigma}(|\nabla f|_{\overline{g}}^2 - f) \circ \Upsilon_\sigma = -e^{-\sigma}(R_{\overline{g}}\circ \Upsilon_\sigma).
    \]
    Therefore, $\partial_\sigma d = \dfrac{e^{-\sigma}(R_{\overline{g}}\circ \Upsilon_\sigma)}{\sqrt{F}}$ and in particular $0 \le \partial_\sigma d \le S_0$. 
\end{proof}

The barrier we construct is 
\begin{equation}
    \phi(y, \sigma) := \exp[-\nu d(y,\sigma) + C_P(1+\nu^2)\sigma]
\label{eqn:bar-in}
\end{equation}
where $C_P := 2 + S_0 + 8(C_f + 1) + C_w$.

\begin{lemma}    \label{lem:bar-supsol}
    For any $\nu \ge 1$, 
    \[
        (\partial_\sigma -  e^{-\sigma}\widetilde{a}^{ij} \nabla^{\widetilde{g}}_i \nabla^{\widetilde{g}}_j - C_w)\phi \ge 0
    \]
    on $M_{\rho_0} \times [0, 8]$ and $\phi \ge 1$ on $\partial M_{\rho_0} \times [0,8]$.
\end{lemma}
\begin{proof}
    Write $\phi = e^{P}$ with $P = -\nu d + C_P(1+\nu^2)\sigma$. Then 
    \[
        \phi^{-1}(\partial_\sigma - e^{-\sigma}\widetilde{a}^{ij} \nabla^{\widetilde{g}}_i \nabla^{\widetilde{g}}_j - C_w)\phi 
        = \partial_\sigma P - e^{-\sigma}\widetilde{a}^{ij} \nabla^{\widetilde{g}}_i \nabla^{\widetilde{g}}_j P - e^{-\sigma}\widetilde{a}(\nabla^{\widetilde{g}}P, \nabla^{\widetilde{g}}P) - C_w.
    \]
    Since 
    \begin{align*}
        \partial_\sigma P &= -\nu \partial_\sigma d + C_P(1+\nu^2) \ge -\nu S_0 + C_P(1+\nu^2),\\
        |e^{-\sigma}\widetilde{a}^{ij} \nabla^{\widetilde{g}}_i \nabla^{\widetilde{g}}_j P| &\le 8\nu (C_f + 1),\\
        |e^{-\sigma}\widetilde{a}(\nabla^{\widetilde{g}}P, \nabla^{\widetilde{g}}P)| &\le 2\nu^2 |\nabla^{\widetilde{g}}d|_{\widetilde{g}}^2 \le 2\nu^2,
    \end{align*}
    we get 
    \[
        \phi^{-1}(\partial_\sigma - e^{-\sigma}\widetilde{a}^{ij} \nabla^{\widetilde{g}}_i \nabla^{\widetilde{g}}_j - C_w)\phi \ge C_P(1+\nu^2) - \nu S_0 - 8(C_f + 1)\nu - 2\nu^2 - C_w \ge 0.
    \]
\end{proof}

\begin{proposition} \label{prop:c0-im}
    For any $K'' > 0$, there exists $K' > 0$ with the following property.
    Let $R_0 \ge \underline{R_0}$ for large enough $\underline{R_0} = \underline{R_0}(K'') \gg 1$, $0 < \eta \le \overline{\eta}$ for small enough $0 < \overline{\eta} \ll 1$, $\Gamma \ge \underline{\Gamma}$ for large enough $\underline{\Gamma} = \underline{\Gamma}(R_0) \gg 1$, and $0 < \epsilon \le \overline{\epsilon}$ for small enough $0 < \overline{\epsilon} = \overline{\epsilon}(R_0, \eta, \Gamma) \ll 1$. Suppose that $(N, G_0)$ is a Riemannian manifold admitting an open set $U \subset N$ and a diffeomorphism $\iota \colon U \to M_{100\Gamma}$ satisfying \eqref{eqn:cl1} and \eqref{eqn:cl2} with $\epsilon$. Given $\tau_0 \ge 8$, suppose that 
    \[
        (1-8\eta)\beta_p(t) \le A_p(t) \le (1+8\eta)\beta_p(t)
    \]
    for $t \in [t_p(\tau_0-8), t_p(\tau_0)]$, where $t_p(\tau)$ is given by \eqref{eqn:tptau}. Denote by $k_\tau := w_{R(\tau)}h_p(\tau)$, where we suppress $p$ for simplicity. If
    \[
        \norm{k_\tau}_{C^0(M_{\theta'R(\tau) - 1})} \le \dfrac{\eta}{K'} \quad \text{for } \tau \in [\tau_0-8, \tau_0],
    \]
    then 
    \begin{equation}
        \norm{h_p(\tau)}_{C^0(M_{4R(\tau)})} \le \dfrac{\eta}{K''} \quad \text{for } \tau \in [\tau_0-1, \tau_0].
        \label{eqn:h-c0}
    \end{equation}
\end{proposition}

\begin{proof}
    We take the barrier $\phi$ as in \eqref{eqn:bar-in} with $\nu = \dfrac{c_d\rho_0}{16C_P}$. By taking $\underline{R_0}$ large enough, $\nu \ge 1$. Hence, $\phi$ is a supersolution by \zcref[S]{lem:bar-supsol}, and at $(y, \sigma)$ where $y = \Upsilon_{-\sigma}x$, $\sigma \in [7,8]$, and $x \in M_{4R(\tau_0)}$,
    \begin{align*}
        \phi & \le \exp\left[ -\dfrac{c_d\rho_0}{16C_P} \cdot c_d\rho_0 + C_P \left(1 + \left( \dfrac{c_d\rho_0}{16C_P}\right)^2\right) \cdot 8 \right]
        = \exp \left[ 8C_P - \dfrac{c_d^2\rho_0^2}{32C_P} \right] \le \exp \left[ 8C_P - \dfrac{c_d^2(\theta')^2}{128C_P}R_0^2 \right].
    \end{align*}

    Now we apply the maximum principle with the barrier 
    \[
        \widetilde{\phi} := e^{C_w\sigma}\left( \dfrac{\eta}{K'} \right)^2 + 256 \eta^2 \phi.
    \]
    Note that $\widetilde{w} \le \widetilde{\phi}$ on $\partial M_{\rho_0} \times [0,8]$ due to \zcref[S]{prop:c0-in} and \zcref[S]{lem:d-neg}. 
    More precisely, the initial face is bounded by the first term of $\widetilde{\phi}$ while the lateral face is bounded by the second term of $\widetilde{\phi}$.
    Therefore, $\widetilde{w} \le \widetilde{\phi}$ on $M_{\rho_0} \times [0,8]$. By taking $K'$ large enough so that 
    \[
        e^{C_w \cdot 8} \cdot \dfrac{1}{(K')^2} < \dfrac{1}{2(K'')^2}
    \]
    and taking $\underline{R_0}$ large enough so that 
    \[
        256\exp \left[ 8C_P - \dfrac{c_d^2(\theta')^2}{128C_P} R_0^2 \right] < \dfrac{1}{2(K'')^2},
    \]
    we get the desired estimate \eqref{eqn:h-c0}.
\end{proof}

\begin{proposition}\label{prop:c12-im}
    For any $K''' > 0$, there exists $K'' > 0$ with the following property. Let $R_0 \ge \underline{R_0}$ for large enough $\underline{R_0} \gg 1$, and $0 < \eta \le \overline{\eta}$ for small enough $0 < \overline{\eta} \ll 1$. Given $\tau_0 \ge 8$, suppose that for $\tau \in [\tau_0-8, \tau_0]$, $h(\tau)$ satisfies 
    \[
        \partial_\tau h = \Delta_{\overline{g}}h - \nabla^{\overline{g}}_{\nabla^{\overline{g}}f}h + 2\Rm_{\overline{g}}(h) + Q(h)
    \]
    where $Q(h)$ is defined by \eqref{eqn:qh} and 
    \[
        \sup_{\tau \in [\tau_0-1, \tau_0]} \norm{h(\tau)}_{C^0(M_{4R(\tau)})} < \dfrac{\eta}{K''}.
    \]
    Then 
    \[
        \norm{h(\tau_0)}_{C^{12}(M_{2R(\tau_0)})} < \dfrac{\eta}{K'''}.
    \]
\end{proposition}
\begin{proof}
    Define the undrifted tensor 
    \[
        \widetilde{h} = e^{-\sigma} \Upsilon_\sigma^\ast h(\tau_0-8+\sigma).
    \]
    The tensor satisfies the following evolution equation.
    \begin{align*}
        \partial_\sigma \widetilde{h}_{ij} 
        &= e^{-\sigma} \left[ \widetilde{a}^{ab} \nabla_a^{\widetilde{g}}\nabla_b^{\widetilde{g}}\widetilde{h}_{ij} + 2\Rm_{\widetilde{g}}(\widetilde{h})_{ij} -\widetilde{g}^{ab}(\Ric_{\widetilde{g}})_{ai} \widetilde{h}_{bj} - \widetilde{g}^{ab}(\Ric_{\widetilde{g}})_{aj}\widetilde{h}_{bi} + \widetilde{Q}_0(\widetilde{h}) \right],
    \end{align*}
    where $\widetilde{Q}_0$ is $Q_0$ as in \eqref{eqn:q0h} but with $\overline{g}$, $h$, $\nabla^{\overline{g}}$ replaced by $\widetilde{g}$, $\widetilde{h}$, $\nabla^{\widetilde{g}}$.
    For any $x \in M_{2R(\tau_0)}$, define $y := \Upsilon_{-8}x$.
    For the parabolic regularity, we define the following two domains:
    \begin{align*}
        \Omega &:= B_r^{\widetilde{g}(8)}(y) \times [8-r^2, 8],\\
        \Omega' &:= B_{2r}^{\widetilde{g}(8)}(y) \times [8-4r^2, 8].
    \end{align*}
    We claim that 
    \[
        B_{2r}^{\widetilde{g}(8)}(y) \subset \Upsilon_{-\sigma}M_{4R(\tau_0-8+\sigma)}
    \]
    for every $\sigma \in [8-4r^2, 8]$. Note that for every $z \in B_{2r}^{\widetilde{g}(8)}(y)$, we have $d_{\overline{g}}(\Upsilon_8 z, x) < 2e^4 r$.
    Therefore, 
    \[
        \overline{b}(\Upsilon_\sigma z) \le \overline{b}(\Upsilon_8 z) < 2R(\tau_0) + 2e^4 r.
    \]
    By increasing $\underline{R_0} \gg 1$, we have that $2R(\tau_0) + 2e^4r < 4R(\tau_0-1)$
    for any $\tau_0 \ge 8$, and the claim is proven. Therefore by the assumption, 
    \[
        \sup_{\Omega'} |\widetilde{h}|_{\widetilde{g}(\sigma)} \le \dfrac{\eta}{K''}.
    \]

    Since $(M, \overline{g})$ is asymptotically conical, the background metrics $\widetilde{g}(\sigma)$ and their derivatives for $0 \le \sigma \le 8$ are uniformly bounded. In particular, by writing 
    \[
        \widetilde{a} = \widetilde{g}^{-1} + [(\widetilde{g} + \widetilde{h})^{-1} - \widetilde{g}^{-1}],
    \]
    the linear operator has uniformly bounded $C^{10,\alpha}$ coefficients and the rest terms can be included in the nonlinearity so that \cite[Proposition C.3]{Sto26} applies to give
    \[
        \norm{h(\tau_0)}_{C^{12}(M_{2R(\tau_0)})} \le \dfrac{C\eta}{K''}.
    \]
    By choosing $K''$ large enough, we get the desired result.
\end{proof}

Combining the three propositions, we get a bootstrap from $L_f^2$ norm to $C^{12}$ norm.
\begin{corollary}[{\ref{box:B3}} cannot cause an exit at later times]
\label{cor:c12-im}
    There exists $K > 0$ with the following property.
    Let $R_0 \ge \underline{R_0}$ for large enough $\underline{R_0} = \underline{R_0}(K'') \gg 1$, $0 < \eta \le \overline{\eta}$ for small enough $0 < \overline{\eta} \ll 1$, $\Gamma \ge \underline{\Gamma}$ for large enough $\underline{\Gamma} = \underline{\Gamma}(R_0) \gg 1$, and $0 < \epsilon \le \overline{\epsilon}$ for small enough $0 < \overline{\epsilon} = \overline{\epsilon}(R_0, \eta, \Gamma) \ll 1$. 
    Define 
    \[
        \mu = \dfrac{\eta}{K}e^{-\theta^2 R_0^2/8}.
    \]
    Suppose that $(N, G_0)$ is a Riemannian manifold admitting an open set $U \subset N$ and a diffeomorphism $\iota \colon U \to M_{100\Gamma}$ satisfying \eqref{eqn:cl1} and \eqref{eqn:cl2} with $\epsilon$. Given $\tau_0 \ge 8$, suppose that 
    \[
        (1-8\eta)\beta_p(t) \le A_p(t) \le (1+8\eta)\beta_p(t)
    \]
    for $t \in [t_p(\tau_0-8), t_p(\tau_0)]$, where $t_p(\tau)$ is given by \eqref{eqn:tptau}. Denote by $k_\tau := w_{R(\tau)}h_p(\tau)$, where we suppress $p$ for simplicity. We also assume that $\norm{k_\tau^+}_{L_f^2} \le \mu e^{-\gamma \tau}$, $\norm{k_\tau^-}_{L_f^2} \le \mu e^{-\gamma \tau}$, and $\norm{h_p(\tau)}_{C^{12}(M_{2R(\tau)})} \le \eta$ for $\tau \in [0, \tau_0]$, then
    \[
        \norm{h_p(\tau_0)}_{C^{12}(M_{2R(\tau_0)})} \le \dfrac{\eta}{2}.
    \] 
\end{corollary}

\begin{proof}
    This is the combination of \zcref[S]{prop:c0-in}, \zcref[S]{prop:c0-im}, and \zcref[S]{prop:c12-im} with $K''' \ge 2$.
\end{proof}

\section{{\ref{box:B2}} cannot cause an exit} \label{sec:B2}
We now compare the two spectral components with the shrinking threshold $\mu e^{-\gamma \tau}$. After multiplying their squared norms by $e^{2\gamma \tau}$, the linearized evolution expands the nonnegative component and contracts the stable component. The nonlinear terms are small when $\eta$ is small, and the cutoff contributes a Gaussian tail error. The estimates below quantify these contributions.

\begin{proposition} \label{prop:pdeodi}
    There exist $C_{uQ} > 0$ and $C_{sQ} > 0$ with the following property.
    Let $K \ge \underline{K}$ for a large enough $\underline{K} \gg 1$, $0 < \theta_0 \le \overline{\theta_0}$ for a small enough $0 < \overline{\theta_0} \ll 1$, $0 < \eta \le \overline{\eta}$ for a small enough $0 < \overline{\eta} = \overline{\eta}(\theta_0) \ll 1$, and $R_0 \ge \underline{R_0}$ for a large enough $\underline{R_0} = \underline{R_0}(K) \gg 1$. Let $h_p(\tau)$ satisfy \eqref{eqn:hp-pe} on $M_{2R(\tau)}$. Denote by 
    \[
        k_\tau := w_{R(\tau)}h_p(\tau), \quad k_\tau^+ := \mathcal{P}_+k_\tau, \quad k_\tau^- := \mathcal{P}_-k_\tau.
    \]
    Set $\mu = \dfrac{\eta}{K}e^{-\frac{\theta^2}{8}R_0^2}$. If $\norm{k_\tau^+}_{L_f^2} \le \mu e^{-\gamma \tau}$, $\norm{k_\tau^-}_{L_f^2} \le \mu e^{-\gamma \tau}$, and $\norm{k_\tau}_{C^2} \le 2\eta$, then 
    \begin{align}
    \dfrac{1}{2}\partial_\tau (e^{2\gamma \tau}\norm{k_\tau^+}_{L_f^2}^2)
    &\ge e^{2\gamma \tau} \norm{k_\tau^+}_{L_f^2} \left( (\gamma - C_{uQ}\eta)\norm{k_\tau^+}_{L_f^2} - C_{uQ} \eta \norm{k_\tau^-}_{L_f^2} - \norm{E_\tau}_{L_f^2} \right),\label{eqn:pdeodi-u}\\
    \dfrac{1}{2}\partial_\tau (e^{2\gamma \tau} \norm{k_\tau^-}_{L_f^2}^2)
    &\le e^{2\gamma\tau} \norm{k_\tau^-}_{L_f^2} \left( C_{sQ}\eta \norm{k_\tau^+}_{L_f^2} + (C_{\gamma \theta_0} + C_{sQ}\eta)\norm{k_\tau^-}_{L_f^2} + \norm{E_\tau}_{L_f^2} \right), \label{eqn:pdeodi-s}
\end{align}
where $C_{\gamma\theta_0} := \gamma + ((1-\theta_0)\lambda_{10} + 2\theta_0 C_{\Rm}) < 0$. 
\end{proposition}

\begin{corollary}]\label{cor:pdeodi}
    With the same setup as in \zcref[S]{prop:pdeodi} and $\norm{h_p(\tau)}_{C^2(M_{2R(\tau)})} \le \eta$,
    \begin{itemize}
        \item (Unstable mode exits immediately)\\ If $\norm{k_\tau^+}_{L_f^2} = \mu e^{-\gamma \tau}$, $\norm{k_\tau^-}_{L_f^2} \le \mu e^{-\gamma \tau}$, and $\norm{k_\tau}_{C^2} \le 2\eta$, then $\partial_\tau (e^{2\gamma \tau}\norm{k_\tau^+}_{L_f^2}^2) > 0$.
        \item (Stable mode cannot exit)\\ If $\norm{k_\tau^+}_{L_f^2} \le \mu e^{-\gamma \tau}$, $\norm{k_\tau^-}_{L_f^2} = \mu e^{-\gamma \tau}$, and $\norm{k_\tau}_{C^2} \le 2\eta$, then $\partial_\tau (e^{2\gamma\tau}\norm{k_\tau^-}_{L_f^2}^2) < 0$.
    \end{itemize}
\end{corollary}

Multiplication by the moving cutoff gives
\begin{equation} \label{eqn:k-pe}
    \partial_\tau (w_{R(\tau)}h_p(\tau)) - L_f(w_{R(\tau)}h_p(\tau)) = Q(w_{R(\tau)}h_p(\tau)) + E_\tau,
\end{equation}
where the error term is 
\[
E_\tau = (\partial_\tau w_{R(\tau)}) h_p(\tau) - (\Delta_f w_{R(\tau)})h_p(\tau) - 2\nabla_{\nabla w_{R(\tau)}} h_p(\tau) + [w_{R(\tau)}Q(h_p(\tau)) - Q(w_{R(\tau)}h_p(\tau))].
\]
Write $h = h_p(\tau)$ for simplicity.
Suppose $\underline{R_0} \gg 1$ is large enough so that for $R_0 \ge \underline{R_0}$, 
\[
    |w_{R(\tau)}h_p(\tau)|_{C^2} \le 2\eta
\] 
provided $|h_p(\tau)|_{C^2(M_{2R(\tau)})} \le \eta$. Then there exists $C_k > 0$ such that the nonlinear commutator term of \eqref{eqn:k-pe} satisfies 
\[
    |w_{R(\tau)}Q(h) - Q(w_{R(\tau)}h)| \le C_k \eta^2.
\]

Note that $\supp E_\tau \subset \{ R(\tau) \le \overline{b} \le 2R(\tau) \}$.

We first show that the $L_f^2$ norm of the error $E_\tau$ is small.
\begin{lemma} \label{lem:ee}
    There exists $C_E > 0$ with the following property. Let $0 < \eta \le \overline{\eta}$ with a small enough $0 < \overline{\eta} \ll 1$ and $R_0 \ge \underline{R_0}$ for a large enough $\underline{R_0} \gg 1$. If $\norm{h_p(\tau)}_{C^2(M_{2R(\tau)})} \le \eta$, then
    \[
        \norm{E_\tau}_{L_f^2} \le C_E \eta \left(1 + \dfrac{|R'(\tau)|}{R(\tau)} \right) R(\tau)^2 e^{-R(\tau)^2/8}.
    \]
\end{lemma}
\begin{proof}
    By the construction, there exists a constant $C>0$ such that 
    \begin{align*}
        |\partial_\tau w_{R(\tau)}| &\le C \dfrac{|R'(\tau)|}{R(\tau)},\\
        |\Delta_f w_{R(\tau)}| &\le C,\\
        |\nabla_{\nabla w}h| &\le C\eta,\\
        |wQ(h) - Q(wh)| &\le C\eta.
    \end{align*}
    Hence, there exists such a constant $C_E > 0$.
\end{proof}

Write $k_\tau := w_{R(\tau)}h_p(\tau)$ and decompose it as 
\[
    k_\tau^+ := \mathcal{P}_+(w_{R(\tau)}h_p(\tau)) \quad \text{and} \quad k_\tau^- := \mathcal{P}_-(w_{R(\tau)}h_p(\tau)).
\]
Projecting the truncated equation and pairing with the corresponding component yields the following spectral inequalities.
\begin{align}
    \dfrac{1}{2}\partial_\tau \norm{k_\tau^+}_{L_f^2}^2 &= \langle \partial_\tau k_\tau^+, k_\tau^+ \rangle_{L_f^2} \label{eqn:si-u}\\
    &= \langle L_f k_\tau^+, k_\tau^+ \rangle_{L_f^2} + \langle \mathcal{P}_+ Q(k_\tau), k_\tau^+ \rangle_{L_f^2} + \langle E_\tau, k_\tau^+ \rangle_{L_f^2}\nonumber\\
    &\ge \langle \mathcal{P}_+ Q(k_\tau), k_\tau^+ \rangle_{L_f^2} + \langle E_\tau, k_\tau^+ \rangle_{L_f^2}\nonumber\\
    \dfrac{1}{2}\partial_\tau \norm{k_\tau^-}_{L_f^2}^2 &= \langle \partial_\tau k_\tau^-, k_\tau^- \rangle_{L_f^2} \label{eqn:si-s}\\
    &= \langle L_f k_\tau^-, k_\tau^- \rangle_{L_f^2} + \langle \mathcal{P}_- Q(k_\tau), k_\tau^- \rangle_{L_f^2} + \langle E_\tau, k_\tau^- \rangle_{L_f^2}\nonumber\\
    &\le \lambda_{10} \norm{k_\tau^-}_{L_f^2}^2 + \langle \mathcal{P}_- Q(k), k_\tau^- \rangle_{L_f^2} + \langle E_\tau, k_\tau^- \rangle_{L_f^2}. \nonumber
\end{align}

We henceforth choose $\underline{R_0} \gg 1$ sufficiently large that $\norm{h_p(\tau)}_{C^2(M_{2R(\tau)})} \le \eta$ implies $\norm{k_\tau}_{C^2} \le 2\eta$. When time is fixed, we suppress the subscript $\tau$.

\subsection{Unstable Mode}
\begin{lemma} \label{lem:pq-u}
    There exists $C_{uQ} > 0$ with the following property. Let $0 < \eta \le \overline{\eta}$ for a small enough $0 < \overline{\eta} \ll 1$, $R_0 \ge \underline{R_0}$ for a large enough $\underline{R_0} \gg 1$ and let $k$ be a compactly supported symmetric 2-tensor on $M$. Denote by 
    \[
        k^+ := \mathcal{P}_+k \quad \text{and} \quad k^- := \mathcal{P}_-k.
    \]
    If $\norm{k}_{C^2} \le 2\eta$, then
    $\norm{\mathcal{P}_+ Q(k)}_{L_f^2} \le C_{uQ}\eta (\norm{k^+}_{L_f^2} + \norm{k^-}_{L_f^2})$.
\end{lemma}
\begin{proof}
    In this proof, $\lesssim$ denotes an inequality up to a fixed multiplicative constant.
    Note that 
    \[
        \norm{\mathcal{P}_+ Q(k)}_{L_f^2} \lesssim \sum_{i=1}^9 |\langle Q(k), \overline{h}_i \rangle_{L_f^2}|.
    \]
    Hence, it is enough to bound each $|\langle Q(k), \overline{h}_i \rangle_{L_f^2}|$. Since $\norm{k}_{C^2} \le 2\eta$, we have 
    \begin{align*}
        |\langle Q_A(k) \ast \nabla^2 k, \overline{h}_i \rangle_{L_f^2}|
        \lesssim C_Q \norm{k}_{C^2} \norm{\overline{h}_i}_{C^0} \int |k|e^{-f}
        \lesssim C_Q \eta \norm{k}_{L_f^1}
        \lesssim C_Q \eta \norm{k}_{L_f^2},
    \end{align*}
    where the last inequality follows from the Cauchy-Schwarz inequality.

    For the next term, we have 
    \[
        |\langle Q_B(k) \ast \nabla k \ast \nabla k, \overline{h}_i \rangle_{L_f^2}| \lesssim C_Q \int_M |\nabla k|^2|\overline{h}_i|e^{-f} \lesssim C_Q \norm{\overline{h}_i}_{C^0} \int_M |\nabla k|^2 e^{-f}.
    \]
    By integration by parts, 
    \[
        \int_M |\nabla k|^2 e^{-f} = \int_M \langle \nabla k, \nabla k \rangle e^{-f} = -\int_M \langle k, \Delta_{f}k \rangle e^{-f} = -\int_M \langle k, \Delta k \rangle e^{-f} + \int_M \langle k, \nabla_{\nabla f}k \rangle e^{-f}.
    \]
    Each term is bounded as follows.
    \begin{align*}
        \int_M \langle k, \Delta k \rangle e^{-f} &\lesssim \norm{\Delta k}_{C^0} \int_M |k|e^{-f} \lesssim \norm{\Delta k}_{C^0} \left( \int_M |k|^2e^{-f} \right)^{1/2} \left( \int_M e^{-f} \right)^{1/2} \lesssim \eta \norm{k}_{L_f^2},\\
        \int_M \langle k, \nabla_{\nabla f}k \rangle e^{-f} &\lesssim \int |k||\nabla k| |\nabla f| e^{-f} \le \norm{\nabla k}_{C^0} \left( \int_M |k|^2 e^{-f} \right)^{1/2} \left( \int_M |\nabla f|^2e^{-f} \right)^{1/2} \lesssim \eta \norm{k}_{L_f^2}.
    \end{align*}
    Therefore,
    \[
        |\langle Q_B(k) \ast \nabla k \ast \nabla k, \overline{h}_i \rangle_{L_f^2}| \lesssim C_Q \eta \norm{k}_{L_f^2}.
    \]

    For the last term, 
    \[
        |\langle Q_C(k) \ast \Rm_{\overline{g}} \ast k \ast k, \overline{h}_i \rangle_{L_f^2}| \lesssim C_Q \cdot C_{\Rm} \cdot \eta\norm{k}_{L_f^2},
    \]
    where $C_{\Rm} > 0$ is a constant so that $|\Rm_{\overline{g}}|_{\overline{g}} \le C_{\Rm}$.
    Adding all three terms, we get the desired estimate.
\end{proof}

\subsection{Stable Mode}
\begin{lemma} \label{lem:pq-s}
    Given $C_+>0$ and $C_->0$, there exists $C_{sQ} > 0$ with the following property. Let $k$ be a compactly supported symmetric 2-tensor on $M$. Denote by 
    \[
        k^+ := \mathcal{P}_+ k \quad \text{and} \quad k^- := \mathcal{P}_- k.
    \]
    If $\norm{k^+}_{C^2} \le C_+\eta$, $\norm{k^-}_{C^2} \le C_-\eta$, and $\norm{k}_{C^2} \le 2\eta$, then
    \[
        |\langle Q(k), k^-\rangle| \le C_{sQ}\eta \left(\norm{k^+}_{L_f^2} \norm{k^-}_{L_f^2} + \norm{k^-}_{L_f^2}^2 + \norm{\nabla k^-}_{L_f^2}^2\right).
    \]
\end{lemma}
\begin{proof}
    As before, we bound each term. 
    \begin{align*}
        \norm{\langle Q_A(k) \ast \nabla^2 k, k^- \rangle_{L_f^2}}
        \lesssim C_Q \eta \int \left( |k^+||k^-| + |k^-|^2 \right)e^{-f}
        \lesssim C_Q \eta ( \norm{k^+}_{L_f^2} \norm{k^-}_{L_f^2} + \norm{k^-}_{L_f^2}^2 ).
    \end{align*}
    By Cauchy--Schwarz inequality,
    \begin{align*}
        |\langle Q_B(k) \ast \nabla k \ast \nabla k, k^- \rangle_{L_f^2}|
        &\lesssim \int |k^-||\nabla k^-|^2 e^{-f} + \int |k^-||\nabla k^+|^2 e^{-f}\\
        &\lesssim \eta \norm{\nabla k^-}_{L_f^2}^2 + \norm{\nabla k^+}_{C^0}^2 \norm{k^-}_{L_f^1}\\
        &\lesssim \eta \norm{\nabla k^-}_{L_f^2}^2 + \eta\norm{k^+}_{L_f^2} \norm{k^-}_{L_f^2},
    \end{align*}
    where the last inequality follows from the finite dimensionality of $V_+$ and Cauchy--Schwarz inequality applied to $k^-$.
    For the last term, 
    \begin{align*}
        |\langle Q_C(k) \ast \Rm_{\overline{g}} \ast k \ast k, k^- \rangle_{L_f^2}|
        &\lesssim \int |k^-|^3 e^{-f} + \int |k^+|^2|k^-|e^{-f}\\
        &\lesssim \eta \norm{k^-}_{L_f^2}^2 + \norm{k^+}_{C^0}^2\norm{k^-}_{L_f^1}\\
        &\lesssim \eta \norm{k^-}_{L_f^2}^2 + \eta\norm{k^+}_{L_f^2}\norm{k^-}_{L_f^2}.
    \end{align*}
    Adding all three terms, we get the desired estimate.
\end{proof}

\begin{lemma} \label{lem:ulfi}
    Given any $0 < \theta_0 < 1$ and $k^- \in V_- \cap H_f^1$,
    \[
        \langle k^-, L_f k^- \rangle_{L_f^2} \le ((1-\theta_0)\lambda_{10} + 2C_{\Rm}\theta_0) \norm{k^-}_{L_f^2}^2 - \theta_0 \norm{\nabla k^-}_{L_f^2}^2.
    \]
\end{lemma}
\begin{proof}
    By the spectral gap, 
    \begin{equation}
        \langle k^-, L_f k^- \rangle_{L_f^2} \le \lambda_{10} \norm{k^-}_{L_f^2}^2. \label{eqn:sp1}
    \end{equation}
    But also by integration by parts, 
    \begin{align}
        \langle k^-, L_f k^- \rangle_{L_f^2}
        &= -\norm{\nabla k^-}_{L_f^2}^2 + 2\langle \Rm_{\overline{g}}(k^-), k^- \rangle_{L_f^2}\nonumber\\
        &\le -\norm{\nabla k^-}_{L_f^2}^2 + 2C_{\Rm}\norm{k^-}_{L_f^2}^2. \label{eqn:sp2}
    \end{align}
    Then $(1-\theta_0) \times \text{\eqref{eqn:sp1}} + \theta_0 \times \text{\eqref{eqn:sp2}}$ gives the desired estimate.
\end{proof}

\begin{proof}[Proof of {\zcref[S]{prop:pdeodi}}]
    The first inequality \eqref{eqn:pdeodi-u} for unstable mode $k_\tau^+$ follows from \eqref{eqn:si-u} and \zcref[S]{lem:pq-u}.
    \begin{align*}
        \dfrac{1}{2}\partial_\tau (e^{2\gamma \tau}\norm{k_\tau^+}_{L_f^2}^2)
        &= \gamma e^{2\gamma \tau}\norm{k_\tau^+}_{L_f^2}^2 + e^{2\gamma \tau} \langle k_\tau^+, L_f k_\tau^+ \rangle_{L_f^2} + e^{2\gamma \tau} \langle k_\tau^+, Q(k)\rangle_{L_f^2} + e^{2\gamma \tau} \langle k_\tau^+, E_\tau \rangle_{L_f^2}\\
        &\ge e^{2\gamma \tau} \left( \gamma \norm{k_\tau^+}_{L_f^2}^2 - C_{uQ}\eta (\norm{k_\tau^+}_{L_f^2}^2 + \norm{k_\tau^+}_{L_f^2} \norm{k_\tau^-}_{L_f^2}) - \norm{E_\tau}_{L_f^2} \norm{k_\tau^+}_{L_f^2} \right)\\
        &= e^{2\gamma \tau} \norm{k_\tau^+}_{L_f^2} \left( (\gamma - C_{uQ}\eta)\norm{k_\tau^+}_{L_f^2} - C_{uQ} \eta \norm{k_\tau^-}_{L_f^2} - \norm{E_\tau}_{L_f^2} \right).
    \end{align*}
    For the stable mode $k_\tau^-$, we first note that since the space of unstable modes is finite-dimensional, $\norm{k^+}_{C^2}$ is bounded by its $L_f^2$ norm. Therefore, there exists $C_+ > 0$ such that $\norm{k^+}_{C^2} \le C_+ \eta$. Then by taking $C_- = C_+ + 2$, we also get $\norm{k_\tau^-}_{C^2} \le \norm{k}_{C^2} + \norm{k_\tau^+}_{C^2} \le C_-\eta$. Therefore, the assumption of \zcref[S]{lem:pq-s} holds, and we may take the corresponding $C_{sQ} > 0$. Then with \zcref[S]{lem:ulfi}, we obtain 
    \begin{align}
        \dfrac{1}{2}\partial_\tau (e^{2\gamma \tau} \norm{k^-}_{L_f^2}^2)
        &= \gamma e^{2\gamma \tau} \norm{k^-}_{L_f^2}^2 + e^{2\gamma \tau} \langle k^-, L_f k^- \rangle_{L_f^2} + e^{2\gamma \tau} \langle k^-, Q(k) \rangle_{L_f^2} + e^{2\gamma \tau} \langle k^-, E \rangle_{L_f^2}\label{eqn:si-s2}\\
        &\le e^{2\gamma \tau} \left( \gamma \norm{k^-}_{L_f^2}^2 + ((1-\theta_0)\lambda_{10} + 2\theta_0C_{\Rm})\norm{k^-}_{L_f^2}^2 - \theta_0 \norm{\nabla k^-}_{L_f^2}^2 \right.\nonumber\\
        & \left. \quad\quad\quad\quad+ C_{sQ}\eta ( \norm{k^+}_{L_f^2} \norm{k^-}_{L_f^2} + \norm{k^-}_{L_f^2}^2 + \norm{\nabla k^-}_{L_f^2}^2 ) + \norm{k^-}_{L_f^2}\norm{E}_{L_f^2} \right). \nonumber
    \end{align}
    First, take $\overline{\theta_0}$ small enough so that 
    \[
        C_{\gamma\theta_0} := \gamma + ((1-\overline{\theta_0})\lambda_{10} + 2\overline{\theta_0}C_{\Rm}) < 0.
    \]
    Such $\overline{\theta_0}$ exists because $\gamma$ was chosen to satisfy \eqref{eqn:gamma}. Then for $0 < \theta_0 \le \overline{\theta_0}$, we take $\overline{\eta} = \overline{\eta}(\theta_0) > 0$ small enough so that $C_{sQ} \overline{\eta} < \theta_0$. Then on the right hand side of the inequality \eqref{eqn:si-s2}, the coefficient of $\norm{\nabla k^-}_{L_f^2}^2$ has the favorable sign. The second inequality \eqref{eqn:pdeodi-s} now follows.
\end{proof}

\begin{proof}[Proof of {\zcref[S]{cor:pdeodi}}]
    Set $\nu \colon = \min \{\gamma, -C_{\gamma \theta_0} \} > 0$. We first take $\overline{\eta}$ sufficiently small so that 
    \[
        \max\{C_{uQ}, C_{sQ} \} \cdot \overline{\eta} \le \dfrac{\nu}{8}.
    \]
    Note that to prove the statement, by \zcref[S]{prop:pdeodi} and \zcref[S]{lem:ee}, it is enough to show 
    \[
        \dfrac{\nu}{2} \mu e^{-\gamma \tau} \ge 2C_E \eta R(\tau)^2 e^{-R(\tau)^2/8},
    \]
    where we have chosen $\underline{R_0}$ large enough so that $|R'(\tau)|/R(\tau) \le 1$.

    By \eqref{eqn:rtau}, it is enough to show 
    \[
        \dfrac{\nu}{2} \mu e^{-\gamma \tau} \ge 2C_E \eta \left( R_0^2 + \dfrac{8\gamma}{\theta^2}\tau \right) \exp \left[ -\dfrac{R_0^2}{8} - \dfrac{\gamma}{\theta^2}\tau \right].
    \]
    Now recall that $\mu = \dfrac{\eta}{K}e^{-\frac{\theta^2}{8}R_0^2}$. So after a simplification the inequality is equivalent to 
    \begin{equation} \label{eqn:corpdeodi1}
        \dfrac{\nu}{2K} \exp \left[ \dfrac{1}{8}(1-\theta^2)R_0^2 \right] \ge 2C_E \left( R_0^2 + \dfrac{8\gamma}{\theta^2}\tau \right) \exp \left[ -\left( \dfrac{1}{\theta^2} - 1 \right) \gamma \tau \right].
    \end{equation}
    Note that 
    \begin{align*}
        &\partial_\tau \left( (R_0^2 + \frac{8\gamma}{\theta^2}\tau) \exp \left[ -\left( \dfrac{1}{\theta^2} - 1\right)\gamma\tau \right] \right)\\
        &= \dfrac{8\gamma}{\theta^2} \exp \left[ -\left( \dfrac{1}{\theta^2} - 1 \right) \gamma \tau \right] - \left( \dfrac{1}{\theta^2} - 1 \right) \gamma \left( R_0^2 + \dfrac{8\gamma}{\theta^2}\tau \right) \exp \left[ -\left( \dfrac{1}{\theta^2}-1 \right)\gamma \tau \right]\\
        &= \left[ \dfrac{8\gamma}{\theta^2} - \left( \dfrac{1}{\theta^2} - 1 \right) \gamma \left( R_0^2 + \dfrac{8\gamma}{\theta^2}\tau \right) \right] \exp \left[ -\left( \dfrac{1}{\theta^2} - 1 \right) \gamma \tau \right] < 0
    \end{align*}
    provided $\underline{R_0} = \underline{R_0}(K) \gg 1$ is taken large enough. Therefore, it is enough to prove the inequality \eqref{eqn:corpdeodi1} at $\tau = 0$. Now the inequality reduces to 
    \[
        \dfrac{\nu}{2K} \exp \left[ \dfrac{1}{8} (1-\theta^2) R_0^2 \right] \ge 2C_E R_0^2.
    \]
    It is then clear that the inequality holds provided that $\underline{R_0} = \underline{R_0}(K) \gg 1$ is sufficiently large.
\end{proof}

\section{Parameter Space} \label{sec:parameters}
The topological argument requires the initial unstable component, rather than the parameter $p$ itself, to range over a ball. We therfore use the map $F(p) = \cP_+(w_{R_0}h_p(0))$ as a coordinate map. The next proposition proves that $F$ is invertible on the relevant parameter neighborhood and that the preimage of a small closed ball is a closed nine-dimensional ball. We then show that its initial stable component lies strictly inside the stable bound.

\begin{proposition}\label{prop:pb}
    Let $R_0 \ge \underline{R_0}$ for a large enough $\underline{R_0} \gg 1$, $\Gamma \ge \underline{\Gamma}$ for a large enough $\underline{\Gamma} = \underline{\Gamma}(R_0) \gg 1$, $0 < p_0 \le \overline{p_0}$ for a small enough $0 < \overline{p_0} = \overline{p_0}(R_0) \ll 1$, and $0 < \epsilon \le \overline{\epsilon}$ for a small enough $0 < \overline{\epsilon} = \overline{\epsilon}(R_0, \Gamma, p_0) \ll 1$. Then there exists $\overline{\mu} = \overline{\mu}(R_0, p_0) > 0$ with the following property.
    Suppose that $(N, G_0)$ is a Riemannian manifold admitting an open set $U \subset N$ and a diffeomorphism $\iota \colon U \to M_{100\Gamma}$ satisfying \eqref{eqn:cl1} and \eqref{eqn:cl2} with $0 < \epsilon \le \overline{\epsilon}$. Let $h_p(0)$ be given by \eqref{eqn:hp} with $\tau = 0$.
    For every $0 < \mu \le \overline{\mu}$, 
    \[
        \mathfrak{P}_\mu := \{p \in \overline{B}_{p_0}(0) \subset \mathbb{R}^9 : \norm{\mathcal{P}_+(w_{R_0}h_p(0))}_{L_f^2} \le \mu \}
    \]
    is diffeomorphic to the closed 9-ball, and its boundary in $\mathbb{R}^9$ is precisely 
    \[
        \partial \mathfrak{P}_\mu := \{p \in \overline{B}_{p_0}(0) \subset \mathbb{R}^9 : \norm{\mathcal{P}_+(w_{R_0}h_p(0))}_{L_f^2} = \mu \}.
    \]
\end{proposition}

\begin{proof}
    Write $(\iota^{-1})^\ast G_0 = \overline{g} + q$ on $M_{100\Gamma}$. By the condition, $\norm{q}_{C^{12}(M_{100\Gamma}); \overline{g}} \le \epsilon$. We may choose $\underline{\Gamma}$ large enough and $\overline{p_0}$ small enough so that for $\Gamma \ge \underline{\Gamma}$ and $p_0 \le \overline{p_0}$, 
    \[
        \Xi_{p}^{-1}(\supp w_{R_0}) \subset M_{4\Gamma}
    \]
    for every $p \in \mathbb{R}^9$ with $|p| \le p_0$. Note that since $\Phi_{p,0} = \mathrm{Id}_N$, 
    \[
        h_p(0) = T_p^{-1} (\Xi_{p}^{-1})^\ast (\overline{g} + q) - \overline{g}.
    \]
    Define $F \colon \overline{B}_{p_0}(0) \to V_+$ by $F(p) = \mathcal{P}_+(w_{R_0}h_p(0))$. 
    We first show that $F$ is a diffeomorphism from $B_{p_0}(0)$ onto its image if $\overline{p_0}$ is sufficiently small.

    If $q = 0$, \zcref[S]{prop:geolc} gives 
    \[
        T_p^{-1}(\Xi_{p}^{-1})^\ast \overline{g} - \overline{g} = \sum_{j=1}^9 p_j \overline{h}_j + O(|p|^2).
    \]
    Define operators $J_{R_0} \colon \mathbb{R}^9 \to V_+$ and $J \colon \mathbb{R}^9 \to V_+$ by 
    \begin{equation} \label{eqn:jp}
        Jp := \sum_{j=1}^9 p_j \overline{h}_j \quad \text{and} \quad
        J_{R_0}p := \mathcal{P}_+ (w_{R_0} Jp).
    \end{equation}
    We first note that $J_{R_0}$ is an isomorphism for sufficiently large $\underline{R_0} \gg 1$ and $R_0 \ge \underline{R_0}$ since 
    \[
        \norm{J_{R_0} - J}_{\mathbb{R}^9 \to V_+} \le \left( \sum_{j=1}^9 \norm{(1-w_{R_0})\overline{h}_j}_{L_f^2}^2 \right)^{1/2} \to 0
    \]
    as $R_0 \to \infty$.

    Since the flow $\Xi_{p}$ depends smoothly on $p$, we have the following Taylor expansion 
    \[
        F(p) = J_{R_0}p + \mathcal{R}(p)
    \]
    where for small enough $\overline{p_0}$ and $|p| \le p_0 \le \overline{p_0}$,
    \begin{align}
        |\mathcal{R}(p)| \le C_{\mathcal{R}}(|p|^2 + \epsilon),\label{eqn:crte1}\\
        \norm{D\mathcal{R}(p)}_{\mathbb{R}^9 \to (V_+, L_f^2)} \le C_{\mathcal{R}}(|p| + \epsilon), \label{eqn:crte2}
    \end{align}
    with $C_{\mathcal{R}} > 0$ independent of $q$.

    Decrease $\overline{p_0}$ and $\overline{\epsilon}$ so that 
    \[
        \overline{p_0} \le \dfrac{1}{8C_{\mathcal{R}}}\norm{J_{R_0}^{-1}}^{-1}_{V_+ \to \mathbb{R}^9} \quad \text{and} \quad \overline{\epsilon} \le \dfrac{p_0}{8C_{\mathcal{R}}}\norm{J_{R_0}^{-1}}^{-1}_{V_+ \to \mathbb{R}^9}.
    \]
    Define $\zeta(p) := J_{R_0}^{-1}F(p) - p$. By \eqref{eqn:crte1} and \eqref{eqn:crte2}, 
    \[
        \sup_{|p| \le p_0} \norm{D\zeta(p)} \le \dfrac{1}{4} \quad \text{and} \quad |\zeta(0)| \le \dfrac{1}{4}p_0.
    \]
    In particular, 
    \begin{equation} \label{eqn:jrf-lip}
        |J_{R_0}^{-1}(F(p) - F(p'))| \ge |p-p'| - |\zeta(p) - \zeta(p')| \ge \dfrac{3}{4}|p-p'|.
    \end{equation}
    This proves the injectivity of $J_{R_0}^{-1}F$ on $B_{p_0}(0)$. Also, $J_{R_0}^{-1}DF = I + D\zeta$ is invertible on $B_{p_0}(0)$. Therefore, by the inverse function theorem $F$ is a diffeomorphism from $B_{p_0}(0)$ onto its image.

    Let $\overline{\mu} = p_0/(4 \norm{J_{R_0}^{-1}}_{V_+ \to \mathbb{R}^9})$. For any $h \in V_+$ with $\norm{h}_{L_f^2} \le \mu$, consider a map 
    \[
        \mathcal{T}_h(p) = J_{R_0}^{-1}h - \zeta(p).
    \]
    It is Lipschitz with Lipschitz constant $1/4$, and therefore 
    \[
        |\mathcal{T}_h(p)| \le \norm{J_{R_0}^{-1}}_{V_+ \to \mathbb{R}^9}\norm{h}_{L_f^2} + |\zeta(0)| + \dfrac{1}{4}|p| \le \dfrac{3}{4}p_0.
    \]
    Therefore, $\mathcal{T}_h$ is a contraction of $\overline{B}_{p_0}(0)$, and its unique fixed point solves $F(p) = h$. Therefore, 
    \[
        \overline{B}_{\overline{\mu}}^{V_+}(0) \subset F(B_{p_0}(0)).
    \]

    Consequently, for any $0 < \mu \le \overline{\mu}$, the map $\overline{B}_1^{\mathbb{R}^9}(0) \to \mathfrak{P}_\mu$ given by 
    \[
        (s_1, \ldots, s_9) \mapsto F^{-1} \left( \mu \sum_{j=1}^9 s_j \dfrac{\overline{h}_j}{\norm{\overline{h}_j}_{L_f^2}} \right)
    \]
    is a diffeomorphism and its inverse sends $p$ to $(s_1, \ldots, s_9)$ with 
    \[
        s_j = \mu^{-1} \langle F(p), \dfrac{\overline{h}_j}{\norm{\overline{h}_j}_{L_f^2}} \rangle_{L_f^2}.
    \]
    The desired result follows.
\end{proof}

\begin{proposition}\label{prop:pbs}
    Let $K \ge \underline{K}$ for large enough $\underline{K} \gg 1$, $R_0 \ge \underline{R_0}$ for large enough $\underline{R_0} = \underline{R_0}(K) \gg 1$, $0 < \mu \le \overline{\mu}$ for small enough $0 < \overline{\mu} \ll 1$, $\Gamma \ge \underline{\Gamma}$ for large enough $\underline{\Gamma} = \underline{\Gamma}(R_0) \gg 1$, and $0 < \epsilon \le \overline{\epsilon}$ for small enough $0 < \overline{\epsilon} = \overline{\epsilon}(\mu) \ll 1$. 
    Suppose that $(N, G_0)$ is a closed Riemannian manifold admitting an open set $U \subset N$ and a diffeomorphism $\iota \colon U \to M_{100\Gamma}$ satisfying \eqref{eqn:cl1} and \eqref{eqn:cl2} with $\epsilon$. Let $h_p(0)$ be given by \eqref{eqn:hp} with $\tau = 0$. Then
    \[
        \sup_{p \in \mathfrak{P}_\mu} \norm{\mathcal{P}_-(w_{R_0}h_p(0))}_{L_f^2} \le \dfrac{\mu}{2}.
    \]
\end{proposition}
\begin{proof}
    Recall that $F(p) := \mathcal{P}_+(w_{R_0}h_p(0))$ and that the maps $J$ and $J_{R_0}$ are given by \eqref{eqn:jp}. By \eqref{eqn:jrf-lip} with $p' = 0$, we have
    \begin{equation} \label{eqn:initstab1}
        |p| \le \dfrac{4}{3} \norm{J_{R_0}^{-1}} \left( \norm{F(p)}_{L_f^2} + \norm{F(0)}_{L_f^2} \right) \le C_p(\mu + \epsilon)
    \end{equation}
    where $C_p = C_p(\underline{R_0}) > 0$. This is because for all sufficiently large $R_0$, the norms  $\norm{J_{R_0}^{-1}}$ are uniformly bounded and we have that $\norm{F(p)}_{L_f^2} \le C\mu$ from the definition of $\mathfrak{P}_\mu$ and $\norm{F(0)}_{L_f^2} \le C\epsilon$ from the $\epsilon$-closeness of $G_0$ to $\overline{g}$. 

    By \zcref[S]{prop:geolc}, 
    \begin{equation} \label{eqn:initstab2}
        \norm{\mathcal{P}_-(w_{R_0}\mathcal{R}(p))}_{L_f^2} \le \norm{w_{R_0}\mathcal{R}(p)}_{L_f^2} \le C_{w\mathcal{R}} (|p|^2 + \epsilon)
    \end{equation}
    for a constant $C_{w\mathcal{R}} > 0$. We also have 
    \begin{equation} \label{eqn:initstab3}
        \norm{\mathcal{P}_-(w_{R_0}J)}_{\mathbb{R}^9 \to V_-} =: \alpha_{R_0} \to 0
    \end{equation}
    as $\norm{\mathcal{P}_-(w_{R_0}J)} = \norm{\mathcal{P}_-((1-w_{R_0})J)}$, which is bounded by $(R_0+1)e^{-R_0^2/8}$ up to a multiplicative factor. Combining \eqref{eqn:initstab1}, \eqref{eqn:initstab2}, and \eqref{eqn:initstab3}, we get 
    \[
        \sup_{p \in \mathfrak{P}_\mu} \norm{\mathcal{P}_-(w_{R_0}h_p(0))}_{L_f^2} \le C_s[\alpha_{R_0}(\mu+\epsilon) + (\mu + \epsilon)^2 + \epsilon]
    \]
    for a uniform $C_s = C_s(\underline{R_0}) > 0$.
    We may increase $\underline{R_0}$ so that
    \[
        \alpha_{R_0} \le \dfrac{1}{12C_s}, \quad \overline{\mu} \le \dfrac{1}{24C_s}, \quad \overline{\epsilon} = \min\left(1, \dfrac{1}{6C_s} \right) \mu.
    \]
    Then 
    \[
        \sup_{p \in \mathfrak{P}_\mu} \norm{\mathcal{P}_-(w_{R_0}h_p(0))}_{L_f^2} \le C_s \left[ \dfrac{1}{12C_s} \cdot 2\mu + \dfrac{1}{6C_s} \cdot \mu + \dfrac{1}{6C_s} \cdot \mu \right] = \dfrac{\mu}{2}. \qedhere
    \]
\end{proof}

\section{{\ref{box:B3}} cannot cause an exit at earlier times} \label{sec:B3e}
The later-time argument uses an interval of length 8 and therefore does not apply at the beginning of the flow. For $0 \le \tau \le 8$, we instead use the small initial comparison error. A barrier based on the larger radius $\Gamma$ keeps the lateral boundary contribution small on the soliton region, and initial-time parabolic estimates give the required derivative bounds.
\subsection{$C^0$ estimate on the initial time tube}
\begin{proposition}[{\ref{box:B3}} cannot cause an exit at earlier times] \label{prop:c0-im-init}
    Let $K \ge \underline{K}$ for large enough $\underline{K} \gg 1$, $R_0 \ge \underline{R_0}$ for large enough $\underline{R_0} = \underline{R_0}(K) \gg 1$, $0 < \eta \le \overline{\eta}$ for small enough $0 < \overline{\eta} = \overline{\eta}(K, R_0) \ll 1$ and $\Gamma \ge \underline{\Gamma}$ for $\underline{\Gamma} = \underline{\Gamma}(K, R_0) \gg 1$. Define 
    \[
        \mu = \dfrac{\eta}{K}e^{-\theta^2 R_0^2/8}.
    \]
    Let $0 < \epsilon \le \overline{\epsilon}$ for small enough $0 < \overline{\epsilon} = \overline{\epsilon}(K, R_0, \eta) \ll 1$.
    Suppose that $(N, G_0)$ is a closed Riemannian manifold admitting an open set $U \subset N$ and a diffeomorphism $\iota \colon U \to M_{100\Gamma}$ satisfying \eqref{eqn:cl1} and \eqref{eqn:cl2}. Given $p \in \mathfrak{P}_\mu$ and $0 \le \tau_0 \le 8$, if 
    \[
        (1-8\eta)\beta_p(t) \le A_p(t) \le (1+8\eta)\beta_p(t)
    \]
    for $t \in [0, t_p(\tau_0)]$, then 
    \[
        \sup_{\tau \in [0, \tau_0]} \norm{h_p(\tau)}_{C^{12}(M_{2R(\tau)})} \le \dfrac{\eta}{2}.
    \]
\end{proposition}
\begin{proof}
    The proof is very similar to the arguments in Section \ref{sec:B3l}. Define
    \[
        d_\Gamma(y, \sigma) := \Gamma - 2\sqrt{e^{-\sigma}f(\Upsilon_\sigma y) + S_0 + 1}
    \]
    and 
    \[
        \phi_\Gamma(y, \sigma) := \exp \left[ -\nu d_\Gamma + C_P(1+\nu^2)\sigma \right]
    \]
    on $M_\Gamma \times [0,\tau_0]$. As in the proof of \zcref[S]{lem:d-neg}, $d_\Gamma < 0$ on $\partial M_\Gamma \times [0, \tau_0]$. Also similar to \zcref[S]{lem:bar-supsol}, 
    \[
        (\partial_\sigma - e^{-\sigma} \widetilde{a}^{ij}\nabla_i^{\widetilde{g}}\nabla_j^{\widetilde{g}}-C_w)\phi_\Gamma \ge 0.
    \]
    Since $\phi_\Gamma \ge 1$ on $\partial M_\Gamma$, by the maximum principle,
    \[
        |\widetilde{h}|_{\widetilde{g}}^2 \le C_h^2(\mu + \epsilon)^2e^{C_w\sigma} + 256\eta^2\phi_\Gamma,
    \]
    for some $C_h > 0$, where the first term on the right hand side comes from the fact that 
    \begin{equation} \label{eqn:initsm}
        \norm{h_p(0)}_{C^{13}} \le C_h(\mu + \epsilon)
    \end{equation}
    given by \zcref[S]{prop:geolc} and \eqref{eqn:initstab1}.

    For $x \in M_{4R(8)}$ and $y = \Upsilon_{-\sigma}x$ where $0 \le \sigma \le 8$, 
    \[
        d_\Gamma(y,\sigma) = \Gamma - 2\sqrt{e^{-\sigma}f(x) + S_0 + 1} \ge \dfrac{\Gamma}{2}
    \]
    by taking $\Gamma^2 \ge 16(4R(8)^2 + S_0 + 1)$. Choose $\nu = \dfrac{\Gamma}{32C_P} \ge 1$. Then 
    \[
        \phi_\Gamma(y,\sigma) \le \exp \left( 8C_P - \dfrac{\Gamma^2}{128C_P} \right).
    \]
    Therefore, 
    \[
        \sup_{0 \le \tau \le \tau_0} \norm{h_p(\tau)}_{C^0(M_{4R(8)})}^2 \le C_h^2(\mu+\epsilon)^2 e^{8C_w} + 256\eta^2 \exp \left( 8C_P - \dfrac{\Gamma^2}{128C_P} \right).
    \]

    Recall that we choose $\mu = \dfrac{\eta}{K}e^{-\theta^2R_0^2/8}$. We also take $\overline{\epsilon}$ sufficiently small so that $\overline{\epsilon} \le \mu$. Then $0 < \dfrac{\mu}{\eta} \ll 1$ can be made arbitrarily small so that for $0 \le \tau \le 8$,
    \[
        \norm{h_p(\tau)}_{C^0(M_{4R(8)})} \le \dfrac{\eta}{K''}
    \]
    where $K''$ can be chosen arbitrarily large by taking $K$ arbitrarily large. Then for a fixed small enough $r>0$, when $\tau \ge 4r^2$, \cite[Proposition C.3]{Sto26}, and when $0 \le \tau \le 4r^2$, \cite[Proposition C.4]{Sto26} with \eqref{eqn:initsm} gives 
    \[
        \norm{h_p(\tau)}_{C^{12}(M_{2R(\tau)})} \le \dfrac{\eta}{2}
    \]
    by choosing $K'' = K''(r) > 0$ large enough.
\end{proof}

\section{Proof of the Main Theorem} \label{sec:proofmain}
We now prove the main result.
\subsection{Ricci Flow cannot cause an exit}
\begin{proposition}[Ricci flow cannot cause an exit] \label{prop:rfce}
    Let $0 < p_0 \le \overline{p_0}$ for small enough $0 < \overline{p_0} \ll 1$, $0 < \eta \le \overline{\eta}$ for small enough $0 < \overline{\eta} \ll 1$ and $0 < \epsilon \le \overline{\epsilon}$ for small enough $0 < \overline{\epsilon} \ll 1$. 
    Suppose that $(N, G_0)$ is a closed Riemannian manifold admitting an open set $U \subset N$ and a diffeomorphism $\iota \colon U \to M_{100\Gamma}$ satisfying \eqref{eqn:cl1} and \eqref{eqn:cl2}. Let $(N, G(t))_{t \in [0, T_G)}$ be its maximal Ricci flow. Given $p \in \mathbb{R}^9$ with $|p| \le p_0$, define 
    \[
        \tau^\ast(p) = \sup \left\{ \sigma \ge 0 : G(t_p(\tau)) \text{ and } \Phi_{p, t_p(\tau)} \text{ exist and the box conditions } \mathfrak{B} \text{ hold for } 0 \le \tau \le \sigma \right\}.
    \]
    If $\tau^\ast(p) < \infty$, then $t^\ast(p) := t_p(\tau^\ast(p)) < T_G$.
\end{proposition}
\begin{proof}
    For simplicity, we write $A = A_p$, $\beta = \beta_p$, $H = H_p$, $S = S_p$, and $\nabla^\beta = \nabla^{\beta_p}$.
    First, note that from the pseudolocality, $T_G > 0$. Hence, we may assume $\tau^\ast(p) > 0$, which gives $0 < t^\ast(p) < T_p$ since $\tau^\ast(p) < \infty$.
    Note that by the construction, $H = A-\beta$ satisfies 
    \begin{equation} \label{eqn:rfce1}
        \partial_t H = A^{ab} \nabla_a^\beta \nabla_b^\beta H + \Rm_\beta \ast H + A^{-1} \ast A^{-1} \ast \nabla^\beta H \ast \nabla^\beta H - S.
    \end{equation}
    We define two annuli
    \[
        \cA := \iota^{-1}(\{4\Gamma < \overline{b} < 8\Gamma\}) \quad \text{and} \quad \cA' := \iota^{-1}(\{3\Gamma < \overline{b} < 9\Gamma\}).
    \]
    By pseudolocality on $\cA'$, the curvature $\Rm_{\iota^\ast \overline{G}_p(t)}$ and $\Rm_{G(t)}$ are uniformly bounded on $\cA' \times [0, t^\ast(p))$. Then by Shi's estimate on $[t^\ast(p)/4, t^\ast(p))$, 
    \[
        \norm{\iota^\ast \overline{G}_p(t)}_{C^4(\cA')} + \norm{G(t)}_{C^4(\cA')} + \norm{\overline{G}_p(t)^{-1}}_{C^0(\cA')} + \norm{G(t)^{-1}}_{C^0(\cA')}
    \]
    is uniformly bounded for $t \in [t^\ast(p)/4, t^\ast(p))$. This in particular implies that $\norm{\Ric_{\beta_p(t)}}_{C^1(\cA')}$ is uniformly bounded on $[t^\ast(p)/4, t^\ast(p))$, and
    therefore, $|S_p|_{C^1(\cA')}$ is also uniformly bounded on the same time interval.

    Noting that $A^{-1} = (\beta + H)^{-1}$, \cite[Proposition C.3]{Sto26} is applicable to \eqref{eqn:rfce1}. By taking $\overline{\eta}$ and $r>0$ sufficiently small, we have that $\norm{H}_{C^0(N)} + r^2\norm{S_p}_{C^{0,\alpha}(N)}$ is below the threshold for all $t \in [t^\ast(p)/4, t^\ast(p))$, then \cite[Proposition C.3]{Sto26} gives that $\norm{H}_{C^{2,\alpha}(N)}$ is uniformly bounded for $t \in [t^\ast(p)/2, t^\ast(p))$.

    Since 
    \[
        |\Rm_A|_A \le C(|\Rm_\beta|_{\beta} + |\nabla^{\beta,2}H|_{\beta} + |\nabla^\beta H|_{\beta}^2),
    \]
    we get that $|\Rm_A|_A$ is uniformly bounded for $t \in [t^\ast(p)/2, t^\ast(p))$. Hence, $T_G > t^\ast(p)$.
\end{proof}

\subsection{Harmonic map heat flow cannot cause an exit}
\begin{proposition}[Harmonic map heat flow cannot cause an exit] \label{prop:hmhfce}
    Let $0 < p_0 \le \overline{p_0}$ for small enough $0 < \overline{p_0} \ll 1$, $0 < \eta \le \overline{\eta}$ for small enough $0 < \overline{\eta} \ll 1$ and $0 < \epsilon \le \overline{\epsilon}$ for small enough $0 < \overline{\epsilon} \ll 1$. Suppose that $(N, G_0)$ is a closed Riemannian manifold admitting an open set $U \subset N$ and a diffeomorphism $\iota \colon U \to M_{100\Gamma}$ satisfying \eqref{eqn:cl1} and \eqref{eqn:cl2}. Let $(N, G(t))_{t \in [0, T_G)}$ be its maximal Ricci flow. Given $p \in \mathbb{R}^9$ with $|p| \le p_0$, let $\Phi_p \colon N \times [0, T_\Phi) \to N$ be its maximal harmonic map heat flow to the target $\beta_p$ whose each time slice is a diffeomorphism. Define
    \[
        \tau^\ast(p) = \sup \left\{ \sigma \ge 0 : G(t_p(\tau)) \text{ and } \Phi_{p, t_p(\tau)} \text{ exist and the box conditions } \mathfrak{B} \text{ hold for } 0 \le \tau \le \sigma \right\}.
    \]
    If $\tau^\ast(p) < \infty$, then $t^\ast(p) := t_p(\tau^\ast(p)) < T_\Phi$.
\end{proposition}
\begin{proof}
    We take $\overline{\eta}$ so that $16\overline{\eta} < \eta_0 < \eta_1 < \overline{\eta}_4$ for some $\eta_0$ and $\eta_1$, where $\overline{\eta}_4$ is taken from \cite[Proposition A.9]{BK22} with $n=4$. 

    By \eqref{eqn:cl1} and \eqref{eqn:cl2}, $T_\Phi > 0$. Hence, we may assume $\tau^\ast(p) > 0$, which gives $0 < t^\ast(p) < T_p$ since $\tau^\ast(p) < \infty$. As in the proof of \zcref[S]{prop:rfce}, $|\Rm_A|_A$ is uniformly bounded for $t \in [t^\ast(p)/2, t^\ast(p)]$. Also, with the same proof but with Shi's estimate for higher derivatives, we have that $\norm{S_p}_{C^{12}}$ is uniformly bounded for $t \in [t^\ast(p)/2, t^\ast(p))$. Therefore, by \cite[Proposition A.9]{BK22}, $T_\Phi > t^\ast(p)$.
\end{proof}

\subsection{Not every parameter can exit through {\ref{box:B1}}}

\begin{proof}[Proof of {\zcref[S]{thm:A} and \zcref[S]{thm:mainthm}}]
    Recall that \zcref[S]{thm:A} is equivalent to \zcref[S]{thm:mainthm} by fixing the scale. Therefore, we only prove \zcref[S]{thm:mainthm}.
    We choose the constants as in \zcref[S]{app:consdep}. 

    For $p \in \mathfrak{P}_\mu$, we first show that the box conditions $\mathfrak{B}$ hold at $\tau = 0$.
    \begin{itemize}
        \item \ref{box:B1} holds by the construction of $\mathfrak{P}_\mu$,
        \item \ref{box:B2} holds by \zcref[S]{prop:pbs},
        \item \ref{box:B3} holds by the choice of $K$, $\mu$, $\epsilon$, and \zcref[S]{prop:geolc}, and
        \item \ref{box:B4} holds because \zcref[S]{prop:geolc} gives $A_p$ and $\beta_p$ are close on $M_{8\Gamma}$ and $A_p$ and $\beta_p$ are identical on $M \setminus M_{8\Gamma}$.
    \end{itemize}

    Now by \zcref[S]{prop:rfce} and \zcref[S]{prop:hmhfce}, we have short time existence of Ricci flow and harmonic map heat flow. In particular, they cannot be the exit conditions. 

    Then by \zcref[S]{prop:aprhmhf}, \zcref[S]{cor:c12-im}, \zcref[S]{cor:pdeodi}, and \zcref[S]{prop:c0-im-init}, if $\tau^\ast(p) < \infty$, the only possible way to exit is via \ref{box:B1} condition of the box $\mathfrak{B}$. 
    For the sake of contradiction, suppose that $\tau^\ast(p) < \infty$ for all $p \in \mathfrak{P}_\mu$. We first show that the exit time $\tau^\ast(p)$ is continuous. Define 
    \[
        Q(p,\tau) := \mu^{-2}e^{2\gamma\tau}\norm{\cP_+(w_{R(\tau)}h_p(\tau))}_{L_f^2}^2.
    \]
    Note that \zcref[S]{prop:rfce} and \zcref[S]{prop:hmhfce} give short time existence of Ricci flow and harmonic map heat flow at $\tau^\ast(p)$ for each $p$. Moreover, we have $C^1$-continuous dependence of Ricci flow and harmonic map heat flow. Also, \zcref[S]{cor:pdeodi} gives strict outward crossing of the box condition $\mathfrak{B}$. Therefore, for each $p \in \mathfrak{P}_\mu$ with $\tau^\ast(p) < \infty$, $Q(p,\tau)$ is continuous on $B_{\epsilon_p}^{\RR^9}(p) \times (\tau^\ast(p) - \epsilon_p, \tau^\ast(p) + \epsilon_p)$ if $\tau^\ast(p) > 0$ and on $(B_{\epsilon_p}^{\RR^9}(p) \cap \mathfrak{P}_\mu) \times [0, \tau^\ast(p) + \epsilon_p)$ if $\tau^\ast(p) = 0$ for sufficiently small $\epsilon_p = \epsilon_p(p) > 0$. Therefore, $\tau^\ast(p)$ is continuous on $\mathfrak{P}_\mu$.

    Define a map 
    \[
        \cE(p) := \dfrac{e^{\gamma \tau^\ast(p)}}{\mu} \cP_+(w_{R(\tau^\ast(p))}h_p(\tau^\ast(p))).
    \]
    By the earlier argument, $\tau^\ast(p)$ is continuous, so this map $\cE \colon \mathfrak{P}_\mu \to S(V_+)$ is also continuous. Then with the map $\overline{B}^{\RR^9}_1(0) \to \mathfrak{P}_\mu$ constructed in the proof of \zcref[S]{prop:pb}, we get a continuous map from $B^9$ to $\partial B^9$ whose restriction to $\partial B^9$ is the identity. This is impossible by the no retraction theorem.
    Hence, there exists $p^\ast \in \mathfrak{P}_\mu$ for which $\tau^\ast(p^\ast) = \infty$.  

    Then for any compact sets $K \subset K' \subset M$, since $R(\tau) \to \infty$, the cutoff $w_{R(\tau)}$ is 1 on $K'$ for large enough $\tau$. Thus $\norm{h_{p^\ast}(\tau)}_{L^2(K')} \le C_{K'}\mu e^{-\gamma \tau}$. Combining this with (B3), we get an interpolation 
    \[
        \norm{h_{p^\ast}(\tau)}_{C^0(K)} \le C_{K,K'} \norm{h_{p^\ast}(\tau)}_{L^2(K')}^{1/2} \norm{h_{p^\ast}(\tau)}_{C^2(K')}^{1/2} \le C_{K,K'} (\mu \eta)^{1/2} e^{-\gamma \tau/2} \to 0
    \]
    as $\tau \to \infty$. Then by parabolic regularity, $h_{p^\ast}(\tau) \to 0$ in $C_{\loc}^{\infty}(M)$. 

    Define 
    \begin{equation} \label{eqn:Ftau}
    \mathcal{F}_\tau := \Phi_{p^\ast, t}^{-1} \circ \iota^{-1} \circ \Psi_{p^\ast, t}^{-1} \colon M_{4R(\tau)} \to N.
    \end{equation}
    By the construction, 
    \[
        \mathcal{F}_\tau^\ast ((T_{p^\ast}-t)^{-1}G(t)) = \overline{g} + h_{p^\ast}(\tau).
    \]
    Therefore, 
    \[
        (N, (T_{p^\ast}-t)^{-1}G(t), x_\tau) \to (M, \overline{g}, o)
    \]
    in the pointed smooth Cheeger--Gromov sense, where $o \in M$ is a point at which $|\Rm_{\overline{g}}(o)| \ne 0$, and $x_\tau := \mathcal{F}_\tau(o)$. 

    Finally, we argue that $T_G = T_{p^\ast}$. The box $\mathfrak{B}$ already gives $T_G \ge T_{p^\ast}$. If $T_G > T_{p^\ast}$, then from the pointed Cheeger--Gromov convergence, 
    \[
        (T_{p^\ast}-t)|\Rm_{G(t)}|_{G(t)}(\mathcal{F}_\tau(o)) \to |\Rm_{\overline{g}}|_{\overline{g}}(o) > 0.
    \]
    But since $|\Rm_{G(t)}|_{G(t)}$ is bounded up to time $T_{p^\ast} < T_G$, the sequence actually goes to zero, giving a contradiction. 
    Therefore, $T_G = T_{p^\ast} < \infty$.
\end{proof}

\section{Strong Rigidity of the blowdown soliton} \label{sec:rigidity}

In this section, we prove the strong rigidity of the blowdown soliton, which states that if a shrinking Ricci soliton is sufficiently close enough to the blowdown soliton on a large but compact set, then the shrinking Ricci soliton must be isometric to the blowdown soliton.

\begin{theorem}[Strong Rigidity of the blowdown shrinking Ricci soliton] \label{thm:srfik}
    Let $(M, \overline{g}, f)$ be the blowdown shrinking Ricci soliton. There exist a large compact set $K \subset M$ and $\epsilon > 0$ with the following property. Suppose that $(M', g', f')$ is a complete shrinking Ricci soliton and $\iota \colon K \to M'$ is a diffeomorphism onto its image such that 
    \[
        \norm{\iota^\ast g' - \overline{g}}_{C^{13}(K)}<\epsilon.
    \]
    Then $(M', g')$ is isometric to $(M, \overline g)$.
\end{theorem}

\begin{proposition} \label{prop:sretb}
    Let $(M, g, f)$ be an $n$-dimensional asymptotically conical shrinking Ricci soliton where $f$ is normalized by $R_{g} + |\nabla^g f|_g^2 - f = 0$. 
    Let $h$ be a smooth symmetric 2-tensor satisfying 
    \[
        (\Delta_f + 2\Rm_{g})h = \lambda h
    \]
    for some $\lambda \ge 0$. If $|h| \le C(1+f)^m$ for some $C>0$ and $m>0$, then $h$ is bounded globally on $M$. Moreover, $|\nabla^k h| \le C_k(1+f)^{-k/2}$ for each $k \ge 0$.
\end{proposition}
\begin{proof}
    Since $M$ is asymptotically conical, there exists $C_{\Rm} > 0$ such that 
    \[
        2|\Rm_g(h)| \le \dfrac{C_{\Rm}}{f} |h|.
    \]
    Define an operator $P = \Delta_f + \left( \dfrac{C_{\Rm}}{f} - \lambda \right)$. Let $u := |h|$. Then by Kato's inequality,
    \[
        u\Delta_f u = \langle \Delta_f h, h \rangle + |\nabla h|^2 - |\nabla u|^2 \ge \langle \Delta_f h, h \rangle.
    \]
    Therefore, 
    \[
        Pu = \left( \Delta_f + \dfrac{C_{\Rm}}{f} - \lambda \right) u \ge 0.
    \]

    Define a barrier $v = 1 - \dfrac{C_{\Rm}+1}{f}$. Then we have 
    \[
        Pv = -\lambda v - \dfrac{1}{f} - \dfrac{(C_{\Rm}+1)(C_{\Rm}+2-n/2)}{f^2} + \dfrac{2(C_{\Rm}+1)R}{f^3}.
    \]
    We may choose $L>0$ large enough so that $v \ge 1/2$ and $Pv < 0$ on $\{f \ge L\}$. Note that 
    \[
        Pf^{m+1} = -(m+1+\lambda)f^{m+1} + [(m+1)(m+n/2)+C_{\Rm}]f^m - (m+1)mRf^{m-1}.
    \]
    Increase $L$ if needed so that $Pf^{m+1} < 0$ on $\{f \ge L\}$. 
    Fix $A>0$ so that $u \le Av$ on $\{f = L\}$. For $\epsilon > 0$, define 
    \[
        \xi_\epsilon := \dfrac{u - Av - \epsilon f^{m+1}}{v}.
    \]
    Since $\xi_\epsilon$ is nonpositive on $\{f = L\}$ and negative when $f$ is sufficiently large, if $\xi_\epsilon$ were to attain any positive value, there is a positive maximum attained on $\{f > L\}$. At the point where the maximum is attained, 
    \[
        P(v\xi_\epsilon) = v\Delta_f \xi_\epsilon + 2\langle \nabla v, \nabla \xi_\epsilon \rangle + \xi_\epsilon Pv < 0.
    \]
    On the other hand, from the previous computations, 
    \[
        P(v\xi_\epsilon) = Pu - APv - \epsilon Pf^{m+1} > 0,
    \]
    giving a contradiction. Therefore, $u - Av - \epsilon f^{m+1} \le 0$ on $\{f \ge L\}$.
    By taking $\epsilon \searrow 0$, we get $u \le Av$. In particular, $|h|$ is bounded.

    The bound on the higher derivatives are obtained by induction: commute the equation and apply the same barrier method to $u_k := (1+f)^{k/2}|\nabla^k h|$.
\end{proof}

\begin{proposition} \label{prop:srspfik}
    There exists $\epsilon > 0$ with the following property. Let $(M, g', f')$ be a shrinking Ricci soliton such that $\norm{g'-g}_{C^{2}(M)} < \epsilon$. Then for the operator $L_{f'} := \Delta_{g', f'} + 2\Rm_{g'}$, 
    \[
        \Spec(L_{f'}|_{\mathcal{T}'}) \subset (-\infty, \dfrac{\lambda_{10}}{2}]
    \]
    where $\mathcal{T}' = \ker \div_{g', f'} \cap (\Ric_{g'})^\perp$ and $\lambda_{10} < 0$ is the 10-th eigenvalue of the blowdown shrinking Ricci soliton.
\end{proposition}

\begin{proof}
    We first note that there is an orthogonal decomposition 
    \[
        L_f^2(S^2T^\ast M) = \RR \Ric_{\overline{g}} \oplus \mathcal{G} \oplus \mathcal{T}
    \]
    where $\mathcal{G}$ is the range of Lie derivatives and $\mathcal{T} := \ker \div_f \cap \Ric_{\overline{g}}^\perp$. 
    
    For the sake of contradiction, we assume that there exist shrinking Ricci solitons $(M, g_j, f_j)$ such that $\norm{g-g_j}_{C^2(M)} < \epsilon_j$ with $\epsilon_j \searrow 0$ and $\Spec(L_{f_j}|_{\mathcal{T}_j'}) \not\subset (-\infty, \lambda_{10}/2]$ where $\mathcal{T}_j' := \ker \div_{g_j, f_j} \cap (\Ric_{g_j})^\perp$. We may assume that $(M, g_j, f_j)$ are normalized so that $R_j + |\nabla^{g_j}f_j|_{g_j}^2 - f_j = 0$.
    
    We first show that $f_j \to f$ in $C_{\loc}^2$. If for some $p \in M$, $|df_j(p)|$ is unbounded, then a locally convergent limit of $\dfrac{df_j}{|df_j(p)|}$ would give a nonzero parallel one-form on the FIK shrinking Ricci soliton, which is a contradiction. Hence, $|df_j(p)|$ is uniformly bounded over $j$ for each $p \in M$, and by the shrinker equation and normalization, $f_j \to f$ in $C_\loc^2$.

    Also, the global $C^2$ convergence gives that outside a fixed compact set $K \subset M$, $\norm{\Ric_{g_j}}_{g_j} \le \frac{1}{4}$, which implies 
    \[
        \dfrac{1}{4}g_j \le \nabla^{g_j, 2} f_j \le \dfrac{3}{4}g_j
    \]
    on $M \setminus K$. Integrating along geodesics gives fixed constants $c_1, c_2, C>0$ such that
    \[
        c_1\rho^2 - C \le f_j \le c_2\rho^2 + C \quad \text{and} \quad |\Rm_{g_j}|_{g_j} \le C
    \]
    on $M$ for every large enough $j$, where $\rho := d_g(o, \cdot)$ for a fixed $o \in M$.

    Now we take $h_j \in \mathcal{T}_j$ such that $L_{f_j}h_j = \lambda_j h_j$, $\norm{h_j}_{L_{f_j}^2} = 1$, and $\lambda_j > \lambda_{10}/2$. Define $u_j := e^{-f_j/2}h_j$, then 
    \[
        \lambda_j e^{f_j/2}u_j = L_{f_j}h_j = e^{f_j/2}\left[ \left(\Delta_{g_j} + 2\Rm_{g_j} + \left( \dfrac{1}{2}\Delta_{g_j}f_j - \dfrac{1}{4}|\nabla^{g_j}f_j|_{g_j}^2 \right)\right)u_j \right].
    \]
    By shrinker identity with $n=4$, $\frac{1}{2}\Delta f_j - \frac{1}{4}|\nabla f_j|^2 = 1 - \frac{1}{4}f_j - \frac{1}{4}R_j$. Hence, pairing with $u_j$, 
    \begin{align*}
        \dfrac{\lambda_{10}}{2} < \lambda_j \norm{u_j}_{L^2}^2 &= \int \langle \left( \Delta_{g_j} + 2\Rm_{g_j} + 1 - \frac{1}{4}f_j - \frac{1}{4}R_j \right)u_j, u_j \rangle d\mu_{g_j}\\
        &\le \int (-|\nabla^{g_j} u_j|^2 - \frac{c_1}{4}\rho^2|u|^2 + (C+1)|u|^2)d\mu_{g_j}.
    \end{align*}
    Since $\norm{u_j}_{L^2} = 1$, we have that both $\norm{\nabla u_j}_{L^2}$ and $\norm{\rho u_j}_{L^2}$ are uniformly bounded over $j$. Moreover, $\lambda_j$ are uniformly bounded from above by $C+1$, so there exists a subsequence (whose indices will also be denoted by $j$) such that $\lim_{j \to \infty} \lambda_j = \lambda_{\infty} \ge \dfrac{\lambda_{10}}{2}$.

    Since $\norm{\rho u_j}_{L^2}$ is uniformly bounded, we have 
    \[
        \sup_{j} \int_{M \setminus B_R(o)} |u_j|^2 d\mu_{g_j} \to 0
    \]
    as $R \to \infty$.
    The Rellich--Kondrachov theorem and this tail bound yield a subsequence converging to some $u_\infty \in L^2(M)$ such that $u_j \to u_\infty$ strongly in $L^2(M)$, $\nabla^{g_j} u_j \rightharpoonup \nabla^{\overline{g}} u_\infty$ weakly in $L^2(M)$, and 
    \[
        \left(\Delta_{\overline{g}} + 2\Rm_{\overline{g}} + \left(1 - \frac{1}{4}f - \frac{1}{4}R\right)\right)u_\infty = \lambda_\infty u_\infty.
    \]
    Then $h_\infty = e^{f/2}u_\infty$ satisfies $L_f h_\infty = \lambda_\infty h_\infty$. But moreover, as $\div_{g_j, f_j} h_j = 0$ and $h_j \perp \Ric_{g_j}$, after taking the limit, $\div_f h_\infty = 0$ and $h_\infty \perp \Ric_{\overline{g}}$. This is a contradiction to \zcref[S]{prop:fiksa}, and the desired result follows.
\end{proof}

\begin{lemma} \label{lem:srextc}
    For any $\epsilon > 0$, there exist $l>0$ and $\delta > 0$ with the following property. Let $(M', g', f')$ be a complete shrinking Ricci soliton with $\iota \colon M_l \to M'$ a diffeomorphism onto its image such that 
    \[
        \norm{\iota^\ast g' - \overline{g}}_{C^4(M_l)} < \delta.
    \]
    Then there is a diffeomorphism $\phi \colon M \to M'$ such that 
    \[
        \norm{\phi^\ast g' - \overline{g}}_{C^2(M)} < \epsilon.
    \]
\end{lemma}
\begin{proof}
    For the sake of contradiction, suppose that for some $\epsilon > 0$, there exist shrinking Ricci solitons $(M_j, g_j, f_j)$, $l_j \nearrow \infty$, $\delta_j \searrow 0$, and $\iota_j \colon M_{l_j} \to M_j$ which are diffeomorphisms onto their image such that $\norm{\iota_j^\ast g_j - \overline{g}}_{C^4(M_{l_j})} < \delta_j$ but no global diffeomorphisms $\phi_j \colon M \to M_j$ with $\norm{\phi_j^\ast g_j - \overline{g}}_{C^2(M)} < \epsilon$ exist. We may assume that $(M_j, g_j, f_j)$ are normalized so that $R_j + |\nabla^{g_j} f_j|_{g_j}^2 - f_j = 0$.

    Let $\hat{g}_j = \iota_j^\ast g_j$ and $\hat{f}_j = \iota_j^\ast f_j$. First, note that if $|df_j(p)|$ is unbounded, then by the earlier proof of \zcref[S]{prop:srspfik}, we get a contradiction. Hence $|df_j|$ is bounded over $j$ at each point, so $\Hess_{\hat{f}_j} = \frac{1}{2}\hat{g}_j - \Ric_{\hat{g}_j}$, $\nabla^{\overline{g}, 2}(\hat{f}_j - f)$ is small on $M_{l_j}$. Since $(M, \overline{g}, f)$ does not have nonzero parallel one-forms, $\nabla^{\overline{g}}(\hat{f}_j - f) \to 0$ locally as $j \to \infty$. Then the normalization and the $C^4$-closeness of the metrics show that $\hat{f}_j \to f$ in $C_{\loc}^4(M)$ as $j \to \infty$.

    Let $\epsilon_0$ and $\delta_0$ be the dimensional constants from \cite[Theorem 24]{LW20}. 
    
    Fix a minimum point $o$ of $f$, and set $o_j = \iota_j(o)$. Choose $a_0$ sufficiently large that 
    \[
        \nu (B_{2\epsilon_0^{-1}}^{\overline{g}}(x), \overline{g}, \epsilon_0^{-2}) > -\delta_0/2 \quad \text{whenever } f(x) \ge a_0^2/4,
    \]
    and set $a = 3a_0$. Since $f_j(o_j) \to f(o)$, the estimate by \cite{CZ10} implies that, for all sufficiently large $j$, 
    \[
        \{x : f_j(x) \le a^2/4 \} \subset B_{g_j}(o_j, a+C),
    \]
    where $C$ is a dimensional constant independent of $j$. By metric convergence, we may choose a fixed $L$ sufficiently large that 
    \[
        B_{g_j}(o_j, a+C+2\epsilon_0^{-1}) \subset \iota_j(M_L)
    \]
    for all sufficiently large $j$.

    Consequently, the entire level $\Sigma_j := f_j^{-1}(a^2/4)$ lies in this fixed region. The local $C^4$ convergence $\iota_j^\ast f_j \to f$ implies that $\iota_j^{-1}(\Sigma_j)$ is a normal graph over $\Sigma := f^{-1}(a^2/4)$. In particular, $\Sigma_j$ is connected. Moreover, local metric convergence and the choice of $a_0$ give 
    \[
        \nu(B_{\epsilon_0^{-1}}^{g_j}(y), g_j, \epsilon_0^{-2}) > -\delta_0 \quad \text{for every } y \in \Sigma_j.
    \]
    Applying \cite[Corollary 9]{LW20}, we conclude that $f_j \ge a^2/4$ is diffeomorphic to $\Sigma_j \times [0, \infty)$ and has a unique asymptotically conical end, diffeomorphic to the end of $(M, \overline{g})$.

    The normal graph parametrization gives a diffeomorphism $\theta_j \colon \Sigma \to \Sigma_j$ such that 
    \begin{equation} \label{eqn:srextc1}
        \norm{\theta_j^\ast g_j|_{\Sigma_j} - \overline{g}|_\Sigma}_{C^2(\Sigma)} \to 0
    \end{equation}
    as $j \to \infty$.

     We extend $\theta_j$ to the exterior region by defining 
    \[
        \phi_j(\Upsilon_\sigma x) = \Upsilon_\sigma^{(j)}\theta_j(x)
    \]
    for $x \in \Sigma$. Then we glue this extended map with $\iota$ on a collar region $\{a-1 \le \overline{b} \le a+1\} \subset M$.

    The $C^2$ norm on the inner region $M_{l_j}$ is already controlled by the assumption and \eqref{eqn:srextc1} because for any $b>0$, 
    \[
        \norm{\phi_j^\ast g_j - \overline{g}}_{C^2(\{\overline{b} \le b\})} \to 0.
    \]
    To obtain the global $C^2$ estimate, we define the followings:
    \[
        \widetilde{g}_\sigma^{(j)} := e^{-\sigma}\phi_j^\ast(\Upsilon_\sigma^{(j)})^\ast  g_j, \quad \text{and} \quad \widetilde{g}_\sigma := e^{-\sigma}\Upsilon_\sigma^\ast \overline{g}
    \]
    which are the pullback of the rescaled Ricci flow on $M$. Since \cite[Corollary 8]{LW20} provides a uniform curvature bound, by Shi's estimate, we get uniform bounds on the derivatives of curvature as well for large enough $j$ in the exterior region. More precisely, for small enough $s>0$, define a collar region $U := \{x \in M : d_{\overline{g}}(x, \Sigma) < s\}$. Then by integrating the rescaled Ricci flow equation $\partial_\sigma \widetilde{g}_\sigma^{(j)} = -2e^{-\sigma}\Ric_{\widetilde{g}_\sigma^{(j)}}$,
    \[
        \norm{\widetilde{g}_\sigma - \widetilde{g}_{\sigma_0}}_{C^2(U; \overline{g})} \le Ce^{-\sigma_0} \quad \text{and} \quad \norm{\widetilde{g}_\sigma^{(j)} - \widetilde{g}_{\sigma_0}^{(j)}}_{C^2(U; \overline{g})} \le Ce^{-\sigma_0}.
    \]
    Then by taking $\sigma_0$ large enough so that $Ce^{-\sigma_0} < \epsilon/2$ and take $j$ large enough so that $l_j \ge 2ae^{\sigma_0/2}$, we get 
    \[
        \norm{\phi_j^\ast g_j - \overline{g}}_{C^2(M)} < \epsilon,
    \]
    which is a contradiction. Therefore, the desired lemma follows.
\end{proof}

\begin{proof}[Proof of {\zcref[S]{thm:srfik}}]
    By taking $K= M_l$ with $l$ sufficiently large and $\epsilon > 0$ sufficiently small, \zcref[S]{prop:sretb}, \zcref[S]{prop:srspfik}, \zcref[S]{lem:srextc}, and \zcref[S]{rem:ssct} give that $(M', g', f')$ also satisfies an analogous statement of \zcref[S]{thm:A}.

    Construct a closed Riemannian manifold $(N, G)$ on which large enough compact set the metric is identical to $g'$ and outside of this compact set the metric is sufficiently flat so that $(N, G)$ satisfies both assumptions of \zcref[S]{thm:A} for $(M, \overline{g}, f)$ and \zcref[S]{thm:A} for $(M', g', f')$. Then as a consequence, Ricci flow starting from $(N, G)$ develops a finite-time singularity whose tangent flow is unique and modeled on both $(M, \overline{g}, f)$ and $(M', g', f')$. This is only possible when $(M, \overline{g}, f)$ and $(M', g', f')$ are isometric.
\end{proof}

\section{Proofs of the other results}\label{sec:consequences}

We now deduce the other results of the introduction from \zcref[S]{thm:mainthm} and its proof in \zcref[S]{sec:proofmain}. By parabolic rescaling, we assume that $T_0=1$, and we let $p^\ast$ and $\mathcal F_\tau$ be as in the proof of \zcref[S]{thm:mainthm}, so that $T_G=T_{p^\ast}$ and \eqref{eqn:Ftau} holds.

\begin{proof}[Proof of {\zcref[S]{thm:B}}]
For $t<T_G$, we write $\mathcal F_t$ and $h(t)$ for $\mathcal F_\tau$ and $h_{p^\ast}(\tau)$, where $\tau=-\log\frac{T_G-t}{T_G}$. In terms of $t$, \eqref{eqn:Ftau}, \ref{box:B3}, and the decay of $h_{p^\ast}$ obtained from \ref{box:B1}-\ref{box:B2} in the proof of \zcref[S]{thm:mainthm} rewrite
\begin{equation}\label{eq:controls}
\mathcal F_t^*G(t)=(T_G-t)\big(\overline g+h(t)\big),\quad \norm{h(t)}_{C^{12}(M_{2R(t)})}\le\eta,\quad \text{ and }\quad\norm{h(t)}_{C^0(K)}\le C_K\Big(\frac{T_G-t}{T_G}\Big)^{\frac\gamma2},
\end{equation}
on $M_{4R(t)}$, where $R(t)^2=R_0^2+\frac{8\gamma}{\theta^2}\log\frac{T_G}{T_G-t}$, for every compact set $K\subset M$. Moreover, $h(t)\to0$ in $C^\infty_{\loc}(M)$ as $t\to T_G$.

As in \cite[Lemma 9.1]{CS21}, pseudolocality on annuli propagates these controls up to the singular time. For $r>0$ and $A>0$, let $t_r:=T_G-(r/A)^2$ and $\mathcal A_r:=\mathcal F_{t_r}(M_{3A}\setminus M_A)$. By \eqref{eq:controls}, $r^{-2}G(t_r)$ on $\mathcal A_r$ is close to $A^{-2}\overline g$ on $M_{3A}\setminus M_A$ for $r$ small, hence to an annulus of the cone $\mathcal C$ for $A$ large. Since $\mathcal C$ is smooth away from its tip, the balls of radius $\kappa r$ centered in $\mathcal A_r$ satisfy the hypotheses of \zcref[S]{prop:pseudoloc} at time $t_r$, for some fixed $\kappa>0$. Taking $A$ large depending on $\kappa$, \zcref[S]{prop:pseudoloc} applies on $[t_r,T_G)$, and gives
\begin{equation}\label{eq:curv}
|\Rm_{G(t)}|\le Cr^{-2}\quad\text{on }\mathcal A_r,\ \text{for }t\in[t_r,T_G).
\end{equation}

\textbf{Step 1: uniqueness of the tangent flow.} By \eqref{eq:controls}, $(N,(T_G-t)^{-1}G(t),\mathcal F_t(o))$ converges to $(M,\overline g,o)$ as $t\to T_G$ for every $o\in M$, without taking subsequences. By \eqref{eq:controls} and \eqref{eq:curv}, for $A$ large and $t$ close to $T_G$, the maximum of $(T_G-t)|\Rm_{G(t)}|$ is attained in $\mathcal F_t(M_A)$.

\textbf{Step 2: conical singular time slice.} We follow the strategy of \cite[Corollary 1.2]{CS21} for mean curvature flow, where the uniqueness of the tangent flow is combined with curvature bounds on annuli such as \eqref{eq:curv}.

\emph{Step 2a: the singular time slice.}
Write $A(t)=(\Phi_{p^\ast,t}^{-1})^\ast G(t)$,
and $\Sigma^*=\iota^{-1}\Xi_{p^\ast}^{-1}(\overline\Sigma)$. The metrics $\beta(t)$ converge smoothly on $N\setminus \Sigma^*$ to $\beta_{T_G}$, which is exactly conical near $\Sigma^*$. Writing $\rho$ for its cone radius, set the annular norm
 $\|q\|_{C_r^4}
 :=\sum_{j=0}^4r^j
 \sup_{\{r<\rho<2r\}}|\nabla_{\beta_{T_G}}^jq|_{\beta_{T_G}}.$
By \ref{box:B4}, \eqref{eqn:rfce1}, and \cite[Proposition C.3]{Sto26}, applied on slightly larger parabolically rescaled annuli, we get the control
\[
 \|\partial_tA(t)\|_{C_r^4}
 +\|\partial_t\beta(t)\|_{C_r^4}\leqslant Cr^{-2},
 \qquad {T_G}-r^2/L_0^2\leqslant t<{T_G},
\]
with $C,L_0$ independent of $r$. Thus $A(t)$ converges smoothly to $A_{T_G}$ on $N\setminus \Sigma^*$. For fixed $L\geqslant L_0$ and $t_r={T_G}-r^2/L^2$, \eqref{eq:controls} gives $\|A(t_r)-\beta(t_r)\|_{C_r^4}\to 0$ as $r\to 0$. Integration and triangle inequality yields, as $r\to 0$,
\[
 \|A_{T_G}-\beta_{T_G}\|_{C_r^4}
 \leqslant CL^{-2}+o(1).
\]
Letting $r\to0$ and then $L\to\infty$ proves
\cite[Definition 1.1]{GS18}. Finally, \ref{box:B4} gives vanishing diameters of the inner regions, so the metric completion of $G(T_G)$ adds one point and is compact. Since $A(t)$ and $G(t)$ are isometric, it is also their Gromov-Hausdorff limit.

\emph{Step 2b: continuous dependence.} For $t$ close to $T_G$, \eqref{eq:controls} and \eqref{eq:curv} show that $G(t)$ satisfies \eqref{eq:A1}-\eqref{eq:A2} with $T_0=T_G-t$, $U=\mathcal F_t(M_{100\Gamma})$ and $\iota = \mathcal{F}_t^{-1}$: at distance $d \geqslant \Gamma\sqrt{T_G-t}$ from $x_0$, we have $|\Rm_{G(t)}|\leqslant Cd^{-2}<(\Gamma^2 T_0)^{-1}$, and the balls of radius $\Gamma\sqrt{T_0} \ll d$ are almost Euclidean. These conditions are open, and the Ricci flow depends continuously on its initial data, from time $0$ to the fixed time $t<T_G$. So \zcref[S]{thm:A} applies at time $t$ to the flows of nearby initial data, and their singular times $T'$ satisfy $|T'-t-T_0|<\frac{T_0}{10}$, as does $T_G$. Since $T_0$ can be taken arbitrarily small, $T_G$ depends continuously on $G_0$. The constants in \eqref{eq:controls} and \eqref{eq:curv} are uniform for nearby initial data, so the conical approximation of $G(T_G)$ at small scales is uniform. Together with the continuous dependence at the times $t_r$ for the larger scales, $G(T_G)$ depends continuously on $G_0$.

\textbf{Step 3: blowdown surgery.} By Step 2a, the metric completion of $G(T_G)$ satisfies the hypotheses of \cite[Theorem 1.1]{GS18}, with the cone $\cC$ over a sphere with curvature operator at least $1$. This theorem directly gives the Ricci flows coming out of $G(T_G)$ with the first three properties, obtained by gluing in the expanding soliton with the same asymptotic cone as \cite{Cao97,FIK03}, which lives on $\mathbb C^2\approx\mathbb R^4$. In $N$, this replaces a neighborhood of the sphere $\mathcal F_t(\overline\Sigma)$, which is isotopic to $\iota^{-1}(\overline\Sigma)$, by a ball, so these flows live on $N'$. The constants of \cite[Theorem 1.1]{GS18} only depend on the conical approximation of the singular time slice and on its curvature away from $x_0$, which are uniform for converging singular time slices by Step 2b. Hamilton's compactness theorem \cite{Ham95} then gives the last property.
\end{proof}

\begin{remark}
    The main theorem of \cite{GS18} is applicable since the estimates of \cite[Lemma 3.1 and Theorem 6.1]{DL17} used in \cite[Lemma 4.2]{GS18} remain valid for this model. Indeed, the tensor maximum principle in \cite[Theorem 4.1]{DL17} requires only nonnegative sectional curvature, and the scalar comparison in their Theorem 4.5 uses only $Rg-2\operatorname{Ric}\geqslant0$.
\end{remark}

\begin{proof}[Proof of Corollary \ref{cor:open}]
Let $g$ be a metric whose Ricci flow $G$ develops a blowdown singularity at its first singular time $T$. From \zcref[S]{thm:B}, for some $t<T$, $G(t)$ satisfies the assumptions \eqref{eq:A1}-\eqref{eq:A2} of \zcref[S]{thm:A} with $T_0=T-t$. By continuous dependence of the Ricci flow on its initial data, from time $0$ to the \textit{fixed} time $t<T$, the flows from metrics $C^{14}$-close to $g$ exist on $[0,t]$, and are $C^{13}$-close to $G(t)$ at time $t$. Since \eqref{eq:A1}-\eqref{eq:A2} are open conditions, \zcref[S]{thm:A, thm:B}(1)-(2) apply to them at time $t$. They develop a blowdown singularity at their first singular time $T'$, with $|T'-t-T_0|<\frac{T_0}{10}$. Since $|T-t-T_0|<\frac{T_0}{10}$ as well, and $T_0$ can be taken arbitrarily small, $T'$ depends continuously on the metric.
\end{proof}

\begin{proof}[Proof of {\zcref[S]{cor:entropy}}]
Since the singularity is of Type I, \cite{MM15} gives 
\[
    \lim_{t \nearrow T_G} \cW(G(t), f_{p, T_G}(t), T_G-t) = \mu(\overline{g}, 1).
\]
On the other hand, $\cW(G(t), f_{p, T_G}(t), T_G-t)$ is monotone increasing and is constant only if the flow is self-similar. However, as the flow is compact non developing the blowdown singularity, which is noncompact, the flow cannot be self-similar. Together with \cite{CN09}, we get a strict inequality
\[
    \mu(G_0, T_G) \leqslant \cW(G_0, f_{p, T_G}(0), T_G) < \lim_{t \nearrow T_G} \cW(G(t), f_{p, T_G}(t), T_G-t) = \mu(\overline{g}, 1). \qedhere
\]
\end{proof}

\begin{proof}[Proof of {\zcref[S]{thm:C}}]
Let $I\geqslant1$ be the number of summands and order the spheres $\Sigma_1,\dots,\Sigma_I$ of self-intersection $\pm1$ as prescribed. Fix a metric on $X$ flat on $I$ disjoint unit balls $D_1,\dots,D_I$ at mutual distance at least $1$, and obtain $g_\delta$ by gluing by cut-off $s_i\overline g$ on $M_{1/(2\sqrt{s_i})}$, $s_i:=\delta^{2(I+1-i)}$, into $D_i$, reversing orientation for $\mathbb{CP}^2$ summands, and interpolating with the flat metric for radii in $[1/2,1]$. Let $G$ be its Ricci flow, continued as in \zcref[S]{thm:B}(3) at its singular times $T_1<T_2<\cdots$, with singular points $x_j$. We use the notation $C$ below for constants that depend only on $(M,\overline g)$, the metric on $X$, the interpolation and $I$. They may change from line to line.

For fixed $1\leqslant i\leqslant I$ and $m\geqslant4$, \zcref[S]{prop:pseudoloc}, \cite{Ham95} and \cite[Theorem 1.1]{CZ06} give, 
\[
s_i^{-1}G(s_it)|_{D_i}\longrightarrow(1-t)\Psi_{0,t}^*\overline g\quad\text{in }C^\infty_{\loc}\bigl(M\times[0,1-m^{-1}]\bigr)
\]
uniformly over the continuations of the flow at $T_1,\dots,T_{i-1}$, as $\delta\to0$. In the cone chart of Step 2a of the proof of \zcref[S]{thm:B}, in which $\beta_{T_j}$ is exactly conical with cone radius $\rho$, \cite[Propositions C.3 and C.4]{Sto26} and pseudolocality on annuli give the condition \cite[(1-1)]{GS18} for $\rho\in(0,1/4]$ at each $x_j$, with $k_Z\leqslant C(\eta+\Gamma^{-2})+o_\delta(1)$ in their notation. We will take $C(\eta+\Gamma^{-2})<\kappa/2$ in \zcref[S]{app:consdep}, with $\kappa$ fixed below.

There is $C$ such that, for $j\leqslant k<I$, $t\in(T_k,T_{k+1})$ and $t-T_k\leqslant C^{-1}$,
\begin{equation}\label{eqn:GS}
|\Rm_{G(t)}|\leqslant\frac{C}{t-T_k}\qquad\text{and}\qquad\Vol_{G(t)}B_r>\frac12\omega_4r^4
\end{equation}
on $\{r_M\leqslant C^{-1}\}$, with $r_M$ as in \cite[Section 5.5]{GS18} for the flow coming out of $G(T_j)$ and $B_r$ centered there with $r\leqslant C^{-1}\sqrt{t-T_k}$. For $j=k$, this is \cite[Theorem 3.1]{GS18} for the approximating flows of \cite[Section 5.3]{GS18}, rescaled so that $r_0=1$ in \cite[Definition 1.1]{GS18}, with \cite[(3-15)]{GS18} for the volume. Its constants depend only on the cone $\mathcal C$, and it holds as long as the flows exist, up to a time depending only on $\mathcal C$, so it passes by \cite[(5-6)]{GS18} and continuous dependence on initial data to $G$ on $(T_k,T_{k+1})$. 

For $j<k$, we apply \cite[Theorem 3.1]{GS18} again around $x_j$ to the approximating flows of the case $j=k$. On $D_j$, the approximating metrics of \cite[Section 5.1]{GS18} at $x_k$ coincide with the rescaled $G(T_k)$. By \cite[Theorem 3.1(2)]{GS18} and pseudolocality on annuli at $x_j$ on $(T_{k-1},T_k)$, which give $C^4$ closeness to the expander at scale $t-T_j$ and to $\mathcal C$, respectively, these metrics, rescaled so that $\{r_M\leqslant C^{-1}\}$ becomes $\{r_M\leqslant1\}$, lie in the class of \cite[Definition 3.1]{GS18} with expander scale $s$ equal to the rescaled $T_k-T_j$.  In the notation of \cite[Theorem 3.1(2)]{GS18}, the application at $T_k$ needs the $\varepsilon$ of the application at $T_{k-1}$ to be at most its $\eta_1$, so we choose these $\varepsilon$ iteratively from $T_{I-1}$ back to $T_j$. This fixes $\kappa=\kappa(\mathcal C,I)$ by \cite[Section 5.2]{GS18}.

Fix $m>4+2(1+C)^2c^2\Gamma^2$ with $C$ as in \eqref{eqn:GS}, set $t_i:=(1-m^{-1})s_i$, and let $\delta$ be sufficiently small. We apply \zcref[S]{thm:A} inductively at $t_i$ with $T_0=s_i/m$ and $\iota=\Psi_{0,1-m^{-1}}$. The convergence of $s_i^{-1}G(s_it)|_{D_i}$ at $t=1-m^{-1}$ gives \eqref{eq:A1}, and{\color{blue}, with \zcref[S]{prop:pseudoloc},} \eqref{eq:A2} on $D_i$. Outside $D_i$, \eqref{eqn:GS}, the convergence of $s_j^{-1}G(s_jt)|_{D_j}$ for $j>i$, \zcref[S]{prop:pseudoloc} and, for $i>1$, $t_i-T_{i-1}\geqslant s_i/2$ give
\[
|\Rm_{G(t_i)}|\leqslant
\begin{cases}
2Cs_i^{-1}&\text{on }D_j\text{ for }j<i,\\
C\delta^2s_i^{-1}&\text{on }D_j\text{ for }j>i,\\
C&\text{elsewhere},
\end{cases}
\]
so the choice of $m$, the volume bound in \eqref{eqn:GS} and the same references give \eqref{eq:A2} at radius $c\Gamma\sqrt{s_i/m}$. The same estimates, at every $t$ in place of $t_i$, exclude singular times in $(T_{i-1},t_i]$ for $i>1$, and in $[0,t_1]$ for $i=1$. Hence, by \zcref[S]{thm:A, thm:B}, $G$ blows down a sphere isotopic to $\Sigma_i$ at $T_i$, with $|T_i-s_i|<s_i/(10m)$, and after $T_I$ it lives on $X$. These estimates hold for the Ricci flow of every $g$ with $\|g-g_\delta\|_{C^{14}(g_\delta)}<\delta$, a $C^{14}$-open set.
\end{proof}

\appendix 
\section{Constant Dependence} \label{app:consdep}
In this appendix, we summarize the order in which the constants are chosen and clarify the dependence of each constant.

\subsection{Fix constants $c$, $\gamma$, $\theta$, $\theta'$, and $\theta_0$}
With $\delta_{\ps} = \frac{1}{2}$, we take $\epsilon_{\ps}$ from \zcref[S]{prop:pseudoloc}. We then fix $c>0$ so that for every $l>0$, 
\[
    |\Rm_{\overline{g}}|_{\overline{g}} \leqslant \dfrac{1}{(2cl)^2} \quad \text{on } M \setminus M_l.
\]
We then fix $\gamma$ so that it satisfies \eqref{eqn:gamma}. Then $\theta$, $\theta'$, and $\delta_\theta$ are chosen by \eqref{eqn:theta-exp}.
A constant $\theta_0$ is then chosen sufficiently small so that $C_{\gamma \theta_0}$ defined \zcref[S]{prop:pdeodi} is negative. These are the first constants that are chosen, so we do not clarify whether other constants are depending on $c$, $\gamma$, $\theta$, $\theta'$, $\delta_\theta$, or $\theta_0$.

\subsection{Choose $K$}
We first fix $r>0$ which is a local radius that only depends on $(M, \overline{g})$. We take $r>0$ so that it satisfies 
\[
    0 < r < \min \{ \dfrac{1}{2}, \dfrac{1}{4}e^{-4} \inf_M \inj_{\overline{g}} \}.
\]
Set $K''' = 2$ and choose $K''$ large enough so that \zcref[S]{prop:c12-im} and parabolic regularity in the proof of \zcref[S]{prop:c0-im-init} hold. Then choose $K'$ large enough for \zcref[S]{prop:c0-im}. Finally, take $K$ large enough so that $K \geqslant 1$, \zcref[S]{prop:c0-in, prop:c0-im-init} hold, and 
\[
    \dfrac{4C_h^2 e^{8C_w}}{K^2} \leqslant \dfrac{1}{2(K'')^2} \quad \text{and} \quad \dfrac{2C_h}{K} \leqslant \dfrac{1}{K''}.
\]

\subsection{Choose $R_0$}
Next, we take $R_0 = R_0(K, K'') \gg 1$ large enough so that \zcref[S]{prop:aprhmhf}, \zcref[S]{prop:bar}, \zcref[S]{cor:c12-im}, \zcref[S]{prop:pdeodi}, \zcref[S]{cor:pdeodi}, \zcref[S]{prop:pb}, \zcref[S]{prop:pbs}, and \zcref[S]{prop:c0-im-init} hold. In particular, we take $R_0$ so that 
\[
    \dfrac{\nu}{2K} e^{(1-\theta^2)R_0^2/8} \geqslant 2C_ER_0^2
\]
where $\nu := \min \{\gamma, -C_{\gamma \theta_0} \} > 0$ as defined in the proof of \zcref[S]{cor:pdeodi} and $C_E$ is taken from \zcref[S]{lem:ee}.

\subsection{Choose $p_0$}
We next take $0 < p_0 = p_0(R_0) \ll 1$ small enough so that \zcref[S]{prop:pb}, \zcref[S]{prop:rfce}, and \zcref[S]{prop:hmhfce} hold. We then define 
\[
    \overline{\mu} = \dfrac{p_0}{4\norm{J_{R_0}^{-1}}}.
\]

\subsection{Choose $\eta$, then define $\mu$}
We take $\eta$ small enough so that \zcref[S]{prop:aprhmhf}, \zcref[S]{prop:pi}, \zcref[S]{prop:bar}, \zcref[S]{cor:c12-im}, \zcref[S]{prop:pdeodi}, \zcref[S]{cor:pdeodi}, \zcref[S]{prop:c0-im-init}, \zcref[S]{prop:rfce}, and \zcref[S]{prop:hmhfce} hold.
In particular, $\eta$ is small enough so that $16\eta < \overline{\eta}_4$ where $\overline{\eta}_4$ is taken from \cite[Proposition A.9]{BK22} with $n=4$ for harmonic map heat flow continuation and $\mu$ defined as below satisfies the following inequality.
\[
    \mu := \dfrac{\eta}{K}e^{-\theta^2R_0^2/8} < \min \left\{\overline{\mu}, \dfrac{1}{24C_s} \right\}
\]
where $C_s$ is taken from \zcref[S]{prop:pbs}.

\subsection{Choose $\Gamma$}
We choose $\Gamma = \Gamma(R_0, \eta) \gg 1$ large enough so that \zcref[S]{prop:aprhmhf}, \zcref[S]{prop:bar}, \zcref[S]{cor:c12-im}, \zcref[S]{prop:pb}, \zcref[S]{prop:pbs}, and \zcref[S]{prop:c0-im-init} hold. We also require $(\epsilon_{\ps}c\Gamma)^2 > 2$ so that on the exterior region singularities cannot happen for time at most 2.

\subsection{Choose $\epsilon$}
Finally, we take $0 < \epsilon = \epsilon(R_0, \eta, \mu, \Gamma) \ll 1$ small enough so that $\epsilon \leqslant \mu$ and \zcref[S]{prop:aprhmhf}, \zcref[S]{prop:bar}, \zcref[S]{cor:c12-im}, \zcref[S]{prop:pb}, \zcref[S]{prop:pbs}, \zcref[S]{prop:c0-im-init}, \zcref[S]{prop:rfce}, and \zcref[S]{prop:hmhfce} hold.

\section{Pseudolocality of Ricci flow}

The following is a pseudolocality theorem for the Ricci flow, first established by Perelman \cite[Theorem 10.1]{Per02}. We state here a slightly more general version due to Lu \cite{Lu10}.
Here, $\omega_n$ denotes the volume of the unit ball in $\mathbb{R}^n$.
\begin{proposition}[Perelman's Pseudolocality] \label{prop:pseudoloc}
    Given $n \geqslant 2$ and $0 < \delta_{\ps} < 1$, there exists $\epsilon_{\ps} = \epsilon_{\ps}(n, \delta_{\ps}) > 0$ with the following property.
    Let $(M^n, g(t))_{t \in [0,T)}$ be a Ricci flow on a closed Riemannian manifold $M$. Suppose that for $x_0 \in M$ and $r_0 > 0$, 
    \begin{align}
        |\Rm_{g(0)}|_{g(0)}(x) \leqslant \dfrac{1}{r_0^2} \quad \text{for } x \in B_{r_0}^{g(0)}(x_0), \label{eqn:pseudoloc1}\\
        \Vol_{g(0)} B_{r_0}^{g(0)}(x_0) \ge \omega_n(1-\delta_{\ps})r_0^n, \label{eqn:pseudoloc2}
    \end{align}
    then 
    \begin{equation}
        |\Rm_{g(t)}|_{g(t)}(x) \le \dfrac{1}{\epsilon_{\ps}^2r_0^2} \quad \text{for } x\in B_{\epsilon_{\ps}r_0}^{g(t)}(x_0) \label{eqn:pseudoloc3}
    \end{equation}
    for $t \in [0, \min \{T, (\epsilon_{\ps}r_0)^2\})$.
    
\end{proposition}
\bibliographystyle{amsalpha}
\bibliography{refs}

\end{document}